\documentclass[reqno,11pt]{amsart}
\usepackage{amsmath,amsthm, amscd, amssymb, amsfonts}
\usepackage[all]{xy}
\usepackage{hhline}
\usepackage[dvips, dvipsnames, usenames]{color}

\numberwithin{equation}{section}

\newcommand{\ord}{\operatorname{ord}}

\newcommand{\oct}{\mathfrak O}

\newcommand\oc{{\mathcal O }}
\newcommand\ocs{{\mathcal O_{s} }}
\newcommand\toba{{\mathfrak B }}

\newcommand{\trid}{\triangleright}

\newcommand{\ku}{\mathbb C}

\newcommand{\Z}{{\mathbb Z}}

\newcommand{\G}{{\mathbb G}}

\newcommand{\D}{{\mathcal D}}

\newcommand{\Oc}{{\mathcal O}}

\newcommand\card{\operatorname{card}}

\newcommand\sgn{\operatorname{sgn}}

\theoremstyle{plain}

\newtheorem{maintheorem}{Theorem}
\newtheorem{lema}{Lemma}[section]
\newtheorem{theorem}[lema]{Theorem}

\theoremstyle{definition}

\theoremstyle{remark}
\newtheorem{obs}[lema]{Remark}

\newcommand\id{\operatorname{id}}

\newcommand\Id{\operatorname{Id}}

\newcommand\sn{\mathbb S_n}

\newcommand\an{\mathbb A_n}

\newcommand\s{\mathbb S}

\def\pf{\begin{proof}}
\def\epf{\end{proof}}

\theoremstyle{remark}

\begin{document}

\renewcommand{\baselinestretch}{1.2}

\thispagestyle{empty}

\title[pointed Hopf algebras: Mathieu simple groups]{On pointed Hopf algebras associated \\ with the Mathieu simple groups}

\author[Fernando Fantino]{Fernando Fantino}
\address{\noindent Facultad de Matem\'{a}tica, Astronom\'{i}a y F\'{i}sica\\
Universidad Nacional de C\'{o}rdoba \\ CIEM - CONICET,
(5000) Ciudad Universitaria \\
C\'{o}rdoba \\Argentina} \email{fantino@famaf.unc.edu.ar}

\thanks{This work was partially
supported by CONICET, ANPCyT and Secyt (UNC)}

\date{\today}

\begin{abstract}
Let $G$ be a Mathieu simple group, $s \in G$, $\oc_s$ the
conjugacy class of $s$ and $\rho$ an irreducible representation of
the centralizer of $s$. We prove that either the Nichols algebra
$\toba(\oc_s,\rho)$ is infinite-dimensional or the braiding of the
Yetter-Drinfeld module $M(\oc_s, \rho)$ is negative. We also show
that if $G=M_{22}$ or $M_{24}$, then the group algebra of $G$ is
the only (up to isomorphisms) finite-dimensional complex pointed
Hopf algebra with group-likes isomorphic to $G$.
\end{abstract}

\maketitle

\section*{Introduction}\label {0}

This article contributes to the classification of
finite-dimensional complex pointed Hopf algebras $H$ whose group
of group-like elements $G(H)$ is isomorphic to a Mathieu simple
group: $M_{11}$, $M_{12}$, $M_{22}$, $M_{23}$ or $M_{24}$.

The crucial step, in order to classify finite-dimensional complex
pointed Hopf algebras with a fixed $G(H)=G$, is to determine when
a Nichols algebra of a Yetter-Drinfeld module over $G$ is
finite-dimensional -- see \cite{AS3}.

The irreducible Yetter-Drinfeld modules over $G$ are
determined by a conjugacy class $\oc$ of $G$ and an irreducible
representation $\rho$ of the centralizer $G^{s}$ of a fixed $s \in \oc$.
Let $M(\oc,\rho)$ be the corresponding Yetter-Drinfeld module and
let $\toba(\oc,\rho)$ denote its Nichols algebra.

The classification of finite-dimension Nichols algebras over an
abelian group $G$ follows from \cite{AS2,H1,H2}; this leads to
substantial classification results of pointed Hopf algebras $H$
with abelian $G(H)$ -- see \cite{AS4}. The next problem is to
discard irreducible Yetter-Drinfeld modules over a finite
non-abelian group containing a braided vector subspace with
infinite-dimensional Nichols algebra. It is natural to begin by
simple or almost simple groups; see \cite{AZ,AF1,AF2,AFZ}, for
$\an$ or $\sn$; \cite{FGV} for $\mathbf{GL}(2,\mathbb{F}_q)$ or
${\mathbf{SL}(2,\mathbb{F}_q)}$; and \cite{FV} for
$\mathbf{PGL}(2,\mathbb{F}_q)$ or
${\mathbf{PSL}(2,\mathbb{F}_q)}$. We plan to consider the other
sporadic groups in \cite{AFGV}.

Let us say that $M(\oc,\rho)$ has \emph{negative braiding} if the
Nichols algebra of any braided subspace corresponding to an
abelian subrack is (twist-equivalent to) an exterior algebra.

We summarize the investigation in this paper in the next
statement.

\begin{maintheorem}
Let $G$ be a Mathieu simple group, $s \in G$, $\ocs$ the conjugacy
class of $s$ and $\rho\in \widehat{G^s}$. If
$\dim\toba(\oc_{s},\rho)<\infty$, then $(\oc_{s},\rho)$ is one of
the pairs listed in Table \ref{maintabla}. In particular, any
finite-dimensional complex pointed Hopf algebra $H$ with
$G(H)\simeq M_{22}$ or $M_{24}$ is necessarily isomorphic to the
group algebras $\ku M_{22}$ or $\ku M_{24}$, respectively.
\end{maintheorem}

\begin{table}[t]
\begin{center}
\begin{tabular}{|c|c|c|c|c|c|}
\hline $G$ & $j$ &  $|s_j|$ & {\bf Centralizer} & {\bf
Representation} & $\dim M(\oc_{s_j},\rho)$
\\ \hline \hline $M_{11}$ & $4$& 4 & $\langle x\rangle\simeq \Z_8$, $x^6=s_4$ &
$\nu_2(x):=i$ & 990
\\ \cline{5-5}& & &  & $\nu_6(x):=-i$ &
\\  \cline{2-6}  & $6$ & 8 & $\langle s_6 \rangle\simeq \Z_8$ &
$\chi_{(-1)}$ & 990
\\ \cline{2-6}&  $7$ & 8 & $\langle s_7 \rangle \simeq \Z_8$ &
$\chi_{(-1)}$ & 990
\\ \hline \hline $M_{12}$
& $13$ & 10 & $\langle s_{13} \rangle \simeq \Z_{10}$ &
$\chi_{(-1)}$ & 9504
\\ \hline
\hline $M_{23}$  &$12$& 14 & $\langle s_{12}\rangle\simeq \Z_{14}$
& $\chi_{(-1)}$ & 728640
\\  \cline{2-6}  & $13$ & 14 & $\langle s_{13} \rangle\simeq \Z_{14}$ &
$\chi_{(-1)}$ & 728640
\\ \hline
\end{tabular}
\end{center}
\caption{Cases of negative braiding.}\label{maintabla}
\end{table}

The proof of the theorem, as well as the unexplained notation, is
contained in sections \ref{sec:abelian} and \ref{sec:non-abelian}.
In Section \ref{preliminaries}, we set some notations and collect
preliminary results needed in the sequel.

We use GAP \cite{Sch97} to compute the character tables and other
computations, such as representatives of conjugacy classes,
intersections between centralizers and conjugacy classes, etc.
These computations are available at
\texttt{http://www.mate.uncor.edu/\~{}fantino/GAP/mathieu.htm}. In
main body of the paper the phrase ``we compute" means that we have
performed the computations with the computational algebra system
mentioned above.

\subsection{Notations}\label{subsec:notations}
We will follow the conventions in \cite{AZ,AF1}. We denote by
$\widehat{G}$ the set of isomorphism classes of irreducible
representations of a finite group $G$. We will use the rack
notation $x\trid y:=xyx^{-1}$. We denote by $\G_n$ the group of
$n$-th roots of 1 in $\ku$ and $\omega _n : = e ^ {\frac {2\pi i}
{n}}$, where $i= \sqrt{-1}$. The representation of the cyclic
group $\Z_{2n}=\langle [1] \rangle$ corresponding to
$\rho([1])=\omega_{2n}^n=-1$ will be denoted by $\chi_{(-1)}$.

For $s \in G$ we denote by $\ocs$ (resp. $G^s$) the conjugacy
class (resp. the centralizer) of $s$ in $G$. For $y \in G^s$, we
denote the conjugacy class of $y$ in the group $G^s$ by
$\oc_{y}^{G^s}$.  Also, for $k$, with $1\leq k \leq
|\widehat{G^s}|$, we denote the $k$-th conjugacy class of $G^s$ by
$\oc_{k}^{G^s}$.

In the character tables, that we give in Section
\ref{sec:abelian}, we include the following information: the first
row enumerates the conjugacy classes of the group with the
parameter $j$, the second row gives the order $|s_j|$ of a
representative of each conjugacy class, the third row give the
order of the centralizer $G^{s_j}$ of $s_j$ in the corresponding
Mathieu simple group. Notice that if all the numbers in the column
corresponding to $s_j$ are real, then $s_j$ is \emph{real}, i.~e.
$s_j^{-1} \in \oc_{s_j}$. For simplicity we will omit the cardinal
of the conjugacy classes and the order of the centralizers in some
character tables. Also, for a complex number $z$ we denote the
complex conjugate of $z$ by $z'$ (and not $\bar{z}$) for a better
reading of the tables.

\section{Preliminaries}\label{preliminaries}

\subsection{Yetter-Drinfeld modules over a finite group}\label{conventionsyd}

We recall that the irreducible Yetter-Drinfeld module $M(\oc,
\rho)$, with $\Oc$ a conjugacy class of $G$ and $\rho=(\rho,V)$ in
$\widehat{G^s}$, for a fixed element $s\in \Oc$, is described as
follows. Let $\sigma_1$, \dots, $\sigma_{N}$ be a numeration of
$\oc$ and let $g_j\in G$ such that $g_j \trid s = \sigma_j$ for
all $1\leq j \leq N$. Then $ M(\oc, \rho) = \oplus_{1\leq j \leq
N}\, g_j\otimes V$. We will write $g_jv := g_j\otimes v \in
M(\oc,\rho)$, $1\leq j \leq N$, $v\in V$. If $v\in V$ and $1\leq j
\leq N$, then the action of $g\in G$ is given by $g\cdot (g_j v) =
g_l(\gamma\cdot v)$, where $gg_j = g_l\gamma$, for some $1\leq l
\leq N$ and $\gamma\in G^s$, and the coaction is given by
$\delta(g_jv) = \sigma_j\otimes g_jv$. The Yetter-Drinfeld module
$M(\oc, \rho)$ is a braided vector space with braiding
\begin{equation} \label{yd-braiding}
c(g_jv\otimes g_kw) = \sigma_j\cdot(g_kw)\otimes g_jv =
g_l(\gamma\cdot w)\otimes g_jv,
\end{equation}
for any $1\leq j,k\leq N$, $v,w\in V$, where $\sigma_j g_k =
g_l\gamma$ for unique $l$, $1\leq l \leq N$ and $\gamma \in G^s$.
Since $s\in Z(G^s)$, the center of the group $G^s$, the Schur
Lemma implies that
\begin{equation}\label{schur} s \text{ acts by a scalar $q_{ss}$
on } V.
\end{equation}

\medbreak

A braided vector space $(W,c)$ is of \emph{diagonal type} if there
exists a basis $w_1, \dots, w_{\theta}$ of $W$ and non-zero
scalars $q_{ij}$, $1\le i,j\le \theta$, such that $c(w_i\otimes
w_j) = q_{ij} w_j\otimes w_i$, for all $1\le i,j\le \theta$. The
\emph{generalized Dynkin diagram} associated with $(W,c)$ of
diagonal type as above is the diagram with vertices $\{1,\dots,
\theta\}$, where the vertex $i$ is labelled by $q_{ii}$, and if
$q_{ij}q_{ji}\neq 1$, then the vertices $i$ and $j$ are joined by
an edge labelled by $q_{ij}q_{ji}$, i.~e.

\begin{align*}
\setlength{\unitlength}{1.4cm}
\begin{picture}(1,0)
\put(0,0){\circle*{.15}} \put(2,0){\circle*{.15}}
\put(0,0){\line(2,0){2}} \put(-0.1,0.3){\small{$q_{ii}$}}
\put(0.8,0.4){\small{$q_{ij}q_{ji}$}}
\put(1.9,0.3){\small{$q_{jj}$}}
\end{picture} \qquad \quad \quad \, \, ,
\end{align*}

\noindent see \cite{H2}. A braided vector space  $(W,c)$ of
diagonal type is of \emph{Cartan type} if $q_{ij}$ is a root of 1
for all $i,j$, $1\leq i,j \leq \theta$, $q_{ii}\neq 1$ for all
$i$, $1\le i \le \theta$, and there exists $a_{ij} \in \Z$, with
$-\ord q_{ii} < a_{ij} \leq 0$, such that $q_{ij}q_{ji} =
q_{ii}^{a_{ij}}$ for all $1\le i\neq j\le \theta$ -- see
\cite{AS2}. Set $a_{ii}=2$ for all $1\le i\le \theta$. Then
$(a_{ij})_{1\le i,j\le \theta}$ is a generalized Cartan matrix.

\subsection{Tools.} We state the principal tools that we
will use in Section \ref{sec:abelian}.

\begin{lema}\label{trivialbraiding}\cite[Remark 1.1]{AZ}.
Let $(W,c)$ be a braided vector space, $U$ a subspace of $W$ such
that $c(U\otimes U) = U\otimes U$. If $\dim \toba(U) =\infty$,
then $\dim \toba(W) =\infty$.\qed
\end{lema}

This result implies that if $\oc_{\id}$ is the conjugacy class of
the identity element of $G$ and $\rho\in \widehat{G}$, then $\dim
\toba(\oc_{\id},\rho) =\infty$. Thus, we omit to consider this
conjugacy class in the proofs of the Theorems \ref{teorM11},
\ref{teorM12}, \ref{teorM22}, \ref{teorM23} and \ref{teorM24}.

\begin{theorem}\label{cartantype} \cite[Th. 4]{H1}, see also
\cite[Th. 1.1]{AS2}. Let $(W,c)$ be a braided vector space of
Cartan type. Then $\dim \toba(W) < \infty$ if and only if the
Cartan matrix is of finite type. \qed
\end{theorem}

We say that $s\in G$ is \emph{real} if it is conjugate to
$s^{-1}$; if $s$ is real, then the conjugacy class of $s$ is also
said to be \emph{real}.

\begin{lema}\label{odd}  If $s$ is real and
$\dim\toba(\Oc_s, \rho)< \infty$, then $q_{ss} = -1$.\qed
\end{lema}

If $s^{-1}\neq s$, this is \cite[Lemma 2.2]{AZ}; if $s^2 = \id$,
then $q_{ss} = \pm 1$, but $q_{ss} = 1$ is excluded by Lemma
\ref{trivialbraiding}. Notice that $q_{ss} = -1$ implies that $s$
has even order. The following is a generalization of Lemma
\ref{odd}. See \cite[Lemmata 1.8 and 1.9]{AF2} or \cite[Corollary
2.2]{FGV}.

\begin{lema}\label{lemaB}
Let $G$ be a finite group, $s\in G$, $\Oc$ the conjugacy class of
$s$ and $\rho=(\rho,V) \in \widehat{G^s}$ such that
$\dim\toba(\Oc,\rho)<\infty$. Assume that there exists an integer
$k$ such that $s^k\in \Oc$ and $s^k\neq s$.
\begin{itemize}
\item[(a)] If $\deg \rho >1$, then $q_{ss}=-1$.
\item[(b)] If $\deg \rho =1$, then either $q_{ss}=-1$ or $q_{ss}$
$\in \G_3-1$.
\end{itemize}
On the other hand, if $s^{k^2}\neq s$, then $q_{ss}=-1$.\qed
\end{lema}

The next important tool follows from \cite{H2}.

\begin{lema}\label{Hecke}
Let $W$ be a Yetter-Drinfeld module, $U \subseteq W$ a braided
vector subspace of diagonal type of $W$ such that $q_{ii}$ is a
root of 1 for all $i$, and let $\mathcal G$ be the generalized
Dynkin diagram corresponding to $U$. If $\mathcal G$ contains an
$r$-cycle with $r>3$, or a vertex with valency greater than 3,
then the Nichols algebra $\toba(U)$ is infinite-dimensional.
Hence, $\dim\toba(W)=\infty$. \qed
\end{lema}

\subsection*{Abelian subspaces of a braided vector space}\label{cartan-subsection}
As in \cite{AF1,AF2}, in a first step we look for braided
subspaces $W$ of diagonal type of $M(\oc, \rho)$ whose Nichols
algebra is infinite-dimensional.

Let  $(X, \trid)$ be a rack, see for example \cite{AG2}.  Let $q:
X\times X\to \ku^{\times}$ be a rack 2-cocycle and let $(\ku X,
c_q)$ be the associated braided vector space, that is $\ku X$ is a
vector space with a basis $e_k$, $k\in X$, and $c_q(e_k\otimes
e_l) = q_{k,l} e_{k\trid l}\otimes e_k$, for all $k$, $l\in X$.
Let us say that a subrack $T$ of $X$  is \emph{abelian} if $k\trid
l = l$ for all $k$, $l\in T$. If $T$ is an abelian subrack of $X$
then $\ku T$ is a braided vector subspace of $(\ku X, c_q)$ of
diagonal type.

We shall say that $(\ku X, c_q)$ is \emph{negative}\label{braineg}
if for any abelian subrack $T$ of $X$ $q_{kk} = -1$ and
$q_{kl}q_{lk} = 1$ for all $k,l\in T$ (hence $\toba(\ku T)$ is
twist-equivalent to an exterior algebra and $\dim \toba(\ku T) =
2^{\card T}$).

Let $G$ be a finite group, $\oc$ a conjugacy class in $G$,
$\rho\in \widehat{G^s}$, with $s\in \oc$ fixed. As in subsection
\ref{conventionsyd}, we fix a numeration $\sigma_1 = s$, \dots,
$\sigma_{N}$ of $\oc$ and $g_k\in G$ such that $g_k \trid s =
\sigma_k$ for all $1\leq k \leq N$. Let $I\subset \{1, \dots, N\}$
and $T := \{\sigma_k: k\in I\}$. We characterize when $T$ is an
abelian subrack of $\oc$. Let
\begin{align}\label{eq:gamma}
\gamma_{k,l} := g_l^{-1} \sigma_k g_l, \qquad \text{ $k$, $l\in
I$.}
\end{align}
Then the following are equivalent:
\begin{align*}
\text{(a)\,\,$\sigma_k\trid \sigma_l = \sigma_l$ (i.~e. $\sigma_k$
and $\sigma_l$ commute) \qquad and \qquad (b)\,\, $\gamma_{k,l}
\in G^s$.}
\end{align*}
Assume that (a) (or, equivalently, (b)) holds for all $k$, $l\in
I$; then $\gamma_{k,l}\in \oc_s \cap G^s$. Let $V$ be the vector
space affording $\rho$.  Let $v_1, \dots, v_R$ be simultaneous
eigenvectors of $\gamma_{k,l}$, $k$, $l\in I$. Because of
\eqref{yd-braiding}, we have that
$$
W = \ku-\text{span of } g_k v_j, \quad k\in I, \, 1\le j \le R,
$$
is a braided subspace of diagonal type of dimension $R \,\card T$.
Notice that $R$ depends not only on $T$ but also on the
representation $\rho$; for instance if $\rho$ is a character then
$R = 1 = \dim V$, and $M(\oc, \rho)$ is of rack type.

\begin{lema}\label{lema:nodep}
Assume that $\sigma_k$, $\sigma_l\in \oc$ commute and
$\deg\rho=1$. Then the scalar $\rho(\gamma_{k,l})$ does not depend
on $g_k$ and $g_l$.
\end{lema}
\pf Let $\widetilde{g_k}$, $\widetilde{g_l}\in G$ such that
$\widetilde{g_k}\trid s=\sigma_k$ and $\widetilde{g_l}\trid
s=\sigma_l$. Thus, $\widetilde{g_k}=g_k \eta_k$ and
$\widetilde{g_l}=g_l \eta_l$, with $\eta_k$, $\eta_l\in G^s$. If
we call $\widetilde{\gamma_{k,l}}:=\widetilde{g_l}^{-1}
\widetilde{g_k} \, s \, \widetilde{g_k}^{-1} \widetilde{g_l}$,
then
\begin{align*}
\widetilde{\gamma_{k,l}}=\eta_l^{-1} \,  g_l^{-1} \sigma_k  g_l \,
\eta_l= \eta_l^{-1} \, \gamma_{k,l} \, \eta_l.
\end{align*}
Hence, $\rho(\widetilde{\gamma_{k,l}})=\rho(\eta_l)^{-1}
\rho(\gamma_{k,l})\rho(\eta_l)=\rho(\gamma_{k,l})$, since $\deg
\rho=1$. \epf

In view of this result, we can choose $g_1=\id$, the identity of
the group $G$.

\begin{obs}
If $\deg\rho=1$, then the condition of negative braiding is
equivalent to (i) $\rho(\gamma_{k,k})=-1$ and
\begin{itemize}
\item[(ii)] for every commuting pair $\sigma_k$, $\sigma_l\in \oc$, it holds
$\rho(\gamma_{k,l}\gamma_{l,k})=1$.
\end{itemize}
\end{obs}

The next result is useful in order to prove that $M(\oc, \rho)$
has negative braiding when $\rho$ is an one-dimensional
representation.

\begin{lema}\label{lema:treneg}
The condition \emph{(ii)} given above is equivalent to
$$
\emph{(ii)'} \qquad \text{for every $\sigma_t\in \oc \cap G^s$, it
holds $\rho(\gamma_{1,t}\gamma_{t,1})=1$.}
$$
\end{lema}

\pf Obviously, (ii) implies (ii)'. Reciprocally, assume that (ii)'
holds. Let $\sigma_k$, $\sigma_l\in \oc$ that commute. Then,
$\gamma_{k,l}=g_l^{-1} g_k \, s \, g_k^{-1} g_l$,
$\gamma_{l,k}=g_k^{-1} g_l \, s \, g_l^{-1} g_k$ are in $\oc\cap
G^s$. Hence, $\gamma_{k,l}=\sigma_t$, for some $1 \leq t \leq N$.
By Lemma \ref{lema:nodep},
$\rho(\gamma_{k,l}\gamma_{l,k})=\rho(\gamma_{1,t}\gamma_{t,1})$,
and the result follows.\epf


\subsection{Criterions from non-abelian subracks}\label{subsec:nonabeltech}
We will mention here some criterions that allow to decide the
dimension of the Nichols algebra $\toba(\oc,\rho)$ using
non-abelian subracks of $\oc$. These criterions were developed in
\cite{AF3} using important results of \cite{AHS}.

Let $p>1$ be an integer. A family $(\sigma_i)_{i\in \Z_p}$ of
distinct elements of a group $G$ is of \emph{type $\D_p$} if
$\sigma_i \trid \sigma_j = \sigma_{2i - j}$, $i,j\in \Z_p$. Let
$(\sigma_i)_{i\in \Z_p}$ and $(\tau_i)_{i\in \Z_p}$ be two
families of type $\D_p$ in $G$, such that  $\sigma_i\neq\tau_j$
for all $i,j\in \Z_p$, we say that $(\sigma, \tau) :=
(\sigma_i)_{i\in \Z_p}\cup (\tau_i)_{i\in \Z_p}$ is of \emph{type
$\D_p^{(2)}$} if
\begin{equation}\label{eqn:dp2}
\sigma_i \trid \tau_j = \tau_{2i - j}, \quad\tau_i \trid \sigma_j
= \sigma_{2i - j}, \quad i,j\in \Z_p.
\end{equation}

\begin{lema}\label{coro:dp-cor} \cite[Cor. 2.9]{AF3}
Let $p$ be an odd prime, $(\sigma_i)_{i \in \Z_p}$ a family of
type $\D_p$ in a finite group $G$ and $\rho\in
\widehat{G^{\,\sigma_0}}$. If there exists $k$, such that
$\sigma_0^k\in \oc$, the conjugacy class of $\sigma_0$, and
$q_{\sigma_0\sigma_0}=-1$, then $\dim\toba(\oc,\rho)=\infty$.\qed
\end{lema}

Let $\oc_4^4$ be the conjugacy class of the $4$-cycles in the
symmetric group $\s_4$. We will say that a family
$(\sigma_i)_{1\leq i \leq 6}$ of distinct elements of a group $G$
is of \emph{type $\oct$} if $(\sigma_i)_{1\leq i \leq 6}$ form a
rack isomorphic to $\oc_4^4$. We call such a rack an
\emph{octahedral} rack. Let $\sigma_i$, $\tau_i\in G$, $1\leq i
\leq 6$, all distinct; we say  that $(\sigma, \tau) :=
(\sigma_i)_{i}\cup (\tau_i)_{i}$ is of \emph{type $\oct^{(2)}$} if
$(\sigma_i)_{i}$ and $(\tau_i)_{i}$ are both of type $\oct$, and
\begin{align}\label{R4^2}
\sigma_i\trid \tau_j=\tau_{i\trid j}, \quad  \tau_i\trid
\sigma_j=\sigma_{i\trid j}, \quad 1\leq i,j \leq 6,
\end{align}
where $\trid$ in the subindex denotes the operation of rack in the
octahedral rack.

We state the main tool from this non-abelian rack -- see \cite[Th.
4.11]{AF3}.

\begin{lema}\label{teor:aplicAHSch2}
Let $G$ be a finite group, $(\sigma_i)_{i}\cup (\tau_i)_{i}$ a
family of type $\oct^{(2)}$ in $\oc$ the conjugacy class of
$\sigma_1$, and $g\in G$ such that $g\trid \sigma_1 =\tau_1$. Let
$\rho=(\rho,V)\in \widehat{G^{\sigma_1}}$ with
$q_{\sigma_1\sigma_1}=-1$. If there exist $v$, $w\in V-0$ such
that
\begin{align*}
\rho(\sigma_6) v&=-v, &\rho(g^{-1}\sigma_1g) w&=-w,\\
\rho(\tau_1) v&=-v, &\rho(g^{-1}\sigma_6g) w&=-w,
\end{align*}
then $\dim\toba(\oc,\rho)=\infty$.\qed
\end{lema}

The following is a useful consequence of this result.

\begin{lema}\label{co:especial2}\cite[Cor. 4.12]{AF3}
Let $G$ be a finite group, $(\sigma_i)_{i}\cup (\tau_i)_{i}$ a
family of type $\oct^{(2)}$ in $\oc$ the conjugacy class of
$\sigma_1$, and $\rho\in \widehat{G^{\sigma_1}}$, with
$q_{\sigma_1\sigma_1}=-1$. If $\sigma_6=\sigma_1^d$ and
$\tau_1=\sigma_1^e$, then $\dim\toba(\oc,\rho)=\infty$.\qed
\end{lema}


\section{Using techniques based on abelian
subracks}\label{sec:abelian}

In this section, we will determine the irreducible Yetter-Drinfeld
modules $M(\oc,\rho)$ with $\dim \toba(\oc,\rho)=\infty$ by mean
of abelian subracks of $\oc$. We will consider each simple Mathieu
group separately.


\subsection{The group $M_{11}$}\label{sectionM11}

The Mathieu simple group $M_{11}$ can be given as a subgroup of
$\s_{11}$ in the following form
\begin{align*}
M_{11}:=\langle \, ( 1, 2, 3, 4, 5, 6, 7, 8, 9,10,11), ( 3, 7,11,
8)(4,10, 5, 6) \, \rangle.
\end{align*}
In Table \ref{tablacarM11}, we show the character table of
$M_{11}$, where $A = (-1-i\sqrt{11})/2$, $B =i \sqrt{2}$. We take
the following elements to be the representatives of the conjugacy
classes of $M_{11}$:
\begin{align*}
s_1&:=\id,  &s_2&:= ( 1, 9, 7,10, 8,11, 5, 4, 3, 6, 2),\\
s_3&:= ( 1, 7, 8, 5, 3, 2, 9,10,11,4, 6), &
s_4&:=( 1, 8, 2, 7)( 4, 6,10, 5),\\
s_5&:= ( 1, 2)( 4,10)( 5, 6)( 7, 8),&
s_6&:= ( 1, 3,11, 6, 7,10, 4, 5)( 8, 9), \\
s_7&:=( 1,10,11, 5, 7, 3, 4, 6)( 8, 9),&s_8&:=( 1, 6, 4)( 2, 9, 7)( 8,11,10),\\
s_9&:=( 1, 2, 3, 4, 8)( 5,10, 7,11, 6),& s_{10}&:=( 1, 5, 8, 4, 6,
9)( 2,10, 3)( 7,11).
\end{align*}

\begin{table}[t]
\begin{center}
\begin{tabular}{|c|c|c|c|c|c|c|c|c|c|c|}
\hline {\bf $j$} &    1& 2& 3& 4& 5& 6& 7& 8& 9& 10\\
\hline \hline {\bf $|s_j|$} &    1& 11& 11& 4& 2& 8& 8& 3& 5& 6\\
\hline {\bf $|G^{s_j}|$} &   7920& 11& 11& 8& 48& 8& 8& 18& 5& 6 \\
\hline {\bf $|\Oc_{s_j}|$} &   1& 720& 720& 990& 165& 990& 990& 440& 1584& 1320\\
\hline \hline

\hline {\bf $\chi_1$} &     1& 1& 1& 1& 1& 1& 1& 1& 1& 1\\

\hline {\bf $\chi_2$} &    10& -1& -1& 2& 2& 0& 0& 1&0& -1\\

\hline {\bf $\chi_3$} &    10& -1& -1& 0& -2& $B$& $B'$ & 1& 0& 1 \\

\hline {\bf $\chi_4$} &     10& -1& -1& 0& -2& $B'$& $B$& 1& 0& 1\\

\hline {\bf $\chi_5$} &     11& 0& 0& -1& 3& -1& -1& 2& 1& 0 \\

\hline {\bf $\chi_6$} &   16& $A$&
      $A'$& 0& 0& 0& 0& -2& 1& 0 \\

\hline {\bf $\chi_7$} &    16& $A'$&
      $A$& 0& 0& 0& 0& -2& 1& 0 \\

\hline {\bf $\chi_8$} &     44& 0& 0& 0& 4& 0& 0& -1& -1& 1 \\

\hline {\bf $\chi_9$} &  45& 1& 1& 1& -3& -1& -1& 0& 0& 0  \\

\hline {\bf $\chi_{10}$} & 55& 0& 0& -1& -1& 1& 1& 1& 0& -1 \\

\hline
\end{tabular}
\end{center}
\caption{Character table of $M_{11}$.}\label{tablacarM11}
\end{table}

In the following statement, we summarize our study by mean of
abelian subracks in the group $M_{11}$.

\begin{theorem}\label{teorM11}
Let $\rho\in \widehat{M_{11}^{s_j}}$, with $1\leq j \leq 10$. The
braiding is negative in the cases $j=4$, with $\rho=\nu_2$ or
$\nu_6$, $j=6$, $7$ and $10$, with $\rho=\chi_{(-1)}$. Otherwise,
$\dim\toba(\oc_{s_j},\rho)=\infty$.
\end{theorem}

\pf  From Table \ref{tablacarM11}, we see that for $j=4$, $5$,
$8$, $9$ and $10$, $s_j$ is real.

\emph{CASE: $j=8$, $9$}. By Lemma \ref{odd}, $\dim\toba(\oc_{s_j},
\rho)=\infty$, for all $\rho \in \widehat{M_{11}^{s_j}}$.

\smallbreak

\emph{CASE: $j=2$, $3$}. We compute that $s_j^3$ and $s_j^{9}$ are
in $\oc_{s_j}$, and $s_j^3\neq s_j^{9}$. By Lemma \ref{lemaB},
$\dim\toba(\oc_{s_j}, \rho)=\infty$, for all $\rho \in
\widehat{M_{11}^{s_j}}$, since $|s_j|=11$.

\smallbreak

\emph{CASE: $j=4$}. The element $s_4$ is real and we compute that
$$M_{11}^{s_4}=\langle x:=( 1, 4, 7, 5, 2,10, 8, 6)(
3,11)\rangle\simeq \Z_8,$$ with $x^6=s_4$. We set
$\widehat{M_{11}^{s_4}}=\{\nu_0,\dots,\nu_7\}$, where
$\nu_l(x):=\omega_8^l$, $0 \leq l \leq 7$. Clearly, if $l=0$, $1$,
$4$, $5$ or $7$, then $q_{s_{4}s_{4}}\neq -1$, and
$\dim\toba(\oc_{s_{4}}, \nu_l)=\infty$, by Lemma \ref{odd}. The
remained two cases correspond to $l=2$, $6$. We compute that
$\oc_{s_4}\cap M_{11}^{s_4}=\{s_4,s_4^{-1}\}$. It is easy to see
that the braiding is negative.

\smallbreak

\emph{CASE: $j=6$ or $7$}. We compute that $M_{11}^{s_j}=\langle
s_j \rangle\simeq \Z_8$, $s_j^3$ is in $\oc_{s_j}$. Since $3$ does
not divide $|s_j|=8$ we have that if $q_{s_js_j}\neq -1$, then
$\dim\toba(\oc_{s_j}, \rho)=\infty$, by Lemma \ref{lemaB}. The
remained case corresponds to $\rho(s_j)=\omega_8^4=-1$, which
satisfies $q_{s_js_j}= -1$. We compute that $\oc_{s_j}\cap
M_{11}^{s_j}=\{s_j,s_j^3\}$. It is straightforward to prove that
the braiding is negative.

\smallbreak

\emph{CASE: $j=10$}. The element $s_{10}$ is real and we compute
that $M_{11}^{s_{10}}=\langle s_{10} \rangle\simeq \Z_6$. Now, we
have that if $q_{s_{10}s_{10}}\neq -1$, then
$\dim\toba(\oc_{s_{10}}, \rho)=\infty$, by Lemma \ref{odd}. The
remained case corresponds to $\rho(s_{10})=\omega_6^3=-1$, which
satisfies $q_{s_{10}s_{10}}= -1$. We compute that
$\oc_{s_{10}}\cap M_{11}^{s_{10}}=\{s_{10},s_{10}^{-1}\}$. Then
the braiding is negative.

\smallbreak

\emph{CASE: $j=5$}. We compute that $M_{11}^{s_5}$ is a
non-abelian group of order 48, whose character table is given by
Table \ref{tablacarM11^s5}.

For every $k$, $1 \leq k \leq 8$, we call $\rho_k=(\rho_k,V_k)$
the irreducible representation of $M_{11}^{s_5}$ whose character
is $\mu_k$. From Table \ref{tablacarM11^s5}, we can deduce that if
$k\neq 4$, $5$, $8$, then $q_{s_5s_5}\neq -1$ and
$\dim\toba(\oc_{s_{5}}, \rho_k)=\infty$, by Lemma \ref{odd}.

On the other hand, if $k= 4$, $5$ or $8$, then $q_{s_5s_5}= -1$.
For these three cases we will prove that $\dim\toba(\oc_{s_{5}},
\rho_k)=\infty$. First, we compute that $\oc_{s_5}\cap
M_{11}^{s_5}$ has 13 elements and it contains $\sigma_1:=s_5$,
$\sigma_2:=( 4,10)( 5, 8)( 6, 7)( 9,11)$ and $\sigma_3:=( 1, 2)(
5, 7)( 6,8)( 9,11)$. Notice that these elements commute each other
and $\sigma_2\sigma_3=s_5$. Also, we compute that $\sigma_2$,
$\sigma_3\in \oc_{2}^{M_{11}^{s_5}}$ -- see Subsection
\ref{subsec:notations}. Now, we choose $g_1:=\id$,
\begin{align*}
g_2:=( 1, 9)( 2,11)( 4,10)( 5, 7) \quad \text{ and } \quad g_3:= (
1, 2)( 4, 9)( 6, 7)(10,11).
\end{align*}
These elements are in $M_{11}$ and they satisfy
\begin{align}
\label{eq1} \sigma_1 g_1&=g_1 \sigma_1, \quad & \sigma_1 g_2&=g_2
\sigma_2,\quad &\sigma_1 g_3&=g_3 \sigma_3,\\
\label{eq2} \sigma_2 g_1&=g_1 \sigma_2, \quad & \sigma_2 g_2&=g_2
\sigma_1,\quad &\sigma_2 g_3&=g_3 \sigma_2,\\
\label{eq3} \sigma_3 g_1&=g_1 \sigma_3, \quad & \sigma_3 g_2&=g_2
\sigma_3, \quad &\sigma_3 g_3&=g_3 \sigma_1.
\end{align}

\begin{table}[t]
\begin{center}
\begin{tabular}{|c|c|c|c|c|c|c|c|c|}
\hline {\bf $k$} & 1 & 2& 3& 4& 5& 6& 7& 8  \\
\hline \hline {\bf $|y_k|$} &  1& 2& 3& 2& 6& 4& 8& 8\\
\hline {\bf $|G^{y_k}|$} &   48& 4& 6& 48& 6& 8& 8& 8 \\
\hline {\bf $|\Oc_{y_k}|$} & 1& 12& 8& 1& 8& 6& 6& 6  \\
\hline \hline

\hline {\bf $\mu_1$} &    1& 1&1& 1& 1& 1& 1& 1  \\

\hline {\bf $\mu_2$} &    1& -1& 1& 1& 1& 1& -1& -1 \\

\hline {\bf $\mu_3$} &   2& 0& -1& 2& -1& 2& 0& 0 \\

\hline {\bf $\mu_4$} &   2& 0& -1& -2& 1& 0& $i \sqrt{2}$& -$i \sqrt{2}$ \\

\hline {\bf $\mu_5$} &     2& 0& -1& -2& 1& 0& -$i \sqrt{2}$ & $i \sqrt{2}$\\

\hline {\bf $\mu_6$} &  3& -1& 0& 3& 0& -1& 1& 1  \\

\hline {\bf $\mu_7$} &  3& 1& 0& 3& 0& -1& -1& -1  \\

\hline {\bf $\mu_8$} &  4& 0& 1& -4& -1& 0& 0& 0 \\

\hline
\end{tabular}
\end{center}
\caption{Character table of $M_{11}^{s_5}$.}\label{tablacarM11^s5}
\end{table}

Assume that $k=4$. Since $\sigma_1$, $\sigma_2$ and $\sigma_3$
commute there exists a basis $\{v_1, v_2\}$ of $V_4$, the vector
space affording $\rho_4$, composed by simultaneous eigenvectors of
$\rho_4(\sigma_1)$, $\rho_4(\sigma_2)$ and $\rho_4(\sigma_3)$. Let
us say $\rho_4(\sigma_2)v_l=\lambda_l v_l$ and
$\rho_4(\sigma_3)v_l=\kappa_l v_l$, $l=1$, $2$. Notice that
$\lambda_l$, $\kappa_l=\pm 1$, $l=1$, $2$, due to
$|\sigma_2|=2=|\sigma_3|$. Moreover, since $\sigma_2\sigma_3=s_5$
we have that $\lambda_l\kappa_l=-1$, $l=1$, $2$. From Table
\ref{tablacarM11^s5}, we can deduce that $\lambda_1+\lambda_2=0$
because $\sigma_2\in \oc_{2}^{M_{11}^{s_5}}$ and
$\mu_4(\oc_{2}^{M_{11}^{s_5}})=0$. Reordering the basis we can
suppose that $\lambda_1=1=-\lambda_2$. We define $W:=\ku$ - span
of $\{g_1v_1,g_2v_2,g_3v_2\}$. Hence, $W$ is a braided vector
subspace of $M(\oc_{s_5},\rho)$ of Cartan type. Indeed, it is
straightforward to compute that the matrix of coefficients
$\mathcal Q$ is
\begin{align}\label{matrix:coef:3x3a}
\begin{pmatrix}
-1 &-1 & 1\\
1 &-1 & 1\\
-1 &-1 & -1\\
\end{pmatrix}.
\end{align}
The corresponding Cartan matrix is given by
\begin{align}\label{cartanmatrix:3x3}
\mathcal A =\begin{pmatrix}
2 &  -1 &  -1 \\
-1 & 2 & -1 \\
-1 & -1 & 2 \\
\end{pmatrix}.
\end{align}
By Theorem \ref{cartantype}, $\dim \toba(\oc_{s_5}, \rho_4) =
\infty$.

\smallbreak

The case $k=5$ is analogous to the case $k=4$ because
$\mu_5(\oc_{2}^{M_{11}^{s_5}})=0$.

\smallbreak

Finally, the case $k=8$ can be reduced to the previous cases.
Indeed, since $\sigma_1$, $\sigma_2$ and $\sigma_3$ commute there
exists a basis $\{v_1, v_2, v_3, v_4\}$ of $V_8$, the vector space
affording $\rho_8$, composed by simultaneous eigenvectors of
$\rho_8(\sigma_1)$, $\rho_8(\sigma_2)$ and $\rho_8(\sigma_3)$. Let
us say $\rho_8(\sigma_2)v_l=\lambda_l v_l$ and
$\rho_8(\sigma_3)v_l=\kappa_l v_l$, $1\leq l \leq 4$, where
$\lambda_l$, $\kappa_l=\pm 1$, $1\leq l \leq 4$, due to
$|\sigma_2|=2=|\sigma_3|$. From Table \ref{tablacarM11^s5}, we can
deduce that
$\lambda_1+\lambda_2+\lambda_3+\lambda_4=0=\kappa_1+\kappa_2+\kappa_3+\kappa_4$
because $\mu_8(\oc_{2}^{M_{11}^{s_5}})=0$. This implies that there
exist $r$, $t\in \{1,2,3,4\}$ such that $\lambda_r=1=-\lambda_t$.
Now, if we define $W:=\ku$ - span of $\{g_1v_r,g_2v_t,g_3v_t\}$,
then $W$ is a braided vector subspace of $M(\oc_{s_5},\rho)$ of
Cartan type, with Cartan matrix given by \eqref{cartanmatrix:3x3}.
Therefore, $\dim \toba(\oc_{s_5}, \rho_8) = \infty$. \epf

\begin{obs}
The group $M_{11}^{s_8}$ has 9 conjugacy classes. So, we can
point out the following fact: there are 84 possible
pairs $(\oc,\rho)$ for $M_{11}$; 79 of them have infinite-dimensional Nichols
algebras, and 5 have negative braiding.
\end{obs}


\subsection{The group $M_{12}$}\label{sectionM12}

The Mathieu simple group $M_{12}$ can be given as a subgroup of
$\s_{12}$ in the following form
\begin{align*}
M_{12}:=\langle & \,( 1, 2, 3, 4, 5, 6, 7, 8, 9,10,11), ( 3, 7,11,
8)( 4,10, 5, 6),\\
& ( 1,12) ( 2,11)( 3, 6)( 4, 8)( 5, 9)( 7,10) \, \rangle.
\end{align*}
In Table \ref{tablacarM12}, we show the character table of
$M_{12}$, with $A = (-1-i\sqrt{11})/2$. We take the following
elements to be the representatives of the conjugacy classes of
$M_{12}$: $s_1:=\id,$ {\footnotesize{\begin{align*}
  s_2&:= ( 1, 8,12)( 2, 3,11, 9,10, 6)( 4, 5), &
s_3&:= ( 1,12, 8)( 2,11,10)( 3, 9,6), \\
s_4&:=( 2, 9)( 3,10)( 4, 5)( 6,11),& s_5&:=( 1,12, 7, 4)( 2, 9,10,
5,11, 3, 6, 8), \\
s_6&:=( 1, 7)( 2,10,11, 6)( 3, 8, 9, 5)( 4,12), &
s_7&:=( 1, 9, 4, 2,11, 8)( 3,10,12, 5, 6, 7),\\
s_8&:=( 1, 4,11)( 2, 8, 9)( 3,12, 6)( 5, 7,10),& s_9&:= ( 1, 2)(
3, 5)( 4, 8)( 6,10)( 7,12)( 9,11),\\  s_{10}&:=( 1, 7, 2, 6, 5)(
3, 9,12,10,11),& s_{11}&:= ( 2,10, 6, 4,12, 5, 7, 3, 8,11, 9),\\
 s_{12}&:= ( 2,6,12, 7, 8, 9,10, 4, 5, 3,11),&
s_{13}&:= ( 1, 2, 7,10, 5, 6, 8,12, 9, 4)( 3,11),\\ s_{14}&:=(
1,10,11, 8, 4, 5,12, 3)( 6, 9),& s_{15}&:=  ( 1,11, 4,12)( 3,10,
8, 5).
\end{align*}}}

{\footnotesize{
\begin{table}[t]
\begin{center}
\tiny{\begin{tabular}{|c|c|c|c|c|c|c|c|c|c|c|c|c|c|c|c|} \hline
{\bf $j$} &    1& 2& 3& 4& 5& 6& 7& 8& 9& 10& 11& 12& 13& 14& 15\\
\hline \hline
{\bf $|s_j|$} &    1& 6& 3& 2& 8& 4& 6& 3& 2& 5& 11& 11& 10& 8& 4\\
\hline {\bf $|G^{s_j}|$} &    95040& 6& 54& 192& 8& 32& 12& 36& 240& 10& 11& 11& 10& 8& 32 \\
\hline \hline

\hline {\bf $\chi_1$} &  1& 1& 1& 1& 1& 1& 1& 1& 1& 1& 1& 1& 1& 1& 1\\

\hline {\bf $\chi_2$} &  11& 0& 2& 3& 1& 3& -1& -1& -1& 1& 0& 0& -1& -1& -1\\

\hline {\bf $\chi_3$} &   11& 0& 2& 3& -1& -1& -1& -1& -1& 1& 0& 0& -1& 1& 3\\

\hline {\bf $\chi_4$} &  16& 0& -2& 0& 0& 0& 1& 1& 4& 1& A' & A &
-1& 0& 0 \\

\hline {\bf $\chi_5$} &  16& 0& -2& 0& 0& 0& 1& 1& 4& 1& A& A'& -1& 0& 0\\

\hline {\bf $\chi_6$} &  45& 0& 0& -3& -1& 1& -1& 3& 5& 0& 1& 1& 0& -1& 1\\

\hline {\bf $\chi_7$} &  54& 0& 0& 6& 0& 2& 0& 0& 6& -1& -1& -1& 1& 0& 2\\

\hline {\bf $\chi_8$} & 55& 1& 1& 7& -1& -1& 1& 1& -5& 0& 0& 0& 0& -1& -1\\

\hline {\bf $\chi_9$} & 55& -1& 1& -1& -1& 3& 1& 1& -5& 0& 0& 0& 0& 1& -1\\

\hline {\bf $\chi_{10}$} & 55& -1& 1& -1& 1& -1& 1& 1& -5& 0& 0& 0& 0& -1& 3\\

\hline {\bf $\chi_{11}$} & 66& -1& 3& 2& 0& -2& 0& 0& 6& 1& 0& 0& 1& 0& -2\\

\hline {\bf $\chi_{12}$} &  99& 0& 0& 3& 1& -1& -1& 3& -1& -1& 0& 0& -1& 1& -1\\

\hline {\bf $\chi_{13}$}& 120& 1& 3& -8& 0& 0& 0& 0& 0& 0& -1& -1& 0& 0& 0\\

\hline {\bf $\chi_{14}$} &  144& 0& 0& 0& 0& 0& 1& -3& 4& -1& 1& 1& -1& 0& 0 \\

\hline {\bf $\chi_{15}$}&  176& 0& -4& 0& 0& 0& -1& -1& -4& 1& 0& 0& 1& 0& 0 \\

\hline
\end{tabular}}
\end{center}
\caption{Character table of $M_{12}$.}\label{tablacarM12}
\end{table}}}

In the following statement, we summarize our study by mean of
abelian subracks in the group $M_{12}$.

\begin{theorem}\label{teorM12}
Let $\rho\in \widehat{M_{12}^{s_j}}$, with $1\leq j \leq 15$. The
braiding is negative in the cases $j=2$, $5$, $13$ and $14$, with
$\rho=\chi_{(-1)}$. Otherwise,
$\dim\toba(\oc_{s_j},\rho)\!=\!\infty$.
\end{theorem}

\pf \emph{CASE: $j=3$, $8$, $10$}. From Table \ref{tablacarM12},
we see that $s_j$ is real. By Lemma \ref{odd},
$\dim\toba(\oc_{s_j}, \rho)=\infty$, for all $\rho \in
\widehat{M_{12}^{s_j}}$.

\smallbreak

\emph{CASE: $j=11$, $12$}. We compute that $s_j^3$ and $s_j^{9}$
are in $\oc_{s_j}$, and $s_j^3\neq s_j^{9}$. By Lemma \ref{lemaB},
$\dim\toba(\oc_{s_j}, \rho)=\infty$, for all $\rho \in
\widehat{M_{12}^{s_j}}$, since $|s_j|=11$.

\smallbreak

\emph{CASE: $j=2$}. The element $s_{2}$ is real and we compute
that $M_{12}^{s_{2}}=\langle s_{2} \rangle\simeq \Z_6$. Now, if
$q_{s_{2}s_{2}}\neq -1$, then $\dim\toba(\oc_{s_{2}},
\rho)=\infty$, by Lemma \ref{odd}. The remained case corresponds
to $\rho(s_{2})=\omega_6^3=-1$, which satisfies $q_{s_{2}s_{2}}=
-1$. We compute that $\oc_{s_{2}}\cap
M_{12}^{s_{2}}=\{s_{2},s_{2}^{-1}\}$, and that the braiding is
negative.

\smallbreak

\emph{CASE: $j=5$, $14$}. The element $s_{j}$ is real and
$M_{12}^{s_{j}}=\langle s_{j} \rangle\simeq \Z_8$. If
$q_{s_{j}s_{j}}\neq -1$, then $\dim\toba(\oc_{s_{j}},
\rho)=\infty$, by Lemma \ref{odd}. The remained case corresponds
to $\rho(s_{j})=\omega_8^4=-1$, which satisfies $q_{s_{j}s_{j}}=
-1$. We compute that $\oc_{s_{j}}\cap
M_{12}^{s_{j}}=\{s_{j},s_{j}^3,s_{j}^5,s_{j}^7\}$, and that the
braiding is negative.

\smallbreak

\emph{CASE: $j=13$}. The element $s_{13}$ is real and
$M_{12}^{s_{13}}=\langle s_{13} \rangle\simeq \Z_{10}$. If
$q_{s_{13}s_{13}}\neq -1$, then $\dim\toba(\oc_{s_{13}},
\rho)=\infty$, by Lemma \ref{odd}. The remained case corresponds
to $\rho(s_{13})=\omega_{10}^5=-1$, which satisfies
$q_{s_{13}s_{13}}= -1$. We compute that $\oc_{s_{13}}\cap
M_{12}^{s_{13}}=\{s_{13},s_{13}^3,s_{13}^7,s_{13}^9\}$, and that
the braiding is negative.

\smallbreak

\emph{CASE: $j=7$}. The element $s_7$ is real and we compute that
$$M_{12}^{s_7}=\langle x, s_{7} \rangle\simeq \Z_2 \times \Z_6,$$
with $x:=( 1,12)( 2, 7)( 3,11)( 4, 6)( 5, 9)( 8,10)$. Let us
define $\{\nu_0,\dots,\nu_5\}$, where $\nu_l(s_7):=\omega_6^l$, $0
\leq l \leq 5$. So,
\begin{align*}
\widehat{M_{12}^{s_7}}=\{\epsilon \otimes \nu_l,\, \sgn \otimes
\nu_l \,\, | \,\, 0 \leq l\leq 5 \},
\end{align*}
where $\epsilon$ and $\sgn$ mean the trivial and the sign
representations of $\Z_2$, respectively. Clearly, if $\rho \in
\widehat{M_{12}^{s_7}}$, with $l\neq 3$, then $q_{s_{7}s_{7}}\neq
-1$, and $\dim\toba(\oc_{s_{7}}, \rho)=\infty$, by Lemma
\ref{odd}. The remained two cases are $\epsilon \otimes \nu_3$ and
$\sgn \otimes \nu_3$. We will prove that also the Nichols algebra
$\toba(\oc_{s_{7}}, \rho)$ is infinite-dimensional. We compute
that $\oc_{s_7}\cap M_{12}^{s_7}$ has 6 elements and it contains
$\sigma_1:=s_7$, {\small{\begin{align*} \sigma_2&:=( 1, 5, 4,
7,11,10)( 2, 3, 8,12, 9, 6), &\sigma_3&:=( 1, 3, 4,12,11, 6)( 2,
5, 8, 7, 9,10).
\end{align*}}}

\vspace*{-0.5cm}\noindent Notice that $\sigma_2\sigma_3=s_7^{-1}$.
We take $g_1:=\id$, $g_2:=( 2, 7,12)( 3, 9, 5)(6,8,10)$ and
$g_3:=g_2^{-1}$. We compute that the following relations hold
\begin{align}
\label{eq1bis} \sigma_1 g_1&=g_1 \sigma_1, \quad & \sigma_1
g_2&=g_2 \sigma_3,\quad &\sigma_1 g_3&=g_3 \sigma_2,\\
\label{eq2bis} \sigma_2 g_1&=g_1 \sigma_2, \quad & \sigma_2
g_2&=g_2
\sigma_1,\quad &\sigma_2 g_3&=g_3 \sigma_3,\\
\label{eq3bis} \sigma_3 g_1&=g_1 \sigma_3, \quad & \sigma_3
g_2&=g_2 \sigma_2, \quad &\sigma_3 g_3&=g_3 \sigma_1.
\end{align}
If $W:=\ku$ - span of $\{g_1,g_2,g_3\}$, then $W$ is a braided
vector subspace of $M(\oc_{s_7},\rho)$ of Cartan type, with matrix
of coefficients given by
\begin{align}\label{matrix:coef:Q1:Q2}
\mathcal Q_1:=\begin{pmatrix} -1 & -1 & 1 \\ 1 & -1 & -1
\\ -1 & 1 & -1
\end{pmatrix}, \qquad
\mathcal Q_2:=\begin{pmatrix} -1 & 1 & -1 \\ -1 & -1 & 1
\\ 1 & -1 & -1
\end{pmatrix},
\end{align}
for $\rho=\epsilon \otimes \nu_3$ and $\rho=\sgn \otimes \nu_3$,
respectively. In both cases the associated Cartan matrix is given
by \eqref{cartanmatrix:3x3}. By Theorem \ref{cartantype},
$\dim\toba(\oc_{s_7},\rho)=\infty$.

\smallbreak

\emph{CASE: $j=6$}. The element $s_6$ is real and we compute that
$M_{12}^{s_6}$ is a non-abelian group of order 32, whose character
table is given by Table \ref{tablacarM12^s6}.

\begin{table}[t]
\begin{center}
\footnotesize{\begin{tabular}{|c|c|c|c|c|c|c|c|c|c|c|c|c|c|c|c|c|}
\hline
{\bf $k$} &    1& 2& 3& 4& 5& 6& 7& 8& 9& 10& 11& 12& 13& 14\\
\hline \hline
{\bf $|y_k|$} &  1& 4& 4& 2& 4& 8& 4& 2& 2& 4& 4& 4& 8& 4 \\
\hline {\bf $|G^{y_k}|$} &  32& 8& 16& 32& 16&8& 16& 16& 8& 32& 32& 16& 8& 16\\
\hline {\bf $|\Oc_{y_k}|$} & 1& 4& 2& 1& 2& 4& 2& 2& 4& 1& 1& 2& 4& 2\\
\hline \hline

\hline {\bf $\mu_1$} & 1& 1& 1& 1& 1& 1& 1& 1& 1& 1& 1& 1& 1& 1\\

\hline {\bf $\mu_2$} & 1& -1& 1& 1&-1&1& -1& 1& -1& 1& 1& -1& 1& -1\\

\hline {\bf $\mu_3$} &   1& 1& 1& 1& -1& -1& -1& 1& 1& 1& 1& -1& -1& -1\\

\hline {\bf $\mu_4$} & 1& -1& 1& 1& 1& -1& 1& 1& -1& 1& 1& 1& -1& 1 \\

\hline {\bf $\mu_5$} & 1& 1& 1& 1& -i& -i& -i& -1& -1& -1& -1& i& i& i\\

\hline {\bf $\mu_6$} & 1& -1& 1& 1& -i& i& -i& -1& 1& -1& -1& i& -i& i\\

\hline {\bf $\mu_7$} & 1& 1& 1& 1& i& i& i& -1& -1& -1& -1&
-i& -i & -i\\

\hline {\bf $\mu_8$} & 1& -1& 1& 1& i& -i& i& -1& 1& -1& -1& -i& i& -i\\

\hline {\bf $\mu_9$} & 2& 0& -2& 2& 0& 0& 0& -2& 0& 2& 2& 0& 0& 0\\

\hline {\bf $\mu_{10}$} & 2& 0& -2& 2& 0& 0& 0& 2& 0& -2& -2& 0& 0& 0\\

\hline {\bf $\mu_{11}$} & 2& 0& 0& -2& 1+i& 0& -1-i& 0& 0&
-2i & 2i& 1-i& 0& -1+i\\

\hline {\bf $\mu_{12}$} &    2& 0& 0& -2& -1-i& 0& 1+i& 0& 0& -2i& 2i& -1+i& 0& 1-i\\

\hline {\bf $\mu_{13}$}&    2& 0& 0& -2& -1+i& 0& 1-i& 0& 0& 2i& -2i& -1-i& 0& 1+i\\

\hline {\bf $\mu_{14}$} &  2& 0& 0& -2& 1-i& 0& -1+i& 0& 0& 2i& -2i& 1+i& 0& -1-i\\

\hline
\end{tabular}}
\end{center}
\caption{Character table of $M_{12}^{s_6}$.}\label{tablacarM12^s6}
\end{table}

For every $k$, $1 \leq k \leq 14$, we call $\rho_k=(\rho_k,V_k)$
the irreducible representation of $M_{12}^{s_6}$ whose character
is $\mu_k$. From Table \ref{tablacarM12^s6}, we can conclude that
if $k\neq 5$, $6$, $7$, $8$, $10$, then $q_{s_6s_6}\neq -1$ and
$\dim\toba(\oc_{s_{6}}, \rho_k)=\infty$, by Lemma \ref{odd}. On
the other hand, if $k= 5$, $6$, $7$, $8$ or $10$, then
$q_{s_6s_6}= -1$. For these five cases we will prove that
$\dim\toba(\oc_{s_{6}}, \rho_k)=\infty$. First, we compute that
$\oc_{s_6}\cap M_{12}^{s_6}=\{\sigma_l\ \,| \, 1 \leq l \leq 6\}$,
where
\begin{align*}
\sigma_1&:=( 1, 4, 7,12)( 2, 6,11,10)( 3, 9) ( 5, 8),\\
\sigma_2&:=(1, 4, 7,12)( 2,11)(3, 5, 9, 8) ( 6,10),
\end{align*}
$\sigma_3:=s_6^{-1}$, $\sigma_4:=s_6$, $\sigma_5:=\sigma_1^{-1}$
and $\sigma_6:=\sigma_2^{-1}$. We compute that these elements
commute each other and that $\sigma_2\in \oc_{7}^{M_{12}^{s_6}}$
and $\sigma_5$, $\sigma_6\in \oc_{14}^{M_{12}^{s_6}}$. We take in
$M_{12}$ the following elements:
\begin{align*}
g_1&:=( 1, 8, 6,12, 3, 2)( 4, 9,11, 7, 5,10), & g_2&:= ( 1,
6,12,11, 7,10,4, 2)( 5, 8)
\end{align*}
and  $g_3:=( 1,12)( 4, 7)( 5, 8)( 6,10)$. Then $\sigma_r
g_r=g_rs_6$,  $1 \leq r \leq 3$, and
\begin{align*} \sigma_2g_1&=g_1\sigma_5 ,\quad &\sigma_1g_2&=g_2
\sigma_5, \quad &\sigma_1g_3&=g_3 \sigma_5,\\
\sigma_3g_1&=g_1\sigma_6 ,\quad &\sigma_3g_2&=g_2 \sigma_4, \quad
&\sigma_2g_3&=g_3 \sigma_6.
\end{align*}

Assume that $k= 5$, $6$, $7$ or $8$. We define $W:=\ku$ - span of
$\{g_1,g_2,g_3\}$. Hence, $W$ is a braided vector subspace of
$M(\oc_{s_6},\rho_k)$ of Cartan type. From Table
\ref{tablacarM12^s6}, we can calculate that the matrix of
coefficients $\mathcal Q$ is
\begin{align}\label{matrix:coef:Q3:Q4}
\mathcal Q_3:=\begin{pmatrix} -1 & i & i \\ i & -1 & i
\\ i & i & -1
\end{pmatrix}, \qquad
\mathcal Q_4:=\begin{pmatrix} -1 & -i & -i \\ -i & -1 & -i
\\ -i & -i & -1
\end{pmatrix},
\end{align}
for $k=5$ or $6$, and $k=7$ or $8$, respectively. In all the cases
the associated Cartan matrix is as in \eqref{cartanmatrix:3x3}.
Therefore, $\dim\toba(\oc_{s_{6}}, \rho_k)=\infty$.

Assume that $k=10$. Since $\sigma_5$ and $\sigma_6$ commute there
exists a basis $\{v_1, v_2\}$ of $V_{10}$, the vector space
affording $\rho_{10}$, composed by simultaneous eigenvectors of
$\rho_{10}(\sigma_5)$ and $\rho_{10}(\sigma_6)$. Let us say
$\rho_{10}(\sigma_5)v_l=\lambda_l v_l$ and
$\rho_{10}(\sigma_6)v_l=\kappa_l v_l$, $l=1$, $2$. Notice that
$\kappa_l,\lambda_l=\pm 1,\pm i$, due to
$|\sigma_5|=4=|\sigma_6|$. From Table \ref{tablacarM12^s6}, we can
deduce that $\lambda_1+\lambda_2=0=\kappa_1+\kappa_2$ because
$\sigma_5$, $\sigma_6\in \oc_{14}^{M_{12}^{s_6}}$  and
$\mu_{10}(\oc_{14}^{M_{12}^{s_6}})=0$. Also, since
$\sigma_5\sigma_6=s_6^{-1}$ we have that $\lambda_l\kappa_l=-1$,
$l=1$, $2$. Now, we consider the four possibilities:
$\lambda_1=\pm1$, $\pm i$.
\begin{itemize}
\item[(i)] If $\lambda_1=\pm 1$, then we take $W:=\ku$ - span of
$\{g_1v_1,g_2v_2,g_3v_1\}$.
\item[(ii)] If $\lambda_1=\pm i$, then we take $W:=\ku$ - span of
$\{g_1v_1,g_2v_1,g_3v_1\}$.
\end{itemize}
In both cases, $W$ is a braided vector subspace of
$M(\oc_{s_6},\rho_{10})$ of Cartan type. We calculate that the
matrices of coefficients are given by $\mathcal Q_1$ (resp.
$\mathcal Q_2$) for the case $\lambda_1=1$ (resp. $\lambda_1=-1$)
-- see \eqref{matrix:coef:Q1:Q2}; whereas the matrices of
coefficients are given by $\mathcal Q_3$ (resp. $\mathcal Q_4$)
for the case $\lambda_1=i$ (resp. $\lambda_1=-i$) -- see
\eqref{matrix:coef:Q3:Q4}. In all these cases, the associated
Cartan matrix is given by \eqref{cartanmatrix:3x3}. Therefore,
$\dim \toba(\oc_{s_6}, \rho_{10}) = \infty$.

\smallbreak

\emph{CASE: $j=15$}. We compute that $M_{12}^{s_{15}}\simeq
M_{12}^{s_6}$. This implies that this case is analogous to the
case $j=6$, since $\oc_{s_{15}}\simeq \oc_{s_{6}}$ as racks.

\smallbreak

\emph{CASE: $j=4$}. We compute that $M_{12}^{s_4}$ is a
non-abelian group of order 192, whose character table is given by
Table \ref{tablacarM12^s4}.

For every $k$, $1 \leq k \leq 13$, we call $\rho_k=(\rho_k,V_k)$
the irreducible representation of $M_{12}^{s_4}$ whose character
is $\mu_k$. From Table \ref{tablacarM12^s4}, we can conclude that
if $k\neq 10$, $11$, $13$, then $q_{s_4s_4}\neq -1$ and
$\dim\toba(\oc_{s_{4}}, \rho_k)=\infty$, by Lemma \ref{odd}.

\begin{table}[t]
\begin{center}
\small{\begin{tabular}{|c|c|c|c|c|c|c|c|c|c|c|c|c|c|c|c|} \hline
{\bf $k$} &    1& 2& 3& 4& 5& 6& 7& 8& 9& 10& 11& 12& 13\\
\hline \hline
{\bf $|y_k|$} & 1& 2& 3& 8& 6& 4& 2& 2& 4& 8& 4& 2& 4\\
\hline {\bf $|G^{y_k}|$} & 192& 8& 6& 8& 6& 32& 192& 32& 16& 8&
32& 16& 16
 \\
\hline {\bf $|\Oc_{y_k}|$} &  1& 24& 32& 24& 32& 6& 1& 6& 12& 24& 6& 12& 12 \\
\hline \hline

\hline {\bf $\mu_1$} & 1& 1& 1& 1& 1& 1& 1& 1&1& 1& 1& 1& 1\\

\hline {\bf $\mu_2$} & 1& -1& 1& -1& 1& 1& 1& 1& -1& -1& 1& 1& -1\\

\hline {\bf $\mu_3$} & 2& 0& -1& 0& -1& 2& 2& 2& 0& 0& 2& 2& 0 \\

\hline {\bf $\mu_4$} & 3& -1& 0& 1& 0& -1& 3& -1& 1& -1& 3& -1& 1\\

\hline {\bf $\mu_5$} & 3& 1& 0& -1& 0& -1& 3& -1& -1& 1& 3& -1& -1\\

\hline {\bf $\mu_6$} & 3& -1& 0& -1& 0& 3& 3& -1& 1& 1& -1& -1& 1\\

\hline {\bf $\mu_7$} & 3& 1& 0& -1& 0& -1& 3& 3& 1& -1& -1& -1& 1\\

\hline {\bf $\mu_8$} & 3& 1& 0& 1& 0& 3& 3& -1& -1& -1& -1& -1& -1\\

\hline {\bf $\mu_9$} & 3& -1& 0& 1& 0& -1& 3& 3& -1& 1& -1& -1& -1\\

\hline {\bf $\mu_{10}$} & 4& 0& 1& 0& -1& 0& -4& 0& 2& 0& 0& 0& -2\\

\hline {\bf $\mu_{11}$} & 4& 0& 1& 0& -1& 0& -4& 0& -2& 0& 0& 0& 2\\

\hline {\bf $\mu_{12}$} & 6& 0& 0& 0& 0& -2& 6& -2& 0& 0& -2& 2& 0\\

\hline {\bf $\mu_{13}$}&  8& 0& -1& 0& 1& 0& -8& 0& 0& 0& 0& 0& 0\\
\hline
\end{tabular}}
\end{center}
\caption{Character table of $M_{12}^{s_4}$.}\label{tablacarM12^s4}
\end{table}

Assume that $k= 10$, $11$ or $13$; thus $q_{s_4s_4}= -1$. We will
prove that $\dim\toba(\oc_{s_{4}}, \rho_k)=\infty$. First, we
compute that $\oc_{s_4}\cap M_{12}^{s_4}$ has 31 elements, and
that it contains $\sigma_1:=s_4$,
\begin{align*}
\sigma_2:=( 3, 4)( 5,10)( 6,11)( 7,12) \quad \text{and} \quad
\sigma_3:=( 2, 9)( 3, 5)( 4,10)( 7,12).
\end{align*}
We compute that $\sigma_2$ and $\sigma_3$ commute, $\sigma_3 \in
\oc_{2}^{M_{12}^{s_4}}$ and $\sigma_2\sigma_3=s_4$. Now, we choose
in $M_{12}$ the following elements: $g_1:=\id$,
\begin{align*}
g_2:=( 2, 7)( 3, 5)( 6,11)( 9,12) \quad \text{and} \quad g_3:=( 2,
9)( 5,10)( 6, 7)(11,12).
\end{align*}
It is easy to check that they satisfy the relations given in
\eqref{eq1}, \eqref{eq2} and \eqref{eq3}. From Table
\ref{tablacarM12^s4}, we have that
$\mu_k(\oc_{2}^{M_{12}^{s_4}})=0$. Now, we can proceed as in the
case $j=5$ and $k=8$ of the proof of Theorem \ref{teorM11}. Thus,
we can obtain a braided vector subspace of $M(\oc_{s_4},\rho_k)$
of Cartan type whose associated Cartan matrix is not of finite
type. Therefore, $\dim\toba(\oc_{s_{4}}, \rho_k)=\infty$, for
$k=10$, $11$, $13$.

\smallbreak

\emph{CASE: $j=9$}. We compute that $M_{12}^{s_9}$ is a
non-abelian group of order 240, whose character table is given by
Table \ref{tablacarM12^s9}.

\begin{table}[t]
\begin{center}
\small{\begin{tabular}{|c|c|c|c|c|c|c|c|c|c|c|c|c|c|c|c|c|} \hline
{\bf $k$} &    1& 2& 3& 4& 5& 6& 7& 8& 9& 10& 11& 12& 13& 14\\
\hline \hline
{\bf $|y_k|$} & 1& 2& 4& 5& 2& 4& 10& 2& 2& 6& 6& 3& 6& 2\\
\hline {\bf $|G^{y_k}|$} & 240& 16& 8& 10&16& 8& 10& 240& 24& 12& 12& 12& 12& 24\\
\hline {\bf $|\Oc_{y_k}|$} & 1& 15& 30& 24& 15& 30&24& 1& 10& 20& 20& 20& 20& 10\\
\hline \hline

\hline {\bf $\mu_1$} & 1& 1& 1& 1& 1& 1& 1& 1& 1& 1& 1& 1& 1& 1\\

\hline {\bf $\mu_2$} & 1& 1& -1& 1& -1& 1& -1& -1& 1& -1& 1& 1& -1& -1\\

\hline {\bf $\mu_3$} & 1& 1& 1& 1& -1& -1& -1& -1& -1& -1& -1& 1& 1& 1\\

\hline {\bf $\mu_4$} & 1& 1& -1& 1& 1& -1& 1& 1& -1& 1& -1& 1& -1& -1 \\

\hline {\bf $\mu_5$} & 4& 0& 0& -1& 0& 0& -1& 4& -2& 1& 1& 1& 1& -2 \\

\hline {\bf $\mu_6$} & 4& 0& 0& -1& 0& 0& 1& -4& 2& -1& -1& 1& 1& -2\\

\hline {\bf $\mu_7$} & 4& 0& 0& -1& 0& 0& -1& 4& 2& 1& -1& 1& -1& 2\\

\hline {\bf $\mu_8$} & 4& 0& 0& -1& 0& 0& 1& -4& -2& -1& 1& 1& -1& 2\\

\hline {\bf $\mu_9$} & 5& 1& -1& 0& 1& -1& 0& 5& 1& -1& 1& -1& 1& 1\\

\hline {\bf $\mu_{10}$} & 5& 1& -1& 0& -1& 1& 0& -5& -1& 1& -1& -1& 1& 1\\

\hline {\bf $\mu_{11}$} & 5& 1& 1& 0& 1& 1& 0& 5& -1& -1& -1& -1& -1& -1\\

\hline {\bf $\mu_{12}$} & 5& 1& 1& 0& -1& -1& 0& -5& 1& 1& 1& -1& -1& -1 \\

\hline {\bf $\mu_{13}$}&  6& -2& 0& 1& -2& 0& 1& 6& 0& 0& 0& 0& 0& 0 \\

\hline {\bf $\mu_{14}$} &  6& -2& 0& 1& 2& 0& -1& -6& 0& 0& 0& 0& 0& 0\\

\hline
\end{tabular}}
\end{center}
\caption{Character table of $M_{12}^{s_9}$.}\label{tablacarM12^s9}
\end{table}

For every $k$, $1 \leq k \leq 14$, we call $\rho_k=(\rho_k,V_k)$
the irreducible representation of $M_{12}^{s_9}$ whose character
is $\mu_k$. From Table \ref{tablacarM12^s9}, we can conclude that
if $k\neq 2$, $3$, $6$, $8$, $10$, $12$, $14$, then
$q_{s_9s_9}\neq -1$ and $\dim\toba(\oc_{s_{9}}, \rho_k)=\infty$,
by Lemma \ref{odd}. On the other hand, if $k= 2$, $3$, $6$, $8$,
$10$, $12$ or $14$, then $q_{s_9s_9}= -1$. For these cases we will
prove that $\dim\toba(\oc_{s_{9}}, \rho_k)=\infty$. First, we
compute that $\oc_{s_9}\cap M_{12}^{s_9}$ has 36 elements, and
that it contains $\sigma_1:=s_9$,
\begin{align*}
\sigma_2&:=( 1, 3)( 2, 5)( 4, 6)( 7, 9)( 8,10)(11,12),\\
\sigma_3&:=( 1, 5)( 2, 3)( 4,10)( 6, 8)( 7,11)( 9,12).
\end{align*}
We set $g_1:=\id$, $g_2:=( 2, 3, 4)( 5, 6, 8)( 7,11, 9)$, $g_3:=(
2, 5, 3)( 4, 8, 6)( 7,12, 9)$; these elements are in $M_{12}$ and
they satisfy $\sigma_r g_r =g_r s_9$, $1 \leq r \leq 3$, $\sigma_3
g_2 =g_2 \sigma_3$ and
\begin{align*} \sigma_2g_1&=g_1\sigma_2 ,\quad &\sigma_1g_2&=g_2
\gamma_{1,2}, \quad &\sigma_1g_3&=g_3 \sigma_2,\\
\sigma_3g_1&=g_1\sigma_3 ,\quad &\sigma_3g_2&=g_2 \gamma_{3,2},
\quad &\sigma_2g_3&=g_3 \sigma_3,
\end{align*}
where $\gamma_{1,2}=( 1, 4)( 2, 8)( 3, 6)( 5,10)( 7,11)( 9,12)$
and $\gamma_{3,2}=s_9 \gamma_{1,2}^{-1}$. Also, we compute that
$\sigma_2\sigma_3=s_9=\sigma_3\sigma_2$, and that $\sigma_3$,
$\gamma_{1,2} \in \oc_{14}^{M_{12}^{s_9}}$ and $\gamma_{3,2} \in
\oc_{9}^{M_{12}^{s_9}}$.

Assume that $k=2$ or $3$. Let us define $W:=\ku$ - span of
$\{g_1,g_2,g_3\}$. Hence, $W$ is a braided vector subspace of
$M(\oc_{s_9},\rho_k)$ of Cartan type. From Table
\ref{tablacarM12^s9}, it is straightforward to calculate that the
associated Cartan matrix is as in \eqref{cartanmatrix:3x3}. Thus,
$\dim\toba(\oc_{s_{9}}, \rho_k)=\infty$.

Assume that $k=6$. Since $\sigma_1$, $\sigma_2$, $\sigma_3$,
$\gamma_{1,2}$ and $\gamma_{3,2}$ commute there exists a basis
$\{v_1, v_2, v_3, v_4\}$ of $V_6$, the vector space affording
$\rho_6$, composed by simultaneous eigenvectors of
$\rho_6(\sigma_1)=-\Id$, $\rho_6(\sigma_2)$, $\rho_6(\sigma_3)$
$\rho_6(\gamma_{1,2})$ and $\rho_6(\gamma_{3,2})$. Let us call
$\rho_6(\sigma_2)v_l=\lambda_l v_l$, $\rho_6(\sigma_3)v_l=\kappa_l
v_l$,  $1\leq l \leq 4$, where $\lambda_l$, $\kappa_l= \pm 1$,
$1\leq l \leq 4$, due to $|\sigma_2|=|\sigma_3|=2$. From Table
\ref{tablacarM12^s9}, we have that
$\lambda_1+\lambda_2+\lambda_3+\lambda_4=-2=-\kappa_1-\kappa_2-\kappa_3-\kappa_4$. 
So, we have that $\lambda_1$, $\lambda_2$, $\lambda_3$ and
$\lambda_4$ are not all equal to $1$ or $-1$. On the other hand,
since $\sigma_2\sigma_3=s_9$ and $q_{s_9s_9}=-1$, we have that
$\lambda_l\kappa_l=-1$, $1\leq l\leq 4$. Now, if $W:=\ku$ - span
of $\{g_1v_l,g_2v_l,g_3v_l \, | \, 1\leq l\leq 4 \}$, then $W$ is
a braided vector subspace of $M(\oc_{s_9},\rho_6)$ of Cartan type,
whose associated Cartan matrix $\mathcal A$ has at least two row
with three $-1$ or more. This means that the corresponding Dynkin
diagram has at least two vertices with three edges or more; thus,
$\mathcal A$ is not of finite type. Hence, $\dim \toba(\oc_{s_9},
\rho_6) = \infty$.

For the cases $k=8$, $9$, $10$, $12$ or $14$, we proceed in an
analogous way. \epf

\begin{obs}
We compute that the groups $M_{12}^{s_3}$, $M_{12}^{s_8}$ and
$M_{12}^{s_{10}}$, have 10, 12 and 10 conjugacy classes,
respectively. Hence, there are 168 possible pairs $(\oc,\rho)$ for
$M_{12}$. We conclude that 164 of them lead to
infinite-dimensional Nichols algebras, and 4 have negative
braiding.
\end{obs}


\subsection{The group $M_{22}$}\label{sectionM22}

The Mathieu simple group $M_{22}$ can be given as a subgroup of
$\s_{22}$ in the following form
\begin{align*}
M_{22}:=\langle & ( 1, 2, 3, 4, 5, 6, 7, 8,9,10,11)(12,13,14,15,16,17,18,19,20,21,22), \,\\
& ( 1, 4, 5, 9, 3)( 2,8,10, 7, 6)(12,15,16,20,14)(13,19,21,18,17),  \\
&  ( 1,21) ( 2,10, 8,6)( 3,13, 4,17)( 5,19,
9,18)(11,22)(12,14,16,20) \, \rangle.
\end{align*}
In Table \ref{tablacarM22}, we show the character table of
$M_{22}$, with $A = (-1-i\sqrt{11})/2$ and $C=(-1-i\sqrt{7})/2$.
The representatives of the conjugacy classes of $M_{22}$ are
$s_1:=\id,$
\begin{align*}   s_2&:= ( 1,10,13,17)( 2,
3,14,15)( 4,20,18, 7)( 5,21)( 6,22)(9,11,12,16),\\
s_3&:= ( 1,13)( 2,14)( 3,15)( 4,18)( 7,20)( 9,12)(10,17)(11,16),\\
s_4&:=( 1, 8,17, 5,11,15, 3, 7)( 2,14, 9,16)( 4,20)( 6,21,13,22,19,18,12,10),\\
s_5&:=  ( 1,12,16,15,19,11,18)( 2, 7, 9,14,13,10, 6)( 3,22, 4,17, 5,21, 8),\\
s_6&:= ( 1,15,18,16,11,12,19)( 2,14, 6, 9,10, 7,13)( 3,17, 8, 4,21,22, 5), \\
s_7&:= ( 1, 4, 2, 6, 3)( 5,15,12,22,18)( 7, 8,11,19,20)( 9,17,10,14,21),\\
s_8&:= ( 1,18, 4,12,15, 8, 3,17,19, 7, 6)( 2, 9,16,11,13,22,20, 5,10,14,21),\\
s_9&:= ( 1, 4,15, 3,19, 6,18,12, 8,17, 7)( 2,16,13,20,10,21,
9,11,22, 5,14),\\
s_{10}&:=  ( 1, 6, 5,17)( 3, 8)( 4,11)( 7,13,16,14)( 9,12,22,15)(10,20,18,19),\\
s_{11}&:= ( 1, 7,22)( 2,13, 6,14, 5, 3)( 4,10)( 8,16, 9,20,19,17)(11,15,21)(12,18),\\
s_{12}&:=( 1,22, 7)( 2, 6, 5)( 3,13,14)( 8,
9,19)(11,21,15)(16,20,17).
\end{align*}

\begin{table}[t]
\begin{center}
{\small{\begin{tabular}{|c|c|c|c|c|c|c|c|c|c|c|c|c|}
\hline {\bf $j$} &    1& 2& 3& 4& 5& 6& 7& 8& 9& 10 & 11 & 12\\
\hline \hline {\bf $|s_j|$} &  1& 4& 2& 8& 7& 7& 5& 11& 11& 4& 6& 3 \\
\hline {\bf $|G^{s_j}|$} & 443520& 32& 384& 8& 7& 7& 5& 11& 11& 16& 12& 36  \\
\hline \hline

\hline {\bf $\chi_1$} &    1& 1& 1& 1& 1& 1& 1& 1& 1& 1& 1& 1\\

\hline {\bf $\chi_2$} &  21& 1& 5& -1& 0& 0& 1& -1& -1& 1& -1& 3 \\

\hline {\bf $\chi_3$} &  45& 1& -3& -1& C& C'& 0& 1& 1& 1& 0& 0 \\

\hline {\bf $\chi_4$} & 45& 1& -3& -1& C'& C& 0& 1& 1& 1& 0& 0  \\

\hline {\bf $\chi_5$} & 55& 3& 7& 1& -1& -1& 0& 0& 0& -1& 1& 1  \\

\hline {\bf $\chi_6$} &  99& 3& 3& -1& 1& 1& -1& 0& 0& -1& 0& 0  \\

\hline {\bf $\chi_7$} &  154& -2& 10& 0& 0& 0& -1& 0& 0& 2& 1& 1 \\

\hline {\bf $\chi_8$} &  210& -2& 2& 0& 0& 0& 0& 1& 1& -2& -1& 3 \\

\hline {\bf $\chi_9$} & 231& -1& 7& -1& 0& 0& 1& 0& 0& -1& 1& -3\\

\hline {\bf $\chi_{10}$} & 280& 0& -8& 0& 0& 0& 0& A'& A& 0& 1& 1\\

\hline {\bf $\chi_{11}$} & 280& 0& -8& 0& 0& 0& 0& A & A'& 0& 1& 1 \\

\hline {\bf $\chi_{12}$} & 385& 1& 1& 1& 0& 0& 0& 0& 0& 1& -2& -2\\

\hline
\end{tabular}}}
\end{center}
\caption{Character table of $M_{22}$.}\label{tablacarM22}
\end{table}

In the following statement, we summarize our study by mean of
abelian subracks in the group $M_{22}$.

\begin{theorem}\label{teorM22}
Let $\rho\in \widehat{M_{22}^{s_j}}$, with $1\leq j \leq 12$. If
$j=4$ and $\rho=\chi_{(-1)}$,  then the braiding is negative.
Otherwise, $\dim\toba(\oc_{s_j},\rho)=\infty$.
\end{theorem}

\pf

\emph{CASE: $j=7$, $12$}. From Table \ref{tablacarM22}, we see
that $s_j$ is real. By Lemma \ref{odd}, $\dim\toba(\oc_{s_j},
\rho)=\infty$, for all $\rho \in \widehat{M_{22}^{s_j}}$.

\smallbreak

 \emph{CASE: $j=5$, $6$}. We compute that $s_j^2$ and
$s_j^{4}$ are in $\oc_{s_j}$. Since $|s_j|=7$ we have that
$\dim\toba(\oc_{s_j}, \rho)=\infty$, for all $\rho \in
\widehat{M_{22}^{s_j}}$, by Lemma \ref{lemaB}.

\smallbreak

\emph{CASE: $j=8$, $9$}. We compute that $s_j^3$ and $s_j^{9}$ are
in $\oc_{s_j}$. Since $|s_j|=11$ we have that
$\dim\toba(\oc_{s_j}, \rho)=\infty$, for all $\rho \in
\widehat{M_{22}^{s_j}}$, by Lemma \ref{lemaB}.

\smallbreak

\emph{CASE: $j=4$}. We compute that $M_{22}^{s_4}=\langle s_4
\rangle\simeq \Z_8$. From Table \ref{tablacarM22}, we see that
$s_4$ is real. Thus, if $q_{s_js_j}\neq -1$, then
$\dim\toba(\oc_{s_j}, \rho)=\infty$, by Lemma \ref{odd}. The
remained case corresponds to $\rho(s_4)=\omega_8^4=-1$, which
satisfies $q_{s_js_j}= -1$. We compute that $\oc_{s_4}\cap
M_{22}^{s_4}=\{s_4,s_4^3,s_4^5,s_4^7\}$. It is straightforward to
prove that the braiding is negative.

\smallbreak

\emph{CASE: $j=11$}. We compute that $M_{22}^{s_{11}}=\langle
x,s_{11} \rangle\simeq \Z_2 \times \Z_6$, where
$$
x:=( 2, 9)( 3,16)( 4,12)( 5, 8)( 6,19)(10,18)(13,20)(14,17).
$$
Let us define $\{\nu_0,\dots,\nu_5\}$, where
$\nu_l(s_{11}):=\omega_6^l$, $0 \leq l \leq 5$. So,
\begin{align*}
\widehat{M_{22}^{s_{11}}}=\{\epsilon \otimes \nu_l,\, \sgn \otimes
\nu_l \,\, | \,\, 0 \leq l\leq 5 \},
\end{align*}
where $\epsilon$ and $\sgn$ mean the trivial and the sign
representations of $\Z_2$, respectively. Since $s_{11}$ is real we
have that if $\rho \in \widehat{M_{22}^{s_{11}}}$, with $l\neq 3$,
then $q_{s_{11}s_{11}}\neq -1$, and $\dim\toba(\oc_{s_{11}},
\rho)=\infty$, by Lemma \ref{odd}. The remained two cases are
$\rho=\epsilon \otimes \nu_3$ and $\rho=\sgn \otimes \nu_3$. We
will prove that also the Nichols algebra $\toba(\oc_{s_{11}},
\rho)$ is infinite-dimensional. First, we compute that
$\oc_{s_{11}}\cap M_{22}^{s_{11}}$ contains $\sigma_1:=s_{11}$,
\begin{align*}
\sigma_2&:=( 1, 7,22)( 2, 8, 6, 9, 5,19)( 3,17,13,16,14,20)(
4,12)(10,18)(11,15,21),\\ \sigma_3&:=  ( 1, 7,22)( 2,20, 6,17,
5,16)( 3, 9,13,19,14, 8)( 4,18)(10,12)(11,15,21).
\end{align*}
We compute that $\sigma_2=x s_{11}^4$ and $\sigma_3=x s_{11}$. We
choose $g_1:=\id$,
\begin{align*}
g_2&:=( 1, 7,22)( 3,19,16)( 4,12,10)( 8,20,13)( 9,17,14)(11,21,15)
\end{align*}
and $g_3:=g_2^{-1}$. These elements are in $M_{22}$ and they
satisfy the same relations as in \eqref{eq1bis}, \eqref{eq2bis}
and \eqref{eq3bis}. Now, if $W:=\ku$ - span of $\{g_1,g_2,g_3\}$,
then $W$ is a braided vector subspace of $M(\oc_{s_{11}},\rho)$ of
Cartan type with matrix of coefficients given by $\mathcal Q_1$
(resp. $\mathcal Q_2$) for the case $\rho=\epsilon \otimes \nu_3$
(resp. $\rho=\sgn \otimes \nu_3$) -- see
\eqref{matrix:coef:Q1:Q2}. In both cases the associated Cartan
matrix is as in \eqref{cartanmatrix:3x3}. By Theorem
\ref{cartantype}, $\dim\toba(\oc_{s_{11}},\rho)=\infty$.

\smallbreak

\emph{CASE: $j=10$}. We compute that $M_{22}^{s_{10}}=\langle
x,s_{10} \rangle\simeq \Z_4 \times \Z_4$, where
$$
x:=( 1, 9,13,10)( 3,11)( 4, 8)( 5,22,14,18)(
6,12,16,20)(7,19,17,15).
$$
If we set $\{\nu_0,\dots,\nu_3\}$, where $\nu_l(-):=\omega_4^l$,
$0 \leq l \leq 3$, then
\begin{align*}
\widehat{M_{22}^{s_{10}}}=\{ \nu_l\otimes \nu_t,\, \,\, | \,\, 0
\leq l, t \leq 3 \}.
\end{align*}
If $\rho=\nu_l\otimes \nu_t$, with $t\neq 2$ and $0\leq l \leq 3$,
then $q_{s_{10}s_{10}}=(\nu_l\otimes \nu_t)
(s_{10})=\omega_4^t\neq -1$; since $s_{10}$ is real we have that
$\dim\toba(\oc_{s_{10}},\rho)=\infty$. Assume that
$\rho=\nu_l\otimes \nu_2$, $0\leq l \leq 3$. Let us define
$\sigma_1:=s_{10}$, $\sigma_2:=s_{10}^{-1}$, $\sigma_3:=x$ and
$\sigma_4:=x^{-1}$. We compute that $x\in \oc_{s_{10}}$. Now, we
choose $g_1:=\id$,
\begin{align*}
g_2&:= ( 3, 8)( 4,11)( 6,17)( 7,16)( 9,18)(10,22)(12,20)(15,19),\\
g_3&:=( 3, 8, 4)( 5,13,14)( 6, 9,19)( 7,22,15)(10,12,17)(16,18,20)
\end{align*}
and $g_4:=g_3g_2$. These elements are in $M_{22}$ and satisfy that
$\sigma_r g_r= g_r s_{10}$, $\sigma_r g_1= g_1 \sigma_r$, $1 \leq
r \leq 4$, and
\begin{align*}
\sigma_1g_2&=g_2 \,\,s_{10}^{-1} ,\quad &\sigma_1g_3&=g_3 \,\,
x^{-1}s_{10}^{-1}, \quad &\sigma_1g_4&=g_4 \,\, xs_{10}^{-1}, \\
\sigma_3g_2&=g_2 \,\,x^{-1}s_{10}^2, \quad &\sigma_2g_3&=g_3
\,\,xs_{10}, \quad &\sigma_2g_4&=g_4 \,\,x^{-1}s_{10},\\
\sigma_4g_2&=g_2 \,\,xs_{10}^2, \quad &\sigma_4g_3&=g_3
\,\,s_{10}^{-1}, \quad &\sigma_3g_4&=g_4 \,\,s_{10}^{-1}.
\end{align*}
If we define $W:=\ku$ - span of $\{g_1,g_2,g_3,g_4\}$, then $W$ is
a braided vector subspace of $M(\oc_{s_{10}},\rho)$ of Cartan
type, whose associated Cartan matrix is given by
\begin{align}\label{cartanmatrix:4x4}
\mathcal A=\begin{pmatrix}
2 & 0 & -1 & -1 \\
0 & 2 & -1 & -1 \\
-1 & -1 & 2 & 0 \\
-1 & -1 & 0 & 2
\end{pmatrix}.
\end{align}
By Theorem \ref{cartantype},
$\dim\toba(\oc_{s_{10}},\rho)=\infty$.

\smallbreak

\emph{CASE: $j=2$}. We have that $s_2$ is a real element, it has
order 4 and we compute that $M_{22}^{s_2}$ is a non-abelian group
of order 32.

Let $\rho=(\rho,V)\in \widehat{M_{22}^{s_{2}}}$. We will prove
that the Nichols algebra $\toba(\oc_{s_{2}}, \rho)$ is
infinite-dimensional. If $q_{s_{2}s_{2}}\neq -1$, then the result
follows from Lemma \ref{odd}. Assume that $q_{s_{2}s_{2}}= -1$. We
compute that $\oc_{s_{2}}\cap M_{22}^{s_{2}}$ has 16 elements and
it contains $\sigma_1:=s_2$,
\begin{align*}
\sigma_2&:= ( 1, 7,15,11)( 2,12,10, 4)( 3,16,13,20)( 5, 6)(
9,17,18,14)(21,22),\\
\sigma_3&:=( 1, 9, 3, 4)( 2,
7,17,16)(5,22)(6,21)(10,11,14,20)(12,15,18,13).
\end{align*}
We compute that $\sigma_1$, $\sigma_2$ and $\sigma_3$ commute and
that $\sigma_2\sigma_3=s_2^{-1}$. We choose $g_1:=\id$, $g_2:=(
2,16,12)( 3,13,15)( 4,17,11)( 5, 6,21)( 7,9,10)(14,20,18)$ and
$g_3:=g_2^{-1}$. These elements belong to $M_{22}$ and they
satisfy the relations given by \eqref{eq1bis}, \eqref{eq2bis} and
\eqref{eq3bis}. Now, we define $W:=\ku$ - span of
$\{g_1v,g_2v,g_3v\}$, where $v \in V-0$. Hence, it is
straightforward to check that $W$ is a braided vector subspace of
$M(\oc_{s_{2}},\rho)$ of Cartan type whose associated Cartan
matrix is given by \eqref{cartanmatrix:3x3}. By Theorem
\ref{cartantype}, $\dim\toba(\oc_{s_{2}},\rho)=\infty$.

\smallbreak

\emph{CASE: $j=3$}. We compute that $M_{22}^{s_3}$ is a
non-abelian group of order 384, whose character table is given by
Table \ref{tablacarM22^s3}, where $d:=i\sqrt{3}$.

\begin{table}[t]
\begin{center}
\tiny{\begin{tabular}{|c|c|c|c|c|c|c|c|c|c|c|c|c|c|c|c|c|c|c|c|}
\hline
{\bf $k$} &    1& 2& 3& 4& 5& 6& 7& 8& 9& 10& 11& 12& 13& 14 & 15 & 16 & 17\\
\hline \hline
{\bf $|y_k|$} &  1& 2& 2& 3& 4& 2& 6& 4& 4& 4& 8& 6& 4& 4& 6& 2& 2 \\
\hline {\bf $|G^{y_k}|$} & 384& 16& 32& 12& 8& 64& 12& 16& 32& 16&
8& 12& 16& 16& 12& 48&
384\\
\hline {\bf $|\Oc_{y_k}|$} &  1& 24& 12& 32& 48& 6& 32& 24& 12& 24& 48& 32& 24& 24& 32& 8& 1\\
\hline \hline

\hline {\bf $\mu_1$} & 1& 1& 1& 1& 1& 1& 1& 1& 1& 1& 1& 1& 1& 1& 1& 1& 1\\

\hline {\bf $\mu_2$} & 1& -1& 1& 1& -1& 1& 1& -1& 1& -1& 1& -1& 1& 1& -1& -1& 1\\

\hline {\bf $\mu_3$} & 1& 1& 1& 1& 1& 1& 1& 1& 1& -1& -1& -1& -1& -1& -1& -1& 1\\

\hline {\bf $\mu_4$} & 1& -1& 1& 1& -1& 1& 1& -1& 1& 1& -1& 1& -1& -1& 1& 1& 1\\

\hline {\bf $\mu_5$} & 2& 0& 2& -1& 0& 2& -1& 0& 2& -2& 0& 1& 0& 0& 1& -2& 2\\

\hline {\bf $\mu_6$} & 2& 0& 2& -1& 0& 2& -1& 0& 2& 2& 0& -1& 0& 0& -1& 2& 2\\

\hline {\bf $\mu_7$} & 3& -1& -1& 0& 1& 3& 0& -1& -1& 1& -1& 0& 1& 1& 0& -3& 3\\

\hline {\bf $\mu_8$} &  3& 1& -1& 0& -1& 3& 0& 1& -1& 1& 1& 0& -1& -1& 0& -3& 3\\

\hline {\bf $\mu_9$} & 3& -1& -1& 0& 1& 3& 0& -1& -1& -1& 1& 0& -1& -1& 0& 3& 3\\

\hline {\bf $\mu_{10}$} & 3& 1& -1& 0& -1& 3& 0& 1& -1& -1& -1& 0& 1& 1& 0& 3& 3\\

\hline {\bf $\mu_{11}$} &  6& 0& -2& 0& 0& -2& 0& 0& 2& 0& 0& 0& 2& -2& 0& 0& 6\\

\hline {\bf $\mu_{12}$} & 6& 0& -2& 0& 0& -2& 0& 0& 2& 0& 0& 0& -2& 2& 0& 0& 6\\

\hline {\bf $\mu_{13}$}& 6& -2& 2& 0& 0& -2& 0& 2& -2& 0& 0& 0& 0& 0& 0& 0& 6\\

\hline {\bf $\mu_{14}$} & 6& 2& 2& 0& 0& -2& 0& -2& -2& 0& 0& 0& 0& 0& 0& 0& 6\\

\hline {\bf $\mu_{15}$} & 8& 0& 0& 2& 0& 0& -2& 0& 0& 0& 0& 0& 0& 0& 0& 0& -8\\

\hline {\bf $\mu_{16}$} & 8& 0& 0& -1& 0& 0& 1& 0& 0& 0& 0& d & 0& 0& -d & 0& -8\\

\hline {\bf $\mu_{17}$} & 8& 0& 0& -1& 0& 0& 1& 0& 0& 0& 0& -d& 0& 0& d& 0& -8\\

\hline
\end{tabular}}
\end{center}
\caption{Character table of $M_{22}^{s_3}$.}\label{tablacarM22^s3}
\end{table}

For every $k$, $1 \leq k \leq 17$, we call $\rho_k=(\rho_k,V_k)$
the irreducible representation of $M_{22}^{s_3}$ whose character
is $\mu_k$. From Table \ref{tablacarM22^s3} and the fact that
$s_3$ is real, we have that if $k\neq 15$, $16$, $17$, then
$q_{s_3s_3}\neq -1$ and $\dim\toba(\oc_{s_{3}}, \rho_k)=\infty$,
by Lemma \ref{odd}. On the other hand, if $k= 15$, $16$ or $17$,
then $q_{s_3s_3}= -1$. For these cases we will prove that
$\dim\toba(\oc_{s_{3}}, \rho_k)=\infty$. First, we compute that
$\oc_{s_3}\cap M_{22}^{s_3}$ has 51 elements, and that it contains
$\sigma_1=s_3$,
\begin{align*}
\sigma_2&:=( 3,15)( 4,20)( 6,22)( 7,18)( 8,19)(
9,16)(10,17)(11,12),\\ \sigma_3&:= ( 1,13)( 2,14)( 4, 7)( 6,22)(
8,19)( 9,11)(12,16)(18,20).
\end{align*}
We compute that $\sigma_1$, $\sigma_2$ and $\sigma_3$ commute each
other, $\sigma_2\sigma_3=s_3$ and that
$\sigma_2\in\oc_{3}^{M_{22}^{s_3}}$. We choose in $M_{22}$ the
following elements: $g_1:=\id$,
\begin{align*}
g_2&:=( 1, 6)( 2, 8)( 4,11)( 7, 9)(12,18)(13,22)(14,19)(16,20),\\
g_3&:=( 3, 6)( 4,16)( 7,11)( 8,17)( 9,20)(10,19)(12,18)(15,22).
\end{align*}
They satisfy the same relations as in \eqref{eq1}, \eqref{eq2} and
\eqref{eq3}.

Assume that $k=15$, $16$ or $17$. Since $\sigma_1$, $\sigma_2$ and
$\sigma_3$ commute there exists a basis $\{v_l \, | \, 1\leq l
\leq 8  \}$ of $V_{k}$, the vector space affording $\rho_{k}$,
composed by simultaneous eigenvectors of
$\rho_{k}(\sigma_1)=-\Id$, $\rho_{k}(\sigma_2)$ and
$\rho_{k}(\sigma_3)$. Let us say $\rho_{k}(\sigma_2)v_l=\lambda_l
v_l$ and $\rho_{k}(\sigma_3)v_l=\kappa_l v_l$, $1\leq l \leq 8$.
Notice that $\lambda_l$, $\kappa_l=\pm 1$, $1\leq l \leq 8$, due
to $|\sigma_2|=2=|\sigma_3|$. On the other hand, since
$\sigma_2\sigma_3=s_3$ we have that $\lambda_l\kappa_l=-1$, $1\leq
l \leq 8$. From Table \ref{tablacarM22^s3}, we can deduce that
$\sum_{l=1}^8\lambda_l=0$ because
$\mu_{k}(\oc_{3}^{M_{22}^{s_3}})=0$. Reordering the basis we can
suppose that $\lambda_1=1$ and $\lambda_2=-1$; thus, $\kappa_1=-1$
and $\kappa_2=1$. We define $W:=\ku$ - span of
$\{g_1v_1,g_2v_2,g_3v_2\}$. Hence, $W$ is a braided vector
subspace of $M(\oc_{s_3},\rho_{k})$ of Cartan type whose Cartan
matrix is as in \eqref{cartanmatrix:3x3}, and
$\dim\toba(\oc_{s_{3}}, \rho_{k})=\infty$.\epf

\begin{obs}
We compute that the groups $M_{22}^{s_{2}}$ and $M_{22}^{s_{12}}$
have 14 and 12 conjugacy classes, respectively. Hence, there are
132 possible pairs $(\oc,\rho)$ for $M_{22}$; only one of then has
negative braiding. The other pairs have infinite-dimensional
Nichols algebras.
\end{obs}


\subsection{The group $M_{23}$}\label{sectionM23}

The Mathieu simple group $M_{23}$ can be given as the subgroup of
$\s_{23}$ generated by $\alpha_1$ and $\alpha_2$, where
\begin{align*}
\alpha_1&:= ( 1, 2, 3, 4, 5, 6, 7, 8,9,10,11,12,13,14,15,16,17,18,19,20,21,22,23),\\
\alpha_2&:=( 3,17,10, 7, 9)( 4,13,14,19,
5)(8,18,11,12,23)(15,20,22,21,16).
\end{align*}
The order of $M_{23}$ is $10200960$. In Table \ref{tablacarM23},
we show the character table of $M_{23}$, where $A = (-1+i
\sqrt{7}))/2$, $B =(-1+i\sqrt{11})/2$, $C =(-1+i \sqrt{15})/2$ and
$D=(-1+i \sqrt{23})/2$. We will denote the representatives of the
conjugacy classes of $M_{23}$ by $s_j$, $1\leq j \leq 17$.

\begin{table}[t]
\begin{center}
\tiny{\begin{tabular}{|p{0.53cm}|c|c|c|c|c|c|c|c|c|c|c|c|c|c|c|c|c|c|c|}
\hline
{\bf $j$} &    1& 2& 3& 4& 5& 6& 7& 8& 9& 10& 11& 12& 13& 14 & 15 & 16 & 17\\
\hline \hline
{\bf $|s_j|$} & 1& 2& 3& 4& 5& 6& 7& 7& 8& 11& 11& 14& 14& 15& 15& 23& 23\\
\hline {\bf $|G^{s_j}|$} & $|M_{23}|$& 2688 & 180& 32& 15& 12& 14&
14& 8& 11& 11& 14& 14& 15&
15& 23& 23\\
\hline \hline

\hline {\bf $\chi_1$} & 1& 1& 1& 1& 1& 1& 1& 1& 1& 1& 1& 1& 1& 1& 1& 1& 1\\

\hline {\bf $\chi_2$} & 22& 6& 4& 2& 2& 0& 1& 1& 0& 0& 0& -1& -1& -1& -1& -1& -1\\

\hline {\bf $\chi_3$} & 45& -3& 0& 1& 0& 0& A& A'& -1& 1& 1& -A& -A'& 0& 0& -1& -1\\

\hline {\bf $\chi_4$} &45& -3& 0& 1& 0& 0& A'& A& -1& 1& 1&
      -A'& -A& 0& 0& -1& -1 \\

\hline {\bf $\chi_5$} & 230& 22& 5& 2& 0& 1& -1& -1& 0& -1& -1& 1& 1& 0& 0& 0& 0\\

\hline {\bf $\chi_6$} &  231& 7& 6& -1& 1& -2& 0& 0& -1& 0& 0& 0& 0& 1& 1& 1& 1\\

\hline {\bf $\chi_7$} & 231& 7& -3& -1& 1& 1& 0& 0& -1& 0& 0& 0&
0&  C& C'& 1&  1 \\

\hline {\bf $\chi_8$} & 231& 7& -3& -1& 1& 1& 0& 0& -1& 0& 0& 0&
0& C' & C& 1& 1\\

\hline {\bf $\chi_9$} & 253& 13& 1& 1& -2& 1& 1& 1& -1& 0& 0& -1& -1& 1& 1& 0& 0\\

\hline {\bf $\chi_{10}$} & 770& -14& 5& -2& 0& 1& 0& 0& 0& 0& 0&
0& 0& 0& 0& D& D'\\

\hline {\bf $\chi_{11}$} & 770& -14& 5& -2& 0& 1& 0& 0& 0& 0& 0&
0& 0& 0& 0& D'& D\\

\hline {\bf $\chi_{12}$} & 896& 0& -4& 0& 1& 0& 0& 0& 0& B&
      B'& 0& 0& 1& 1& -1& -1\\

\hline {\bf $\chi_{13}$}& 896& 0& -4& 0& 1& 0& 0& 0& 0& B'&
      B& 0& 0& 1& 1& -1& -1 \\

\hline {\bf $\chi_{14}$} & 990& -18& 0& 2& 0& 0& A& A'& 0& 0& 0&
      A& A'& 0& 0& 1& 1\\

\hline {\bf $\chi_{15}$} & 990& -18& 0& 2& 0& 0& A'& A& 0& 0& 0&
      A'& A& 0& 0& 1& 1\\

\hline {\bf $\chi_{16}$} & 1035& 27& 0& -1& 0& 0& -1& -1& 1& 1& 1& -1& -1& 0& 0& 0& 0 \\

\hline {\bf $\chi_{17}$} & 2024& 8& -1& 0& -1& -1& 1& 1& 0& 0& 0& 1& 1& -1& -1& 0& 0\\

\hline
\end{tabular}}
\end{center}
\caption{Character table of $M_{23}$.}\label{tablacarM23}
\end{table}

In the following statement, we summarize our study by mean of
abelian subracks in the group $M_{23}$.

\begin{theorem}\label{teorM23}
Let $\rho\in \widehat{M_{23}^{s_j}}$, with $1\leq j \leq 17$. The
braiding is negative in the cases $j=9$, $12$ and $13$, with
$\rho=\chi_{(-1)}$. Otherwise, $\dim\toba(\oc_{s_j},\rho)=\infty$.
\end{theorem}

\pf

\emph{CASE: $j=3$, $5$}. From Table \ref{tablacarM23}, we see
that $s_j$ is real. By Lemma \ref{odd}, $\dim\toba(\oc_{s_j},
\rho)=\infty$, for all $\rho \in \widehat{M_{23}^{s_j}}$.

\smallbreak

 \emph{CASE: $j=7$, $8$, $14$, $15$, $16$, $17$ }. We compute that $s_j^2$ and
$s_j^{4}$ are in $\oc_{s_j}$. Since $|s_j|$ is odd we have that
$\dim\toba(\oc_{s_j}, \rho)=\infty$, for all $\rho \in
\widehat{M_{23}^{s_j}}$, by Lemma \ref{lemaB}.

\smallbreak

\emph{CASE: $j=10$, $11$}. We compute that $s_j^3$ and $s_j^{9}$
are in $\oc_{s_j}$. Since $|s_j|=11$ we have that
$\dim\toba(\oc_{s_j}, \rho)=\infty$, for all $\rho \in
\widehat{M_{23}^{s_j}}$, by Lemma \ref{lemaB}.

\smallbreak

\emph{CASE: $j=12$, $13$}. We compute that $M_{23}^{s_j}=\langle
s_j \rangle\simeq \Z_{14}$, and that $s_j^9$ and $s_j^{11}$ are in
$\oc_{s_j}$. Thus, if $q_{s_js_j}\neq -1$, then
$\dim\toba(\oc_{s_j}, \rho)=\infty$, by Lemma \ref{lemaB}. The
remained case corresponds to $\rho(s_j)=\omega_{14}^7=-1$, which
satisfies $q_{s_js_j}= -1$. We compute that $\oc_{s_j}\cap
M_{23}^{s_j}=\{s_j,s_j^9,s_j^{11}\}$. It is straightforward to
prove that the braiding is negative.

\smallbreak

\emph{CASE: $j=9$}. We compute that $M_{23}^{s_9}=\langle s_9
\rangle\simeq \Z_{8}$, and that $s_9$ is real. Thus, if
$q_{s_9s_9}\neq -1$, then $\dim\toba(\oc_{s_9}, \rho)=\infty$, by
Lemma \ref{odd}. The remained case corresponds to
$\rho(s_9)=\omega_{8}^4=-1$, which satisfies $q_{s_9s_9}= -1$. We
compute that $\oc_{s_9}\cap M_{23}^{s_9}= \{ s_9, s_9^3, s_9^5,
s_9^7 \}$. It is easy to check that the braiding is negative.

\smallbreak

\emph{CASE: $j=6$}. The representative is
$$s_6=( 1,19,20)( 2, 9,18,17,14, 5)( 3,21)( 4,13,23,10,11,22)( 6, 8,15)(
7,16),$$ it is real and has order 6. We compute that
$M_{23}^{s_6}=\langle x, s_6 \rangle\simeq \Z_2 \times \Z_{6}$,
with $x:= ( 2,22)( 3, 7)( 4, 9)(
5,11)(10,14)(13,18)(16,21)(17,23)$. Let us define
$\{\nu_0,\dots,\nu_5\}$, where $\nu_l(s_6):=\omega_6^l$, $0 \leq l
\leq 5$. So,
\begin{align*}
\widehat{M_{23}^{s_6}}=\{\epsilon \otimes \nu_l,\, \sgn \otimes
\nu_l \,\, | \,\, 0 \leq l\leq 5 \},
\end{align*}
where $\epsilon$ and $\sgn$ mean the trivial and the sign
representations of $\Z_2$, respectively. If $\rho \in
\widehat{M_{23}^{s_6}}$, with $l\neq 3$, then $q_{s_{6}s_{6}}\neq
-1$, and $\dim\toba(\oc_{s_{6}}, \rho)=\infty$, by Lemma
\ref{odd}. The remained two cases are $\rho=\epsilon \otimes
\nu_3$ and $\rho=\sgn \otimes \nu_3$. We will prove that also the
Nichols algebra $\toba(\oc_{s_{6}}, \rho)$ is
infinite-dimensional. First, we compute that $\oc_{s_6}\cap
M_{23}^{s_6}$ has 6 elements, and it contains $\sigma_1:=s_{6}$,
\begin{align*}
\sigma_2&:=( 1,19,20)( 2, 4,18,23,14,11)( 3,16)( 5,22,
9,13,17,10)(6, 8,15)( 7,21),\\
\sigma_3&:=( 1,19,20)( 2,10,18,22,14,13)( 3, 7)( 4, 5,23,
9,11,17)( 6, 8,15)(16,21).
\end{align*}
Also, we compute that $\sigma_2=xs_{6}$ and $\sigma_3=xs_{6}^4$.
We choose $g_1:=\id$,
\begin{align*}
g_2&:=( 1,20,19)( 3,21, 7)( 4,10, 9)( 5,11,13)( 6,
8,15)(17,23,22),
\end{align*}
and $g_3:=g_2^{-1}$. These elements are in $M_{23}$ and they
satisfy the same relations as in \eqref{eq1bis}, \eqref{eq2bis}
and \eqref{eq3bis}. If $W:=\ku$ - span of $\{g_1,g_2,g_3\}$, then
$W$ is a braided vector subspace of $M(\oc_{s_{6}},\rho)$ of
Cartan type whose Cartan matrix is as in \eqref{cartanmatrix:3x3}.
Therefore, $\dim\toba(\oc_{s_{6}},\rho)=\infty$.

\smallbreak

\emph{CASE: $j=4$}. The representative is
$$s_4=( 1,17,10, 4)( 2, 8)( 3, 6,14,11)( 5,12,13,21)( 7,15)(16,20,22,23),$$
which has order 4 and it is real. We compute that $M_{23}^{s_4}$
is a non-abelian group of order 32. Let $\rho=(\rho,V)\in
\widehat{M_{23}^{s_{4}}}$. We will prove that the Nichols algebra
$\toba(\oc_{s_{4}}, \rho)$ is infinite-dimensional. If
$q_{s_{4}s_{4}}\neq -1$, then the result follows from Lemma
\ref{odd}. Assume that $q_{s_{4}s_{4}}= -1$. We compute that
$\oc_{s_{4}}\cap M_{23}^{s_{4}}$ has 16 elements and it contains
$\sigma_1:=s_4$,
\begin{align*}
\sigma_2&:=( 1,12, 6,20)( 2, 7)( 3,16, 4, 5)(
8,15)(10,21,11,23)(13,14,22,17),\\
\sigma_3&:=( 1,16,11,13)( 2,15)( 3,21,17,20)( 4,23,14,12)(
5,10,22, 6)( 7, 8).
\end{align*}
These elements commute and $\sigma_2\sigma_3=s_4^{-1}$. We choose
$g_1:=\id$,
\begin{align*}
g_2&:=( 2,15, 7)( 3,21, 5)( 4,20,13)( 6,11,10)(12,16,17)(14,23,22)
\end{align*}
and $g_3:=g_2^{-1}$. These elements belong to $M_{23}$ and they
satisfy the relations given by \eqref{eq1bis}, \eqref{eq2bis} and
\eqref{eq3bis}. Now, we define $W:=\ku$ - span of
$\{g_1v,g_2v,g_3v\}$, where $v \in V-0$. Hence, it is
straightforward to check that $W$ is a braided vector subspace of
$M(\oc_{s_{4}},\rho)$ of Cartan type whose associated Cartan
matrix is given by \eqref{cartanmatrix:3x3}. By Theorem
\ref{cartantype}, $\dim\toba(\oc_{s_{4}},\rho)=\infty$.

\smallbreak

\emph{CASE: $j=2$}. The representative is
$$s_2=( 1,10)( 3,14)( 4,17)( 5,13)( 6,11)(12,21)(16,22)(20,23).$$
We compute that $M_{23}^{s_2}$ is a non-abelian group of order
2688, whose character table is given by Table
\ref{tablacarM23^s2}, where $A=(-1+i \sqrt{7})/2$ and $B= i
\sqrt{3}$.

{\tiny{
\begin{table}[t]
\begin{center}
\tiny{\begin{tabular}{|c|c|c|c|c|c|c|c|c|c|c|c|c|c|c|c|c|c|c|}
\hline
{\bf $k$} &    1& 2& 3& 4& 5& 6& 7& 8& 9& 10& 11& 12& 13& 14 & 15 & 16 \\
\hline \hline
{\bf $|y_k|$} &  1& 7& 7& 14& 14& 3& 6& 6& 6& 2& 4& 8& 2& 4& 2& 4\\
\hline {\bf $|G^{y_k}|$} & 2688& 14& 14& 14& 14& 12& 12& 12& 12&
2688& 32& 8& 192& 16& 32&
8\\
\hline \hline

\hline {\bf $\mu_1$} & 1& 1& 1& 1& 1& 1& 1& 1& 1& 1& 1& 1& 1& 1& 1& 1\\

\hline {\bf $\mu_2$} & 3& A& A'& A& A'& 0& 0& 0& 0& 3& -1& 1& 3& -1& -1& 1\\

\hline {\bf $\mu_3$} & 3& A'& A& A'& A& 0& 0& 0& 0& 3& -1& 1& 3& -1& -1& 1\\

\hline {\bf $\mu_4$} & 6& -1& -1& -1& -1& 0& 0& 0& 0& 6& 2& 0& 6& 2& 2& 0\\

\hline {\bf $\mu_5$} & 7& 0& 0& 0& 0& 1& 1& -1& -1& 7& 3& 1& -1& -1& -1& -1\\

\hline {\bf $\mu_6$} & 7& 0& 0& 0& 0& 1& 1& 1& 1& 7& -1& -1& 7& -1& -1& -1\\

\hline {\bf $\mu_7$} & 7& 0& 0& 0& 0& 1& 1& -1& -1& 7& -1& -1& -1& -1& 3& 1\\

\hline {\bf $\mu_8$} & 8& 1& 1& 1& 1& -1& -1& -1& -1& 8& 0& 0& 8& 0& 0& 0\\

\hline {\bf $\mu_9$} & 8& 1& 1& -1& -1& 2& -2& 0& 0& -8& 0& 0& 0& 0& 0& 0\\

\hline {\bf $\mu_{10}$} &  8& 1& 1& -1& -1& -1& 1& B& -B& -8& 0& 0& 0& 0& 0& 0\\

\hline {\bf $\mu_{11}$} & 8& 1& 1& -1& -1& -1& 1& -B& B& -8& 0& 0& 0& 0& 0& 0\\

\hline {\bf $\mu_{12}$} & 14& 0& 0& 0& 0& -1& -1& 1& 1& 14& 2& 0& -2& -2& 2& 0\\

\hline {\bf $\mu_{13}$}& 21& 0& 0& 0& 0& 0& 0& 0& 0& 21& 1& -1& -3& 1& -3& 1\\

\hline {\bf $\mu_{14}$} & 21& 0& 0& 0& 0& 0& 0& 0& 0& 21& -3& 1& -3& 1& 1& -1\\

\hline {\bf $\mu_{15}$} & 24& A'& A& -A'& -A& 0& 0& 0& 0& -24& 0& 0& 0& 0& 0& 0\\

\hline {\bf $\mu_{16}$} & 24& A& A'& -A& -A'& 0& 0& 0& 0& -24& 0& 0& 0& 0& 0& 0\\

\hline
\end{tabular}}
\end{center}
\caption{Character table of $M_{23}^{s_2}$.}\label{tablacarM23^s2}
\end{table}}}

For every $k$, $1 \leq k \leq 16$, we call $\rho_k=(\rho_k,V_k)$
the irreducible representation of $M_{23}^{s_2}$ whose character
is $\mu_k$. From Table \ref{tablacarM23^s2}, we have that if
$k\neq 9$, $10$, $11$, $15$, $16$, then $q_{s_2s_2}\neq -1$ and
$\dim\toba(\oc_{s_{2}}, \rho_k)=\infty$, by Lemma \ref{odd}. On
the other hand, if $k= 9$, $10$, $11$, $15$ or $16$, then
$q_{s_2s_2}= -1$. For these cases we will prove that
$\dim\toba(\oc_{s_{2}}, \rho_k)=\infty$. First, we compute that
$\oc_{s_2}\cap M_{23}^{s_2}$ has 99 elements and it contains
$\sigma_1:=s_2$ and
$$\sigma_2:=( 3, 6)( 5,20)( 7, 9)(11,14)(12,21)(13,23)(15,18)(16,22).$$
We compute that $\sigma_2 \in \oc_{15}^{M_{23}^{s_2}}$. Now, we
choose $g_1:=\id$ and
$$g_2:=( 1, 7)( 3,20)( 4,18)( 5,14)( 6,23)( 9,10)(11,13)(15,17).$$
Then, $g_2$ is in $M_{23}$, and we have $\sigma_r g_r=g_r
\sigma_1$, $r=1$, $2$, $\sigma_2 g_1=g_1 \sigma_2$ and $\sigma_1
g_2=g_2 \sigma_2$.

Assume that $k= 9$, $10$, $11$, $15$ or $16$. From Table
\ref{tablacarM23^s2}, we have that the degree of $\rho_k$ is $8$
or $24$. Since $\sigma_1$ and $\sigma_2$ commute there exists a
basis $\{v_l \, | \, 1\leq l \leq \deg(\rho_k)\}$ of $V_k$, the
vector space affording $\rho_k$, composed by simultaneous
eigenvectors of $\rho_k(\sigma_1)=-\Id$ and $\rho_k(\sigma_2)$.
Let us call $\rho_k(\sigma_2)v_l=\lambda_l v_l$, $1\leq l \leq
\deg(\rho_k)$, where $\lambda_l=\pm 1$, due to $|\sigma_2|=2$.
From Table \ref{tablacarM23^s2}, we have that
$\sum_{l=1}^{\deg(\rho_k)}\lambda_l=0$. Reordering the basis we
can suppose that
$\lambda_1=\cdots=\lambda_{\deg(\rho_k)/2}=1=-\lambda_{1+\deg(\rho_k)/2}=\cdots=-\lambda_{\deg(\rho_k)}$.
It is straightforward to check that if $W:=\ku$ - span of
$\{g_1v_l, g_2v_l \, | \, 1\leq l \leq \deg(\rho_k)\}$, then $W$
is a braided vector subspace of $M(\oc_{s_{2}},\rho)$ of Cartan
type whose associated Cartan matrix $\mathcal A$ has at least two
row with three $-1$ or more. This means that the corresponding
Dynkin diagram has at least two vertices with three edges or more;
thus, $\mathcal A$ is not of finite type. Hence, $\dim
\toba(\oc_{s_2}, \rho_k) = \infty$.\epf

\begin{obs}
We compute that the groups $M_{23}^{s_3}$, $M_{23}^{s_4}$,
$M_{23}^{s_5}$, $M_{23}^{s_7}$ and $M_{23}^{s_{8}}$ have 15, 14,
15, 14 and 14 conjugacy classes, respectively. Hence, there are
251 possible pairs $(\oc,\rho)$ for $M_{23}$; 248 of them lead to
infinite-dimensional Nichols algebras, and 3 have negative
braiding.
\end{obs}


\subsection{The group $M_{24}$}\label{sectionM24}

The Mathieu simple group $M_{24}$ can be given as the subgroup of
$\s_{24}$ generated by $\alpha_1$, $\alpha_2$ and $\alpha_3$,
where
\begin{align*}
\alpha_1&:= ( 1, 2, 3, 4, 5, 6, 7, 8,9,10,11,12,13,14,15,16,17,18,19,20,21,22,23),\\
\alpha_2&:=( 3,17,10, 7, 9)( 4,13,14,19,
5)(8,18,11,12,23)(15,20,22,21,16),\\
\alpha_3&:=( 1,24) ( 2,23)( 3,12)( 4,16)( 5,18)( 6,10)( 7,20)(
8,14)( 9,21)(11,17)\\&  \qquad (13,22) (15,19).
\end{align*}
The order of $M_{24}$ is $244823040$. In Tables \ref{tablacarM24i}
and \ref{tablacarM24ii}, we show the character table of $M_{24}$,
where $A = (-1+i\sqrt{7})/2$, $C=(-1+i\sqrt{15})/2$ and $D=(-1+i
\sqrt{23})/2$. We will denote the representatives of the conjugacy
classes of $M_{24}$ by $s_j$, $1\leq j \leq 26$.

\begin{table}[t]
\begin{center}
\tiny{\begin{tabular}{|c|c|c|c|c|c|c|c|c|c|c|c|c|c|} \hline
{\bf $j$} &    1& 2& 3& 4& 5& 6& 7& 8& 9& 10& 11& 12& 13 \\
\hline \hline
{\bf $|s_j|$} & 1& 2& 2& 3& 3& 4& 4& 4& 5& 6& 6& 7& 7 \\
\hline {\bf $|G^{s_j}|$} & $|M_{24}|$& 21504& 7680& 1080& 504&
384& 128& 96& 60& 24& 24& 42& 42 \\
\hline \hline

\hline {\bf $\chi_1$} & 1& 1& 1& 1& 1& 1& 1& 1& 1& 1& 1& 1& 1\\

\hline {\bf $\chi_2$} &  23& 7& -1& 5& -1& -1& 3& -1& 3& 1& -1& 2& 2\\

\hline {\bf $\chi_3$} &   45& -3& 5& 0& 3& -3& 1& 1& 0& 0 & -1& A& A' \\

\hline {\bf $\chi_4$} & 45& -3& 5& 0& 3& -3& 1& 1& 0& 0& -1& A'& A\\

\hline {\bf $\chi_5$} & 231& 7& -9& -3& 0& -1& -1& 3& 1& 1& 0& 0& 0\\

\hline {\bf $\chi_6$} &  231& 7& -9& -3& 0& -1& -1& 3& 1& 1& 0& 0& 0\\

\hline {\bf $\chi_7$} & 252& 28& 12& 9& 0& 4& 4& 0& 2& 1& 0& 0& 0 \\

\hline {\bf $\chi_8$} & 253& 13& -11& 10& 1& -3& 1& 1& 3& -2& 1& 1& 1\\

\hline {\bf $\chi_9$} & 483& 35& 3& 6& 0& 3& 3& 3& -2& 2& 0& 0& 0\\

\hline {\bf $\chi_{10}$} & 770& -14& 10& 5& -7& 2& -2& -2& 0& 1& 1& 0& 0\\

\hline {\bf $\chi_{11}$} & 770& -14& 10& 5& -7& 2& -2& -2& 0& 1& 1& 0& 0\\

\hline {\bf $\chi_{12}$} & 990& -18& -10& 0& 3& 6& 2& -2& 0& 0&
-1& A& A'\\

\hline {\bf $\chi_{13}$}&  990& -18& -10& 0& 3& 6& 2& -2& 0& 0&
-1& A'& A\\

\hline {\bf $\chi_{14}$} & 1035& 27& 35& 0& 6& 3& -1& 3& 0& 0& 2& -1& -1\\

\hline {\bf $\chi_{15}$} & 1035& -21& -5& 0& -3& 3& 3& -1& 0& 0& 1& 2 A&  2 A'\\

\hline {\bf $\chi_{16}$} & 1035& -21& -5& 0& -3& 3& 3& -1& 0& 0&
1& 2A'& 2A \\

\hline {\bf $\chi_{17}$} & 1265& 49& -15& 5& 8& -7& 1& -3& 0& 1& 0& -2& -2\\

\hline {\bf $\chi_{18}$} & 1771& -21& 11& 16& 7& 3& -5& -1& 1& 0& -1& 0& 0\\

\hline {\bf $\chi_{19}$} & 2024& 8& 24& -1& 8& 8& 0& 0& -1& -1& 0& 1& 1\\

\hline {\bf $\chi_{20}$} & 2277& 21& -19& 0& 6& -3& 1& -3& -3& 0& 2& 2& 2\\

\hline {\bf $\chi_{21}$}&    3312& 48& 16& 0& -6& 0& 0& 0& -3& 0& -2& 1& 1 \\

\hline {\bf $\chi_{22}$} & 3520& 64& 0& 10& -8& 0& 0& 0& 0& -2& 0& -1& -1\\

\hline {\bf $\chi_{23}$} & 5313& 49& 9& -15& 0& 1& -3& -3& 3& 1& 0& 0& 0\\

\hline {\bf $\chi_{24}$} & 5544& -56& 24& 9& 0& -8& 0& 0& -1& 1& 0& 0& 0\\

\hline {\bf $\chi_{25}$} & 5796& -28& 36& -9& 0& -4& 4& 0& 1& -1& 0& 0& 0\\

\hline {\bf $\chi_{26}$} & 10395& -21& -45& 0& 0& 3& -1& 3& 0& 0& 0& 0& 0\\

\hline
\end{tabular}}
\end{center}
\caption{Character table of $M_{24}$ (i).}\label{tablacarM24i}
\end{table}

\begin{table}[t]
\begin{center}
\tiny{\begin{tabular}{|c|c|c|c|c|c|c|c|c|c|c|c|c|c|} \hline
{\bf $j$} &    14& 15& 16& 17& 18& 19& 20& 21& 22& 23& 24& 25& 26 \\
\hline \hline
{\bf $|s_j|$} & 8& 10& 11& 12& 12& 14& 14& 15& 15&  21& 21& 23& 23\\
\hline {\bf $|G^{s_j}|$} & 16&  20& 11& 12& 12& 14& 14& 15& 15& 21& 21& 23& 23\\
\hline \hline

\hline {\bf $\chi_1$} &  1& 1& 1& 1& 1& 1& 1& 1& 1& 1& 1& 1& 1\\

\hline {\bf $\chi_2$} & 1& -1& 1& -1& -1& 0& 0& 0& 0& -1& -1& 0& 0\\

\hline {\bf $\chi_3$} & -1& 0& 1& 0& 1& -A& -A'& 0& 0& A& A'& -1& -1\\

\hline {\bf $\chi_4$} & -1& 0& 1& 0& 1& -A'& -A& 0& 0& A'& A& -1& -1\\

\hline {\bf $\chi_5$} &  -1& 1& 0& -1& 0& 0& 0& C& C'& 0& 0& 1& 1\\

\hline {\bf $\chi_6$} & -1& 1& 0& -1& 0& 0& 0& C'& C& 0& 0& 1& 1\\

\hline {\bf $\chi_7$} & 0& 2& -1& 1& 0& 0& 0& -1& -1& 0& 0& -1& -1 \\

\hline {\bf $\chi_8$} & -1& -1& 0& 0& 1& -1& -1& 0& 0& 1& 1& 0& 0\\

\hline {\bf $\chi_9$} &  -1& -2& -1& 0& 0& 0& 0& 1& 1& 0& 0& 0& 0\\

\hline {\bf $\chi_{10}$} & 0& 0& 0& -1& 1& 0& 0& 0& 0& 0& 0& D& D'\\

\hline {\bf $\chi_{11}$} & 0& 0& 0& -1& 1& 0& 0& 0& 0& 0& 0& D'& D\\

\hline {\bf $\chi_{12}$} & 0& 0& 0& 0& 1& A& A'& 0& 0& A& A'& 1& 1 \\

\hline {\bf $\chi_{13}$}&  0& 0& 0& 0& 1& A'& A& 0& 0& A'& A& 1& 1\\

\hline {\bf $\chi_{14}$} & 1& 0& 1& 0& 0& -1& -1& 0& 0& -1& -1& 0& 0\\

\hline {\bf $\chi_{15}$} & -1& 0& 1& 0& -1& 0& 0& 0& 0& -A& -A'& 0& 0\\

\hline {\bf $\chi_{16}$} & -1& 0& 1& 0& -1& 0& 0& 0& 0& -A'& -A& 0& 0\\

\hline {\bf $\chi_{17}$} & 1& 0& 0& -1& 0& 0& 0& 0& 0& 1& 1& 0& 0\\

\hline {\bf $\chi_{18}$} & -1& 1& 0& 0& -1& 0& 0& 1& 1& 0& 0& 0& 0\\

\hline {\bf $\chi_{19}$} & 0& -1& 0& -1& 0& 1& 1& -1& -1& 1& 1& 0& 0 \\

\hline {\bf $\chi_{20}$} & -1& 1& 0& 0& 0& 0& 0& 0& 0& -1& -1& 0& 0\\

\hline {\bf $\chi_{21}$}&  0& 1& 1& 0& 0& -1& -1& 0& 0& 1& 1& 0& 0\\

\hline {\bf $\chi_{22}$} & 0& 0& 0& 0& 0& 1& 1& 0& 0& -1& -1& 1& 1\\

\hline {\bf $\chi_{23}$} & -1& -1& 0& 1& 0& 0& 0& 0& 0& 0& 0& 0& 0\\

\hline {\bf $\chi_{24}$} & 0& -1& 0& 1& 0& 0& 0& -1& -1& 0& 0& 1& 1\\

\hline {\bf $\chi_{25}$} & 0& 1& -1& -1& 0& 0& 0& 1& 1& 0& 0& 0& 0\\

\hline {\bf $\chi_{26}$} & 1& 0& 0& 0& 0& 0& 0& 0& 0& 0& 0& -1& -1\\

\hline
\end{tabular}}
\end{center}
\caption{Character table of $M_{24}$ (ii).}\label{tablacarM24ii}
\end{table}

In the following statement, we summarize our study by mean of
abelian subracks in the group $M_{24}$.

\begin{theorem}\label{teorM24}
Let $\rho\in \widehat{M_{24}^{s_j}}$, with $1\leq j \leq 26$. The
braiding is negative in the cases $j=6$, with $\rho=\rho_{2,6}$ or
$\rho_{3,6}$, $j=8$, with $\rho=\rho_{2,8}$ or $\rho_{3,8}$,
$j=14$, with $\rho=\epsilon\otimes \chi_{(-1)}$ or $\sgn\otimes
\chi_{(-1)}$, $j=17$, $18$, $19$ and $20$, with
$\rho=\chi_{(-1)}$. Otherwise, $\dim\toba(\oc_{s_j},\rho)=\infty$.
\end{theorem}

\pf

\emph{CASE: $j=4$, $5$, $9$, $16$}. From Tables \ref{tablacarM24i}
and \ref{tablacarM24ii}, we see that $s_j$ is real. By Lemma
\ref{odd}, $\dim\toba(\oc_{s_j}, \rho)=\infty$, for all $\rho \in
\widehat{M_{24}^{s_j}}$.

\smallbreak

\emph{CASE: $j=12$, $13$, $21$, $22$, $23$, $24$, $25$, $26$}. We
compute that $s_j^2$, $s_j^{4}\in \oc_{s_j}$. By Lemma
\ref{lemaB}, $\dim\toba(\oc_{s_j}, \rho)=\infty$, for all $\rho
\in \widehat{M_{24}^{s_j}}$, since $|s_j|$ is odd.

\smallbreak

\emph{CASE: $j=19$, $20$}. We have that $|s_j|=14$ and
$M_{24}^{s_j}\simeq \Z_{14}$. Although $s_j$ is not real, we
compute that $s_j^{9}$, $s_j^{11}\in \oc_{s_j}$. Thus, if
$q_{s_js_j}\neq -1$, then $\dim\toba(\oc_{s_j}, \rho)=\infty$, by
Lemma \ref{lemaB}. The remained case corresponds to
$\rho(s_j)=\omega_{14}^7=-1$, which satisfies $q_{s_js_j}= -1$. We
compute that $\oc_{s_j}\cap M_{24}^{s_j}=\{s_j,s_j^9,s_j^{11}\}$.
It is straightforward to prove that the braiding is negative.

\smallbreak

\emph{CASE: $j=17$, $18$}. We have that $|s_j|=12$ and
$M_{24}^{s_j}\simeq \Z_{12}$. Also, we compute that $s_j$ is real.
Thus, if $q_{s_js_j}\neq -1$, then $\dim\toba(\oc_{s_j},
\rho)=\infty$, by Lemma \ref{odd}. The remained case corresponds
to $\rho(s_j)=\omega_{12}^6=-1$, which satisfies $q_{s_js_j}= -1$.
We compute that $\oc_{s_j}\cap M_{24}^{s_j}=\{s_j,s_j^5,
s_j^7,s_j^{11}\}$. It is straightforward to prove that the
braiding is negative.

\smallbreak

\emph{CASE: $j=15$}. The representative $s_{15}$ is
$$( 1,11)( 2, 3,14,13, 9,23, 7,17, 4,16)( 5,19,10,20,18, 6,21,15,22, 8)(12,24),$$
it has order 10 and it is real. We compute that
$M_{24}^{s_{15}}=\langle x, s_{15} \rangle\simeq \Z_2 \times
\Z_{10}$, where
\begin{align*}
x&:= ( 1,12)( 2,18)( 3, 6)( 4,10)( 5, 7)( 8,23)(
9,22)(11,24)(13,15)(14,21)\\ & \qquad (16,20) (17,19).
\end{align*}
Let us define $\{\nu_0,\dots,\nu_9\}$, where
$\nu_l(s_{15}):=\omega_{10}^l$, $0 \leq l \leq 9$. So,
\begin{align*} \widehat{M_{24}^{s_{15}}}=\{\epsilon
\otimes \nu_l,\, \sgn \otimes \nu_l \,\, | \,\, 0 \leq l\leq 9 \},
\end{align*}
where $\epsilon$ and $\sgn$ mean the trivial and the sign
representations of $\Z_2$, respectively. If $\rho \in
\widehat{M_{24}^{s_{15}}}$, with $l\neq 5$, then
$q_{s_{15}s_{15}}\neq -1$, and $\dim\toba(\oc_{s_{15}},
\rho)=\infty$, by Lemma \ref{odd}. The remained two cases are
$\rho=\epsilon \otimes \nu_5$ and $\rho=\sgn \otimes \nu_5$. We
will prove that the Nichols algebra $\toba(\oc_{s_{15}}, \rho)$ is
infinite-dimensional. First, we compute that $\oc_{s_{15}}\cap
M_{24}^{s_{15}}$ has 12 elements, and it contains
$\sigma_1:=s_{15}$, {\small{\begin{align*} \sigma_2&:=( 1,12)( 2,
5,14,10, 9,18, 7,21, 4,22)( 3,19,13,20,23, 6,17,15,16, 8)(11,24),\\
\sigma_3&:=( 1,24)( 2, 6,14,15, 9, 8, 7,19, 4,20)( 3,21,13,22,23,
5,17,10,16,18)(11,12).
\end{align*}}}
We compute that $\sigma_2\sigma_3=\sigma_{1}^{7}$. We choose
$g_1=\id$,
\begin{align*}
g_2&:=( 1,11,24)( 3, 5, 6)( 8,23,18)(10,15,13)(16,22,20)(17,21,19)
\end{align*}
and $g_3:=g_2^{-1}$. These elements are in $M_{24}$ and they
satisfy the same relations as in \eqref{eq1bis}, \eqref{eq2bis}
and \eqref{eq3bis}.

Assume that $\rho=\epsilon \otimes \nu_5$ or $\sgn \otimes \nu_5$.
If $W:=\ku$ - span of $\{g_1,g_2,g_3\}$, then $W$ is a braided
vector subspace of $M(\oc_{s_{15}},\rho)$ of Cartan type whose
Cartan matrix is as in \eqref{cartanmatrix:3x3}. Therefore,
$\dim\toba(\oc_{s_{15}},\rho)=\infty$.

\smallbreak

\emph{CASE: $j=14$}. The representative is
$$s_{14}=( 1,22,18,14,11,19,16, 7)( 2,13)( 3,12, 6,21, 8,10,20,23)( 9,17,15,24),$$
it has order 8 and it is real. We compute that
$M_{24}^{s_{14}}=\langle x, s_{14} \rangle\simeq \Z_2 \times
\Z_{8}$, where
\begin{align*}
x&:=  ( 1, 3)( 2,13)( 4, 5)( 6,18)( 7,23)( 8,11)(
9,15)(10,19)(12,22)(14,21)\\ & \qquad (16,20)(17,24).
\end{align*}
Let us define $\{\nu_0,\dots,\nu_7\}$, where
$\nu_l(s_{14}):=\omega_{8}^l$, $0 \leq l \leq 7$. So,
\begin{align*}
\widehat{M_{24}^{s_{14}}}=\{\epsilon \otimes \nu_l,\, \sgn \otimes
\nu_l \,\, | \,\, 0 \leq l\leq 7 \},
\end{align*}
where $\epsilon$ and $\sgn$ mean the trivial and the sign
representations of $\Z_2$, respectively. If $\rho \in
\widehat{M_{24}^{s_{14}}}$, with $l\neq 4$, then
$q_{s_{14}s_{14}}\neq -1$, and $\dim\toba(\oc_{s_{14}},
\rho)=\infty$, by Lemma \ref{odd}. The remained two cases
corresponding to $\rho=\epsilon \otimes \nu_4$ and $\rho=\sgn
\otimes \nu_4$ have negative braiding. Indeed, we compute that
\begin{align*}
\oc_{s_{14}}\cap M_{24}^{s_{14}}=\{s_{14}^{-1},\,\, x
s_{14}^{-3},\,\, x s_{14},\,\, s_{14}^3,\,\, s_{14}^{-3},\,\, x
s_{14}^3,\,\, s_{14},\,\, x s_{14}^{-1} \}.
\end{align*}
For simplicity, we write $\sigma_1:=s_{14}^{-1}$, $\sigma_2:=x
s_{14}^{-3}$, $\sigma_3:=x s_{14}$, $\sigma_4:=s_{14}^3$,
$\sigma_5:=s_{14}^{-3}$, $\sigma_6:=x s_{14}^3$,
$\sigma_7:=s_{14}$ and $\sigma_8:=x s_{14}^{-1}$. We choose in
$M_{24}$ the following elements:
\begin{align*}
g_1&:=( 3, 8)( 4, 5)( 7,22)( 9,15)(10,23)(12,21)(14,19)(16,18),\\
g_2&:=( 2, 4)( 5,13)( 7,21)( 9,17)(10,22)(12,19)(14,23)(15,24),\\
g_3&:=( 2, 5)( 4,13)( 7,23)( 9,24)(10,19)(12,22)(14,21)(15,17),\\
g_4&:=( 2,13)( 3, 8)( 7,19)(10,21)(12,23)(14,22)(16,18)(17,24),\\
g_5&:=( 2,13)( 4, 5)( 7,14)( 9,15)(10,12)(17,24)(19,22)(21,23),\\
g_6&:=( 2, 4,13, 5)( 3, 8)(
7,10,14,12)(9,24,15,17)(16,18)(19,23,22,21),
\end{align*}
$g_7:=\id$ and $g_8:=g_6^{-1}$. We compute that these elements
satisfy $\sigma_k g_7=g_7 \sigma_k$, $1\leq k \leq 8$, $\sigma_7
g_k=g_k \sigma_k$, $1\leq k \leq 5$, $\sigma_7 g_6=g_6 \sigma_8$
and $\sigma_7 g_8=g_8 \sigma_6$. It is easy to see that, if
$\rho=\epsilon \otimes \nu_4$ or $\rho=\sgn \otimes \nu_4$, then
$\rho(\gamma_{k,7}\gamma_{7,k})=1$, for every $1\leq k\leq 8$.
From Lemma \ref{lema:treneg}, the braiding is negative.

\smallbreak

\emph{CASE: $j=10$}. The representative $s_{10}$ is
$$( 1,20)( 3, 4,16)( 5,14,21,19,23,15)( 7,11,12,24,18,13)( 8,22,10)( 9,17),$$
it has order 6 and it is real. We compute that the centralizer
$M_{24}^{s_{10}}$ is a non-abelian group of order 24.

Let $\rho=(\rho,V)\in \widehat{M_{24}^{s_{10}}}$. We will prove
that the Nichols algebra $\toba(\oc_{s_{10}}, \rho)$ is
infinite-dimensional. If $q_{s_{10}s_{10}}\neq -1$, then the
result follows from Lemma \ref{odd}. Assume that
$q_{s_{10}s_{10}}= -1$. We compute that $\oc_{s_{10}}\cap
M_{24}^{s_{10}}$ has 10 elements and it contains
$\sigma_1:=s_{10}$, $\sigma_2:=s_{10}^{-1}$,
\begin{align*}
\sigma_3&:= ( 1,20)( 2, 6)( 3, 8,16,10, 4,22)( 5,14,21,19,23,15)(
7,18,12)(11,13,24)
\end{align*}
and $\sigma_4:=\sigma_3^{-1}$. These elements commute each other.
We choose $g_1:=\id$, $g_2:=( 2, 6)( 3,10)( 4,22)(
8,16)(11,13)(12,18)(14,15)(21,23)$, {\small{\begin{align*} g_3&:=(
1, 6,17)( 2, 9,20)( 3, 7, 5)( 4,18,23)(
8,11,14)(10,24,19)(12,21,16) (13,15,22)
\end{align*}}}
and $g_4:=g_3g_2$. These elements are in $M_{24}$ and they satisfy
$\sigma_r g_r=g_r \sigma_1$, $\sigma_r g_1=g_1 \sigma_r$, $1\leq r
\leq 4$, and
\begin{align*}
\sigma_1 g_2&=g_2 \sigma_2, \quad & \sigma_1 g_3&=g_3
\gamma,\quad &\sigma_1 g_4&=g_4 \gamma^{-1},\\
\sigma_3 g_2&=g_2 \sigma_4, \quad & \sigma_2 g_3&=g_3
\gamma^{-1},\quad &\sigma_2 g_4&=g_4 \gamma,\\
\sigma_4 g_2&=g_2 \sigma_3, \quad & \sigma_4 g_3&=g_3 \sigma_2,
\quad &\sigma_3 g_4&=g_4 \sigma_2,
\end{align*}
where
\begin{align*}
\gamma&:=( 2, 6)( 3, 8,16,10, 4,22)( 5,23,21)( 7,11,12,24,18,13)(
9,17)(14,15,19).
\end{align*}
Also, we compute that $\sigma_3 \gamma=s_{10}^{-1}$. Now, we
define $W:=\ku$ - span of $\{g_1v,g_2v,g_3v,g_4v\}$, where $v \in
V-0$. Hence, it is straightforward to check that $W$ is a braided
vector subspace of $M(\oc_{s_{10}},\rho)$ of Cartan type whose
associated Cartan matrix is given by \eqref{cartanmatrix:4x4}. By
Theorem \ref{cartantype}, $\dim\toba(\oc_{s_{10}},\rho)=\infty$.

\smallbreak

\emph{CASE: $j=11$}. We compute that $M_{24}^{s_{11}}\simeq
M_{24}^{s_{10}}$. This implies that this case is analogous to the
previous case, since $\oc_{s_{11}}\simeq \oc_{s_{10}}$ as racks.

\smallbreak

\emph{CASE: $j=8$}. The representative $s_{8}$ is
$$( 1,22,11,16)( 2, 8,14,17)( 3,23,24, 6)( 4,21,12,13)( 5,15,19, 9)( 7,18,20,10),$$
it has order 4 and it is real. We compute that the centralizer
$M_{24}^{s_{8}}$ is a non-abelian group of order 96 whose
character table is given by Table \ref{tablacarM24^s8}.

{\tiny{\begin{table}[t]
\begin{center}
\tiny{\hspace*{-0.3cm}\begin{tabular}{|c|c|c|c|c|c|c|c|c|c|c|c|c|c|c|c|c|c|c|c|c|c|c|}
\hline
{\bf $k$} &    1& 2& 3& 4& 5& 6& 7& 8& 9& 10& 11& 12& 13& 14 & 15 & 16 & 17& 18 & 19 & 20 \\
\hline \hline
{\bf $|y_k|$} &   1& 2& 2& 4& 4& 12& 4& 4& 12& 4& 3& 4& 6& 4& 2& 2& 2& 4& 4& 4\\
\hline \hline

\hline {\bf $\mu_1$} & 1& 1& 1& 1& 1& 1& 1& 1& 1& 1& 1& 1& 1& 1& 1& 1& 1& 1& 1& 1 \\

\hline {\bf $\mu_2$} & 1& -1& 1& -1& 1& -1& 1& 1& -1& 1& 1& -1& 1& -1& -1& 1& 1& -1& -1& -1\\

\hline {\bf $\mu_3$} & 1& 1& 1& -1& -1& -1& -1& -1& -1& -1& 1& 1& 1& 1& 1& 1& 1& -1& -1& -1\\

\hline {\bf $\mu_4$} & 1& -1& 1& 1& -1& 1& -1& -1& 1& -1& 1& -1& 1& -1& -1& 1& 1& 1& 1& 1\\

\hline {\bf $\mu_5$} & 1& -1& 1& -i& i& i& -i& -i& -i& i& 1& -1&
-1& 1& 1&
      -1& -1& i& -i& i\\

\hline {\bf $\mu_6$} & 1& 1& 1& -i& -i& i& i& i& -i& -i& 1& 1& -1&
-1& -1&
      -1& -1& i& -i& i\\

\hline {\bf $\mu_7$} & 1& -1& 1& i& -i& -i& i& i& i& -i& 1& -1&
-1& 1& 1&
      -1& -1& -i& i& -i\\

\hline {\bf $\mu_8$} & 1& 1& 1& i& i& -i& -i& -i& i& i& 1& 1& -1&
-1& -1&
      -1& -1& -i& i& -i\\

\hline {\bf $\mu_9$} & 2& 0& 2& 2& 0& -1& 0& 0& -1& 0& -1& 0& -1& 0& 0& 2& 2& 2& 2& 2\\

\hline {\bf $\mu_{10}$} &  2& 0& 2& -2& 0& 1& 0& 0& 1& 0& -1& 0& -1& 0& 0& 2& 2& -2& -2& -2\\

\hline {\bf $\mu_{11}$} & 2& 0& 2& -2i& 0& -i& 0& 0& i& 0& -1& 0&
1& 0& 0& -2& -2&
     2i& -2i& 2i\\

\hline {\bf $\mu_{12}$} & 2& 0& 2& 2i& 0& i& 0& 0& -i& 0& -1& 0&
1& 0& 0& -2& -2&
      -2i& 2i& -2i\\

\hline {\bf $\mu_{13}$}& 3& -1& -1& 1& -1& 0& 1& -1& 0& 1& 0& 1& 0& 1& -1& -1& 3& 1& -3& -3\\

\hline {\bf $\mu_{14}$} &3& 1& -1& 1& 1& 0& -1& 1& 0& -1& 0& -1& 0& -1& 1& -1& 3& 1& -3& -3\\

\hline {\bf $\mu_{15}$} & 3& -1& -1& -1& 1& 0& -1& 1& 0& -1& 0& 1& 0& 1& -1& -1& 3& -1& 3& 3\\

\hline {\bf $\mu_{16}$} &3& 1& -1& -1& -1& 0& 1& -1& 0& 1& 0& -1& 0& -1& 1& -1& 3& -1& 3& 3\\

\hline {\bf $\mu_{17}$}& 3& -1& -1& i& -i& 0& -i& i& 0& i& 0& 1&
0& -1& 1& 1& -3&
      -i& -3i& 3i \\

\hline {\bf $\mu_{18}$} & 3& 1& -1& i& i& 0& i& -i& 0& -i& 0& -1&
0& 1& -1& 1& -3&
      -i& -3i& 3i\\

\hline {\bf $\mu_{19}$} & 3& -1& -1& -i& i& 0& i& -i& 0& -i& 0& 1&
0& -1& 1& 1& -3&
      i& 3i& -3i\\

\hline {\bf $\mu_{20}$} & 3& 1& -1& -i& -i& 0& -i& i& 0& i& 0& -1&
0& 1& -1& 1& -3&
      i& 3i& -3i\\

\hline
\end{tabular}}
\end{center}
\caption{Character table of $M_{24}^{s_8}$.}\label{tablacarM24^s8}
\end{table}}}

For every $k$, $1 \leq k \leq 20$, we call $\rho_k=(\rho_k,V_k)$
the irreducible representation of $M_{24}^{s_8}$ whose character
is $\mu_k$. We compute that $s_8\in \oc_{20}^{M_{24}^{s_8}}$.

From Table \ref{tablacarM24^s8}, we have that if $k\neq 2$, $3$,
$10$, $13$, $14$, then $q_{s_8s_8}\neq -1$ and
$\dim\toba(\oc_{s_{8}}, \rho_k)=\infty$, by Lemma \ref{odd}. On
the other hand, we compute that $\oc_{s_8}\cap M_{24}^{s_8}$ has
32 elements and it contains {\small{\begin{align*} \sigma_1&:=( 1,
2,11,14)( 3,18, 9,13)( 4,23,20, 5)( 6, 7,19,12)(
8,16,17,22)(10,15,21,24),\\
\sigma_2&:=( 1, 2,11,14)( 3,21, 9,10)( 4,19,20, 6)( 5, 7,23,12)(
8,16,17,22)(13,15,18,24),\\
\sigma_5&:=( 1,16,11,22)( 2,17,14, 8)( 3, 5,24,19)( 4,10,12,18)(
6, 9,23,15)( 7,13,20,21),
\end{align*}}}

\vspace*{-0.5cm} \noindent $\sigma_3:=\sigma_2^{-1}$,
$\sigma_4:=\sigma_1^{-1}$, $\sigma_6:=s_{8}^{-1}$,
$\sigma_7:=\sigma_5^{-1}$ and $\sigma_8:=s_{8}$. We choose in
$M_{24}$
\begin{align*}
g_1&:= ( 2,16,14,22)( 3, 9,15,24)( 4, 7)( 5,21,19,10)(
6,18,23,13)(8,17) ,\\
g_2&:= ( 2,16,14,22)( 4,12, 7,20)( 5,13,23,21)(
6,10,19,18)( 8,17)( 9,24) ,\\
g_3&:= ( 2,22,14,16)( 3, 9,15,24)(
5,13,19,18)( 6,10,23,21)( 8,17)(12,20) ,\\
g_6&:= ( 2,14)( 3, 4,15, 7)( 5,18, 6,21)( 9,20,24,12)(10,23,13,19)(16,22) ,\\
g_7&:= ( 3, 4)( 5,13)( 6,10)( 7,15)( 9,20)(12,24)(18,23)(19,21),
\end{align*}
$g_4:=g_2^{-1}$, $g_5:=g_2^{2}$ and $g_8:=\id$. We compute that
these elements satisfy $\sigma_k g_l=g_l \gamma_{k,l}$, where
$\gamma_{k,k}=s_{8}$, $\gamma_{k,8}=\sigma_k$, $1\leq k \leq 8$,
$\gamma_{8,1}=\sigma_3$, $\gamma_{8,2}=\sigma_4$,
$\gamma_{8,3}=\sigma_1$, $\gamma_{8,4}=\sigma_2$,
$\gamma_{8,5}=\sigma_5$, $\gamma_{8,6}=\sigma_6$,
$\gamma_{8,7}=\sigma_7$, and
\begin{align*}
\gamma_{1,2}&=\sigma_7,
&\gamma_{1,3}&=\sigma_5,&\gamma_{1,4}&=\sigma_6,
&\gamma_{1,5}&=\sigma_3, &\gamma_{1,6}&=\sigma_3,
&\gamma_{1,7}&=\sigma_1,\\
\gamma_{2,1}&=\sigma_7, &\gamma_{2,3}&=\sigma_6,
&\gamma_{2,4}&=\sigma_5, &\gamma_{2,5}&=\sigma_4,
&\gamma_{2,6}&=\sigma_4, &\gamma_{2,7}&=\sigma_2,\\
\gamma_{3,1}&=\sigma_5, &\gamma_{3,2}&=\sigma_6,
&\gamma_{3,4}&=\sigma_7, &\gamma_{3,5}&=\sigma_1,
&\gamma_{3,6}&=\sigma_1, &\gamma_{3,7}&=\sigma_3,\\
\gamma_{4,1}&=\sigma_6, &\gamma_{4,2}&=\sigma_5,
&\gamma_{4,3}&=\sigma_7, &\gamma_{4,5}&=\sigma_2,
&\gamma_{4,6}&=\sigma_2, &\gamma_{4,7}&=\sigma_4,\\
\gamma_{5,1}&=\sigma_1, &\gamma_{5,2}&=\sigma_2,
&\gamma_{5,3}&=\sigma_3, &\gamma_{5,4}&=\sigma_4,
&\gamma_{5,6}&=\sigma_7, &\gamma_{5,7}&=\sigma_6,\\
\gamma_{6,1}&=\sigma_2, &\gamma_{6,2}&=\sigma_1,
&\gamma_{6,3}&=\sigma_4, &\gamma_{6,4}&=\sigma_3,
&\gamma_{6,5}&=\sigma_7, &\gamma_{6,7}&=\sigma_5,\\
\gamma_{7,1}&=\sigma_4, &\gamma_{7,2}&=\sigma_3,
&\gamma_{7,3}&=\sigma_2, &\gamma_{7,4}&=\sigma_1,
&\gamma_{7,5}&=\sigma_6, &\gamma_{7,6}&=\sigma_5.
\end{align*}

Assume that $k=10$, $13$ or $14$; thus, $q_{s_8s_8}=-1$. We check
\emph{case by case} that always we can construct a braided vector
subspace of $M(\oc_{s_{8}}, \rho_k)$ of diagonal type whose
generalized Dynkin diagram contains an $r$-cycle with $r>3$ or a
vertex with valency greater than 3. By Lemma \ref{Hecke},
$\dim\toba(\oc_{s_{8}}, \rho_k)=\infty$.

Finally, assume that $k= 2$ or $3$. Thus, $q_{s_8s_8}= -1$ and we
compute that $\rho(\gamma_{1,t}\gamma_{t,1})=1$, for every $1\leq
t \leq 32$. By Lemma \ref{lema:treneg}, the braiding is negative.


\smallbreak

\emph{CASE: $j=7$}. The representative is
$$s_{7}=( 1,18,11,16)( 3, 6, 8,20)( 7,22,14,19)( 9,15)(10,23,12,21)(17,24),$$
it has order 4 and it is real. We compute that the centralizer
$M_{24}^{s_{7}}$ is a non-abelian group of order 128.

Let $\rho=(\rho,V)\in \widehat{M_{24}^{s_{7}}}$. We will prove
that the Nichols algebra $\toba(\oc_{s_{7}}, \rho)$ is
infinite-dimensional. If $q_{s_{7}s_{7}}\neq -1$, then the result
follows from Lemma \ref{odd}. Assume that $q_{s_{7}s_{7}}= -1$. We
compute that $\oc_{s_{7}}\cap M_{24}^{s_{7}}$ has 40 elements and
it contains
\begin{align*}
\sigma_1&:= ( 2, 4,13, 5)( 3, 8)( 6,20)( 7,19,14,22)(
9,17,15,24)(10,23,12,21),\\
\sigma_3&:=( 1,16,11,18)( 2, 5,13, 4)( 3, 6, 8,20)(
9,17,15,24)(10,12)(21,23),
\end{align*}
$\sigma_2:=\sigma_1^{-1}$ and $\sigma_4:=\sigma_3^{-1}$. These
elements commute each other and $\sigma_1\sigma_3=s_{7}^{-1}$.
Now, we choose
\begin{align*} g_1&:=( 1,
2,18, 4,11,13,16, 5)(
3, 9,20,24, 8,15, 6,17)(10,21,12,23)(19,22),\\
g_2&:=( 1, 2,11,13)( 3, 9, 8,15)( 4,18, 5,16)(
6,24,20,17)(10,23)(12,21),\\
g_3&:=( 2,19, 4,14,13,22, 5, 7)( 3, 6, 8,20)(
9,23,17,12,15,21,24,10)(16,18),\\
g_4&:=( 2,14,13, 7)( 3, 6)( 4,19, 5,22)( 8,20)(
9,10,15,12)(17,21,24,23).
\end{align*}
These elements are in $M_{24}$ and they satisfy
\begin{align*}
\sigma_1 g_1&=g_1 s_{7}, \quad & \sigma_1 g_2&=g_2s_{7}^{-1},
\quad &\sigma_1 g_3&=g_3 \sigma_1,
\quad &\sigma_1 g_4&=g_4 \sigma_2,\\
\sigma_2 g_1&=g_1 s_{7}^{-1}, \quad & \sigma_2 g_2&=g_2 s_{7},
\quad &\sigma_2 g_3&=g_3 \sigma_2,
\quad &\sigma_2 g_4&=g_4 \sigma_1,\\
\sigma_3 g_1&=g_1 \sigma_3, \quad & \sigma_3 g_2&=g_2 \sigma_4,
\quad &\sigma_3 g_3&=g_3 s_{7},
\quad &\sigma_3 g_4&=g_4 s_{7}^{-1},\\
\sigma_4 g_1&=g_1 \sigma_4, \quad & \sigma_4 g_2&=g_2\sigma_3,
\quad &\sigma_4 g_3&=g_3 s_{7}^{-1} , \quad &\sigma_4 g_4&=g_4
s_{7}.
\end{align*}
We define $W:=\ku$ - span of $\{g_1v,g_2v,g_3v,g_4v\}$, where $v
\in V-0$. Hence, it is straightforward to check that $W$ is a
braided vector subspace of $M(\oc_{s_{7}},\rho)$ of Cartan type
whose associated Cartan matrix is given by
\eqref{cartanmatrix:4x4}. By Theorem \ref{cartantype},
$\dim\toba(\oc_{s_{7}},\rho)=\infty$.


\smallbreak

\emph{CASE: $j=6$}. The representative $ s_{6}$ is
\begin{align*}
( 1, 9,20,17)( 2, 6)( 3,10)( 4, 8)( 5,24,19,
7)(11,14,18,23)(12,21,13,15)(16,22)
\end{align*}
and it is real. We compute that $M_{24}^{s_{6}}$ is a non-abelian
group of order 384 whose character table is given by Tables
\ref{tablacarM24^s6i} and \ref{tablacarM24^s6ii}.

\begin{table}[t]
\begin{center}
\tiny{\begin{tabular}{|c|c|c|c|c|c|c|c|c|c|c|c|c|c|}
\hline {\bf $k$} &    1& 2& 3& 4& 5& 6& 7& 8& 9& 10& 11& 12& 13\\
\hline \hline {\bf $|y_k|$} &   1& 3& 2& 2& 4& 2& 2& 6& 4& 4& 4&
2& 4\\

\hline \hline

\hline {\bf $\mu_1$} & 1& 1& 1& 1& 1& 1& 1& 1& 1& 1& 1& 1& 1\\

\hline {\bf $\mu_2$} & 1& 1& 1& -1& -1& 1& 1& 1& 1& -1& -1& -1&
1\\

\hline {\bf $\mu_3$} & 1& 1& 1& 1& 1& 1& 1& 1& 1& 1& 1& 1& -1\\

\hline {\bf $\mu_4$} & 1& 1& 1& -1& -1& 1& 1& 1& 1& -1& -1& -1&
-1\\

\hline {\bf $\mu_5$} & 2& -1& 2& 0& 0& 2& 2& -1& 2& 0& 0& 0& 0\\

\hline {\bf $\mu_6$} & 2& -1& 2& 0& 0& 2& 2& -1& 2& 0& 0& 0& 0\\

\hline {\bf $\mu_7$} &  3& 0& -1& -1& 1& 3& -1& 0& -1& -1& 1& -1&
-1\\

\hline {\bf $\mu_8$} & 3& 0& -1& 1& -1& 3& -1& 0& -1& 1& -1& 1&
1\\

\hline {\bf $\mu_9$} & 3& 0& -1& -1& 1& 3& -1& 0& -1& -1& 1& -1&
1\\

\hline {\bf $\mu_{10}$} & 3& 0& -1& 1& -1& 3& -1& 0& -1& 1& -1& 1&
-1\\

\hline {\bf $\mu_{11}$} &3& 0& 3& 1& 1& -1& -1& 0& -1& -1& -1& 1&
1\\

\hline {\bf $\mu_{12}$} & 3& 0& -1& 1& -1& -1& 3& 0& -1& -1& 1& 1&
1\\

\hline {\bf $\mu_{13}$}& 3& 0& -1& -1& 1& -1& 3& 0& -1& 1& -1& -1&
1\\

\hline {\bf $\mu_{14}$} & 3& 0& -1& 1& -1& -1& 3& 0& -1& -1& 1& 1&
-1\\

\hline {\bf $\mu_{15}$} & 3& 0& 3& -1& -1& -1& -1& 0& -1& 1& 1&
-1& 1\\

\hline {\bf $\mu_{16}$} & 3& 0& 3& 1& 1& -1& -1& 0& -1& -1& -1& 1&
-1\\

\hline {\bf $\mu_{17}$} & 3& 0& -1& -1& 1& -1& 3& 0& -1& 1& -1&
-1& -1\\

\hline {\bf $\mu_{18}$} &   3& 0& 3& -1& -1& -1& -1& 0& -1& 1& 1&
-1& -1\\

\hline {\bf $\mu_{19}$} & 4& 1& 0& 2& 0& 0& 0& -1& 0& 0& 0& -2&
2i\\

\hline {\bf $\mu_{20}$} & 4& 1& 0& -2& 0& 0& 0& -1& 0& 0& 0& 2&
-2i\\

\hline {\bf $\mu_{21}$} &  4& 1& 0& -2& 0& 0& 0& -1& 0& 0& 0& 2&
2i\\

\hline {\bf $\mu_{22}$} & 4& 1& 0& 2& 0& 0& 0& -1& 0& 0& 0& -2&
-2i\\

\hline {\bf $\mu_{23}$} &  6& 0& -2& 0& 0& -2& -2& 0& 2& 0& 0& 0&
0\\

\hline {\bf $\mu_{24}$} & 6& 0& -2& 0& 0& -2& -2& 0& 2& 0& 0& 0&
0\\

\hline {\bf $\mu_{25}$} & 8& -1& 0& 0& 0& 0& 0& 1& 0& 0& 0& 0& 0\\

\hline {\bf $\mu_{26}$} &  8& -1& 0& 0& 0& 0& 0& 1& 0& 0& 0& 0&
0\\

\hline
\end{tabular}}
\end{center}
\caption{Character table of $M_{24}^{s_6}$
(i).}\label{tablacarM24^s6i}
\end{table}

\begin{table}[t]
\begin{center}
\tiny{\begin{tabular}{|c|c|c|c|c|c|c|c|c|c|c|c|c|c|}
\hline {\bf $k$} & 14 & 15 & 16 & 17& 18& 19& 20& 21& 22& 23 & 24 & 25 & 26 \\
\hline \hline {\bf $|y_k|$} &  4& 4& 4& 4& 12& 4& 12& 2&
4& 4& 4&  4& 2\\

\hline \hline

\hline {\bf $\mu_1$} & 1& 1& 1& 1& 1& 1& 1& 1&
1& 1& 1& 1& 1\\

\hline {\bf $\mu_2$} & 1& 1& 1& 1& -1& -1& -1& -1&
      -1& -1& -1& -1& 1 \\

\hline {\bf $\mu_3$} & -1& -1& -1& -1& -1& -1& -1& -1&
      -1& -1& -1& -1& 1\\

\hline {\bf $\mu_4$} & -1& -1& -1& -1& 1& 1& 1& 1&
      1& 1& 1& 1& 1\\

\hline {\bf $\mu_5$} & 0& 0& 0& 0& 1& -2& 1& -2& -2&
      -2& -2& -2& 2\\

\hline {\bf $\mu_6$} & 0& 0& 0& 0& -1& 2& -1& 2& 2& 2&
      2& 2& 2\\

\hline {\bf $\mu_7$} & -1& -1& 1& 1& 0& 3& 0& -1&
      -1& 3& -1& 3& 3\\

\hline {\bf $\mu_8$} & 1& 1& -1& -1& 0& 3& 0& -1&
      -1& 3& -1& 3& 3\\

\hline {\bf $\mu_9$} & 1& 1& -1& -1& 0& -3& 0& 1&
      1& -3& 1& -3& 3\\

\hline {\bf $\mu_{10}$} & -1& -1& 1& 1& 0& -3& 0& 1&
      1& -3& 1& -3& 3\\

\hline {\bf $\mu_{11}$} & -1& 1& -1& 1& 0& -1& 0& -1&
      -1& 3& 3& 3& 3 \\

\hline {\bf $\mu_{12}$} & -1& 1& 1& -1& 0& -1& 0& -1&
      3& 3& -1& 3& 3\\

\hline {\bf $\mu_{13}$}& -1& 1& 1& -1& 0& 1& 0& 1&
      -3& -3& 1& -3& 3\\

\hline {\bf $\mu_{14}$} & 1& -1& -1& 1& 0& 1& 0& 1&
      -3& -3& 1& -3& 3\\

\hline {\bf $\mu_{15}$} & -1& 1& -1& 1& 0& 1& 0& 1& 1&
      -3& -3& -3& 3\\

\hline {\bf $\mu_{16}$} & 1& -1& 1& -1& 0& 1& 0& 1& 1&
      -3& -3& -3& 3\\

\hline {\bf $\mu_{17}$} & 1& -1& -1& 1& 0& -1& 0& -1&
      3& 3& -1& 3& 3\\

\hline {\bf $\mu_{18}$} & 1& -1& 1& -1& 0& -1& 0& -1&
      -1& 3& 3& 3& 3\\

\hline {\bf $\mu_{19}$} & 0& -2i& 0& 0& i& 0&
      -i& 0& 0& -4i& 0& 4i& -4\\

\hline {\bf $\mu_{20}$} & 0& 2i& 0& 0& i& 0&
      -i& 0& 0& -4i& 0& 4i& -4\\

\hline {\bf $\mu_{21}$} & 0& -2i& 0& 0& -i& 0&
      i& 0& 0& 4i& 0& -4i& -4\\

\hline {\bf $\mu_{22}$} & 0& 2i& 0& 0& -i& 0&
      i& 0& 0& 4i& 0& -4i& -4\\

\hline {\bf $\mu_{23}$} & 0& 0& 0& 0& 0& -2& 0& 2& -2& 6&
      -2& 6& 6\\

\hline {\bf $\mu_{24}$} & 0& 0& 0& 0& 0& 2& 0& -2& 2& -6&
      2& -6& 6\\

\hline {\bf $\mu_{25}$} & 0& 0& 0& 0& -i& 0& i& 0& 0&
      -8i& 0& 8i& -8\\

\hline {\bf $\mu_{26}$} & 0& 0& 0& 0& i& 0& -i& 0& 0&
      8i& 0& -8i& -8 \\

\hline
\end{tabular}}
\end{center}
\caption{Character table of $M_{24}^{s_6}$
(ii).}\label{tablacarM24^s6ii}
\end{table}

For every $k$, $1 \leq k \leq 26$, we call $\rho_k=(\rho_k,V_k)$
the irreducible representation of $M_{24}^{s_6}$ whose character
is $\mu_k$. We compute that $s_6\in \oc_{23}^{M_{24}^{s_6}}$. From
Tables \ref{tablacarM24^s6i} and \ref{tablacarM24^s6ii}, we have
that if $k\neq 2$, $3$, $5$, $9$, $10$, $13$, $14$, $15$, $16$,
$24$, then $q_{s_6s_6}\neq -1$ and $\dim\toba(\oc_{s_{6}},
\rho_k)=\infty$, by Lemma \ref{odd}.

We compute that $\oc_{s_6}\cap M_{24}^{s_6}$ has 80 elements and
it contains {\small{\begin{align*} \sigma_1:&=( 1, 7)( 2, 4,
3,16)( 5, 9)( 6,
8,10,22)(11,14,18,23)(12,15,13,21)(17,19) (20,24),\\
\sigma_4:&=( 1, 9,20,17)( 2,10)( 3, 6)( 4,22)( 5,24,19, 7)(
8,16)(11,23,18,14) (12,15,13,21),\\
\sigma_7:&=( 1,24)( 2, 4, 3,16)( 5,17)( 6, 8,10,22)( 7,20)(
9,19)(11,23,18,14) (12,21,13,15),
\end{align*}}}

\vspace*{-0.4cm}\noindent $\sigma_2:=\sigma_1^{-1}$,
$\sigma_3:=s_6$, $\sigma_5:=s_6^{-1}$ $\sigma_6:=\sigma_4^{-1}$
and $\sigma_8:=\sigma_7^{-1}$. We compute that $\sigma_1$,
$\sigma_2$, $\sigma_7$, $\sigma_8\in \oc_{17}^{M_{24}^{s_6}}$,
$\sigma_4$, $\sigma_6 \in \oc_{24}^{M_{24}^{s_6}}$ and
$s_6^{-1}\in \oc_{25}^{M_{24}^{s_6}}$. We choose in $M_{24}$
{\small{\begin{align*} g_1&:= ( 1, 2,17,16)( 3, 9, 4,20)( 5,
8,24,10)( 6,19,22, 7)(11,12,23,21)(13,14,15,18) ,\\
g_4&:= ( 2, 4)( 3,16)( 6,22)( 8,10)(11,14)(12,21)(13,15)(18,23) ,\\
g_5&:=  ( 2,10)( 3, 6)( 4,22)( 7,24)( 8,16)(
9,17)(14,23)(15,21),\\
g_7&:= ( 1, 2,24,10,20, 3, 7, 6)( 4, 5, 8,17,16,19,22,
9)(12,15,13,21)(14,23) ,
\end{align*}}}

\vspace*{-0.5cm}\noindent $g_3:= \id$, $g_2:=g_1g_5$,
$g_6:=g_4g_5$ and $g_8:=g_7g_5$. We compute that these elements
satisfy $\sigma_k g_l=g_l \gamma_{k,l}$, where
$\gamma_{k,k}=s_{6}$, $\gamma_{k,3}=\sigma_k$, $1\leq k \leq 8$, $
\gamma_{3,1}=\sigma_2$, $\gamma_{3,2}=\sigma_8$,
$\gamma_{3,4}=\sigma_4$, $\gamma_{3,5}=\sigma_5$,
$\gamma_{3,6}=\sigma_6$, $\gamma_{3,7}=\sigma_2$,
$\gamma_{3,8}=\sigma_8$, and
\begin{align*}
\gamma_{1,2}&=\sigma_5,
&\gamma_{1,4}&=\sigma_2,&\gamma_{1,5}&=\sigma_7,
&\gamma_{1,6}&=\sigma_8, &\gamma_{1,7}&=\sigma_4,
&\gamma_{1,8}&=\sigma_6,\\
\gamma_{2,1}&=\sigma_5, &\gamma_{2,4}&=\sigma_1,
&\gamma_{2,5}&=\sigma_8, &\gamma_{2,6}&=\sigma_7,
&\gamma_{2,7}&=\sigma_6, &\gamma_{2,8}&=\sigma_4,\\
\gamma_{4,1}&=\sigma_8, &\gamma_{4,2}&=\sigma_2,
&\gamma_{4,5}&=\sigma_6, &\gamma_{4,6}&=\sigma_5,
&\gamma_{4,7}&=\sigma_8, &\gamma_{4,8}&=\sigma_2,\\
\gamma_{5,1}&=\sigma_1, &\gamma_{5,2}&=\sigma_7,
&\gamma_{5,4}&=\sigma_6, &\gamma_{5,6}&=\sigma_4,
&\gamma_{5,7}&=\sigma_1, &\gamma_{5,8}&=\sigma_7,\\
\gamma_{6,1}&=\sigma_7, &\gamma_{6,2}&=\sigma_1,
&\gamma_{6,4}&=\sigma_5, &\gamma_{6,5}&=\sigma_4,
&\gamma_{6,7}&=\sigma_7, &\gamma_{6,8}&=\sigma_1,\\
\gamma_{7,1}&=\sigma_4, &\gamma_{7,2}&=\sigma_6,
&\gamma_{7,4}&=\sigma_8, &\gamma_{7,5}&=\sigma_1,
&\gamma_{7,6}&=\sigma_2, &\gamma_{7,5}&=\sigma_5,\\
\gamma_{8,1}&=\sigma_6, &\gamma_{8,2}&=\sigma_4,
&\gamma_{8,4}&=\sigma_7, &\gamma_{8,5}&=\sigma_2,
&\gamma_{8,6}&=\sigma_1, &\gamma_{8,7}&=\sigma_5.
\end{align*}

Assume that $k= 5$, $9$, $10$, $13$, $14$, $15$, $16$ or $24$;
thus, $q_{s_6s_6}=-1$. We check \emph{case by case} that always we
can construct a braided vector subspace of $M(\oc_{s_{6}},
\rho_k)$ of diagonal type whose generalized Dynkin diagram
contains an $r$-cycle with $r>3$ or a vertex with valency greater
than 3. By Lemma \ref{Hecke}, $\dim\toba(\oc_{s_{6}},
\rho_k)=\infty$.

Finally, assume that $k= 2$ or $3$. Thus, $q_{s_6s_6}= -1$ and we
compute that $\rho(\gamma_{1,t}\gamma_{t,1})=1$, for every $1\leq
t \leq 80$. By Lemma \ref{lema:treneg}, the braiding is negative.


\smallbreak

\emph{CASE: $j=2$}. The representative is
$$s_{2}=( 1,20)( 5,19)( 7,24)( 9,17)(11,18)(12,13)(14,23)(15,21).$$
We compute that the centralizer $M_{24}^{s_{2}}$ is a non-abelian
group of order 21504 whose character table is given by Tables
\ref{tablacarM24^s2i} and \ref{tablacarM24^s2ii}, where $A = (-1+i
\sqrt{7})/2$.

\begin{table}[t]
\begin{center}
\tiny{\begin{tabular}{|c|c|c|c|c|c|c|c|c|c|c|c|c|c|c|c|}
\hline {\bf $k$} &    1& 2& 3& 4& 5& 6& 7& 8& 9& 10& 11& 12& 13& 14 & 15 \\
\hline \hline {\bf $|y_k|$} &  1& 2& 3& 4& 7& 7& 6& 2& 2& 4& 4& 2&
4& 4& 6 \\

\hline \hline

\hline {\bf $\mu_1$} & 1& 1& 1& 1& 1& 1& 1& 1& 1& 1& 1& 1& 1& 1&
1 \\

\hline {\bf $\mu_2$} & 3& -1& 0& 1& A& A'& 0& -1& 3& 1& -1&
      3& -1& -1& 0 \\

\hline {\bf $\mu_3$} & 3& -1& 0& 1& A'& A& 0& -1& 3& 1& -1&
      3& -1& -1& 0 \\

\hline {\bf $\mu_4$} &  6& 2& 0& 0& -1& -1& 0& 2& 6& 0& 2& 6& 2&
2& 0 \\

\hline {\bf $\mu_5$} & 7& -1& 1& -1& 0& 0& -1& 3& -1& 1& -1& 7&
-1& -1& 1\\

\hline {\bf $\mu_6$} & 7& 3& 1& 1& 0& 0& 1& 3& 7& 1& 3& -1& -1&
-1& 1 \\

\hline {\bf $\mu_7$} & 7& -1& 1& -1& 0& 0& 1& -1& 7& -1& -1& 7&
-1& -1& 1\\

\hline {\bf $\mu_8$} &7& 3& 1& 1& 0& 0& -1& -1& -1& -1& -1& 7& 3&
3& 1\\

\hline {\bf $\mu_9$} & 7& -1& 1& -1& 0& 0& 1& -1& 7& -1& -1& -1&
-1& 3& 1\\

\hline {\bf $\mu_{10}$} & 8& 0& -1& 0& 1& 1& -1& 0& 8& 0& 0& 8& 0&
0& -1\\

\hline {\bf $\mu_{11}$} & 8& 0& 2& 0& 1& 1& 0& 4& 0& 2& 0& 0& 0&
0& -2\\

\hline {\bf $\mu_{12}$} & 14& 2& -1& 0& 0& 0& 1& 2& -2& 0& -2& 14&
2& 2& -1\\

\hline {\bf $\mu_{13}$}& 14& 2& -1& 0& 0& 0& -1& 2& 14& 0& 2& -2&
-2& 2& -1\\

\hline {\bf $\mu_{14}$} & 21& 1& 0& -1& 0& 0& 0& 5& -3& 1& -3& -3&
1& -3& 0\\

\hline {\bf $\mu_{15}$} & 21& -3& 0& 1& 0& 0& 0& 1& -3& -1& 1& 21&
-3& -3& 0\\

\hline {\bf $\mu_{16}$} & 21& 1& 0& -1& 0& 0& 0& 1& 21& -1& 1& -3&
1& -3& 0\\

\hline {\bf $\mu_{17}$} & 21& 1& 0& -1& 0& 0& 0& -3& -3& 1& 1& -3&
-3& 5& 0\\

\hline {\bf $\mu_{18}$} & 21& -3& 0& 1& 0& 0& 0& 1& -3& -1& 1& -3&
1& 1& 0\\

\hline {\bf $\mu_{19}$} & 21& 5& 0& 1& 0& 0& 0& 1& -3& -1& -3& -3&
-3& 1& 0\\

\hline {\bf $\mu_{20}$} & 21& 1& 0& -1& 0& 0& 0& -3& -3& 1& 1& 21&
1& 1& 0\\

\hline {\bf $\mu_{21}$} & 21& -3& 0& 1& 0& 0& 0& -3& 21& 1& -3&
-3& 1& 1& 0\\

\hline {\bf $\mu_{22}$} & 24& 0& 0& 0& A'& A& 0& -4& 0& 2& 0& 0&
      0& 0& 0\\

\hline {\bf $\mu_{23}$} & 24& 0& 0& 0& A& A'& 0& -4& 0& 2& 0& 0&
      0& 0& 0\\

\hline {\bf $\mu_{24}$} & 28& -4& 1& 0& 0& 0& -1& 4& -4& 0& 0& -4&
0& 4& 1\\

\hline {\bf $\mu_{25}$} & 28& 4& 1& 0& 0& 0& -1& -4& -4& 0& 0& -4&
0& -4& 1\\

\hline {\bf $\mu_{26}$} & 42& -2& 0& 0& 0& 0& 0& -2& -6& 0& 2& -6&
2& -2& 0\\

\hline {\bf $\mu_{27}$} & 48& 0& 0& 0& -1& -1& 0& 8& 0& 0& 0& 0&
0& 0& 0\\

\hline {\bf $\mu_{28}$} & 56& 0& -1& 0& 0& 0& 1& 0& -8& 0& 0& -8&
0& 0& -1\\

\hline {\bf $\mu_{29}$} & 56& 0& 2& 0& 0& 0& 0& -4& 0& -2& 0& 0&
0& 0& -2\\

\hline {\bf $\mu_{30}$} & 64& 0& -2& 0& 1& 1& 0& 0& 0& 0& 0& 0& 0&
0& 2\\

\hline
\end{tabular}}
\end{center}
\caption{Character table of $M_{24}^{s_2}$
(i).}\label{tablacarM24^s2i}
\end{table}

\begin{table}[t]
\begin{center}
\tiny{\begin{tabular}{|c|c|c|c|c|c|c|c|c|c|c|c|c|c|c|c|} \hline
{\bf $k$} &   16 & 17& 18& 19& 20& 21& 22& 23 & 24 & 25 & 26& 27&
28& 29 & 30\\ \hline \hline {\bf $|y_k|$} & 8& 14& 6& 14&
4& 2& 2& 4& 12& 4& 4& 4& 2& 4& 2 \\

\hline \hline

\hline {\bf $\mu_1$} & 1 & 1& 1& 1& 1& 1& 1& 1& 1& 1& 1& 1& 1& 1& 1\\

\hline {\bf $\mu_2$} & 1& A& 0& A'& -1& 3& -1& 1& 0& 1& -1& -1& -1& 3& 3\\

\hline {\bf $\mu_3$} &  1& A'& 0& A& -1& 3&
      -1& 1& 0& 1& -1& -1& -1& 3& 3\\

\hline {\bf $\mu_4$} & 0& -1& 0& -1& 2& 6& 2& 0&  0& 0& 2& 2& 2& 6& 6\\

\hline {\bf $\mu_5$} &  -1& 0& 1& 0& 3& -1&
      3& 1& -1& 1& -1& -1& 3& -1& 7\\

\hline {\bf $\mu_6$} &  -1& 0& -1& 0& -1& -1& 3&
      -1& -1& 1& -1& -1& -1& -1& 7\\

\hline {\bf $\mu_7$} &  -1& 0& 1& 0& -1& 7&
      -1& -1& 1& -1& -1& -1& -1& 7& 7\\

\hline {\bf $\mu_8$} &  1& 0& 1& 0& -1& -1& -1&
      -1& -1& -1& -1& -1& -1& -1& 7 \\

\hline {\bf $\mu_9$} &  1& 0& -1& 0& -1& -1&
      -1& 1& -1& -1& -1& 3& 3& -1& 7\\

\hline {\bf $\mu_{10}$} &  0& 1& -1& 1& 0& 8& 0& 0&
      -1& 0& 0& 0& 0& 8& 8\\

\hline {\bf $\mu_{11}$} &   0& -1& 0& -1& 0& 0& -4& 0&
      0& -2& 0& 0& 0& 0& -8 \\

\hline {\bf $\mu_{12}$} &  0& 0& -1& 0& 2& -2& 2&
      0& 1& 0& -2& -2& 2& -2& 14\\

\hline {\bf $\mu_{13}$}&  0& 0& 1& 0& -2& -2&
      2& 0& 1& 0& -2& 2& 2& -2& 14 \\

\hline {\bf $\mu_{14}$} &  1& 0& 0& 0& -3& 5& 5&
      -1& 0& 1& 1& 1& 1& -3& 21\\

\hline {\bf $\mu_{15}$} &  1& 0& 0& 0& 1& -3& 1&
      -1& 0& -1& 1& 1& 1& -3& 21\\

\hline {\bf $\mu_{16}$} &  1& 0& 0& 0& 1& -3& 1&
      1& 0& -1& 1& -3& -3& -3& 21\\

\hline {\bf $\mu_{17}$} &  1& 0& 0& 0& 1& 5& -3&
      -1& 0& 1& 1& -3& 1& -3& 21\\

\hline {\bf $\mu_{18}$} &  -1& 0& 0& 0& -3& 5& 1&
      1& 0& -1& 1& -3& 5& -3& 21\\

\hline {\bf $\mu_{19}$} &  -1& 0& 0& 0& 1& 5& 1&
      1& 0& -1& 1& 1& -3& -3& 21\\

\hline {\bf $\mu_{20}$} &  -1& 0& 0& 0& -3& -3&
      -3& 1& 0& 1& 1& 1& -3& -3& 21\\

\hline {\bf $\mu_{21}$} &  -1& 0& 0& 0& 1& -3&
      -3& -1& 0& 1& 1& 1& 1& -3& 21\\

\hline {\bf $\mu_{22}$} & 0& -A'& 0& -A& 0& 0& 4& 0&
      0& -2& 0& 0& 0& 0& -24\\

\hline {\bf $\mu_{23}$} & 0& -A& 0& -A'& 0& 0& 4& 0&
      0& -2& 0& 0& 0& 0& -24 \\

\hline {\bf $\mu_{24}$} &  0& 0& -1& 0& 0& -4& 4&
      0& 1& 0& 0& 0& -4& 4& 28\\

\hline {\bf $\mu_{25}$} &  0& 0& -1& 0& 0& -4&
      -4& 0& 1& 0& 0& 0& 4& 4& 28\\

\hline {\bf $\mu_{26}$} &  0& 0& 0& 0& 2& 10& -2&
      0& 0& 0& -2& 2& -2& -6& 42\\

\hline {\bf $\mu_{27}$} &  0& 1& 0& 1& 0& 0& -8& 0&
      0& 0& 0& 0& 0& 0& -48 \\

\hline {\bf $\mu_{28}$} &  0& 0& 1& 0& 0& -8& 0&
      0& -1& 0& 0& 0& 0& 8& 56\\

\hline {\bf $\mu_{29}$} &  0& 0& 0& 0& 0& 0& 4& 0&
      0& 2& 0& 0& 0& 0& -56\\

\hline {\bf $\mu_{30}$} &  0& -1& 0& -1& 0& 0& 0& 0&
      0& 0& 0& 0& 0& 0& -64\\

\hline
\end{tabular}}
\end{center}
\caption{Character table of $M_{24}^{s_2}$
(ii).}\label{tablacarM24^s2ii}
\end{table}

For every $k$, $1 \leq k \leq 30$, we call $\rho_k=(\rho_k,V_k)$
the irreducible representation of $M_{24}^{s_2}$ whose character
is $\mu_k$. We will prove that the Nichols algebra
$\toba(\oc_{s_{2}}, \rho_k)$ is infinite-dimensional, for every
$k$, $1 \leq k \leq 30$. We compute that $s_2\in
\oc_{30}^{M_{24}^{s_2}}$. If $k\neq 11$, $22$, $23$, $27$, $29$,
$30$, then the result follows from Lemma \ref{odd}, because
$q_{s_2s_2}\neq -1$. Assume that $k= 11$, $22$, $23$, $27$, $29$
or $30$; thus, $q_{s_2s_2}= -1$. From Tables \ref{tablacarM24^s2i}
and \ref{tablacarM24^s2ii}, we have that $8 \leq\deg(\rho_k)$
even. We compute that $\oc_{s_{2}}\cap M_{24}^{s_{2}}$ has 281
elements and it contains $\sigma_1:=s_{2}$ and
$$\sigma_2:=( 2,16)( 3, 4)( 5,19)( 6,22)( 7,24)(
8,10)(11,18)(14,23),$$ with $\sigma_2 \in \oc_{9}^{M_{24}^{s_2}}$.
We choose $g_1:=\id$ and
\begin{align*} g_2&:=( 1, 2)( 3,13,22,17,10,15)( 4,12, 6, 9,
8,21)( 5,11, 7)(16,20)(18,24,19).
\end{align*}
These elements are in $M_{24}$ and they satisfy $\sigma_1 g_1=g_1
\sigma_1$, $\sigma_2 g_1=g_1 \sigma_2$, $\sigma_1 g_2=g_2
\sigma_2$ and $\sigma_2 g_2=g_2 \sigma_1$. Since $\sigma_1$ and
$\sigma_2$ commute there exists a basis $\{v_l \, | \, 1\leq l
\leq \deg(\rho_k)\}$ of $V_k$, the vector space affording
$\rho_k$, composed by simultaneous eigenvectors of
$\rho_k(\sigma_1)=-\Id$ and $\rho_k(\sigma_2)$. Let us call
$\rho_k(\sigma_2)v_l=\lambda_l v_l$, $1\leq l \leq \deg(\rho_k)$,
where $\lambda_l=\pm 1$, due to $|\sigma_2|=2$. From Table
\ref{tablacarM24^s2i}, we have that
$\sum_{l=1}^{\deg(\rho_k)}\lambda_l=0$. Reordering the basis we
can suppose that
$\lambda_1=\cdots=\lambda_{\deg(\rho_k)/2}=1=-\lambda_{1+\deg(\rho_k)/2}=\cdots=-\lambda_{\deg(\rho_k)}$.
It is straightforward to check that if $W:=\ku$ - span of
$\{g_1v_l, g_2v_l \, | \, 1\leq l \leq \deg(\rho_k)\}$, then $W$
is a braided vector subspace of $M(\oc_{s_{2}},\rho)$ of Cartan
type whose associated Cartan matrix $\mathcal A$ has at least two
row with three $-1$ or more. This means that the corresponding
Dynkin diagram has at least two vertices with three edges or more;
thus, $\mathcal A$ is not of finite type. Hence, $\dim
\toba(\oc_{s_2}, \rho_k) = \infty$.


\smallbreak

\emph{CASE: $j=3$}. The representative $s_{3}$ is
{\small{\begin{align*} ( 1,11)( 2,23)( 3, 7)( 4,13)( 5, 6)( 8,18)(
9,16)(10,15)(12,24)(14,17)(19,21) (20,22).
\end{align*}}}

\vspace*{-0.4cm}\noindent We compute that the centralizer
$M_{24}^{s_{3}}$ is a non-abelian group of order 7680 whose
character table is given by Tables \ref{tablacarM24^s3i} and
\ref{tablacarM24^s3ii}, where $A=- i \sqrt{5}$.

{\tiny{
\begin{table}[t]
\begin{center}
\begin{tabular}{|c|c|c|c|c|c|c|c|c|c|c|c|c|c|c|c|}
\hline {\bf $k$} &    1& 2& 3& 4& 5& 6& 7& 8& 9& 10& 11& 12& 13&
14 & 15 \\
\hline \hline {\bf $|y_k|$} &  1& 5& 10& 10& 10& 3& 6& 12& 6& 6&
2& 4& 8& 8& 4\\

\hline \hline

\hline {\bf $\mu_1$} & 1& 1& 1& 1& 1& 1& 1& 1& 1& 1& 1& 1& 1& 1&
1\\

\hline {\bf $\mu_2$} &  1& 1& 1& 1& 1& 1& 1& -1& -1& 1& 1& 1& -1&
-1& 1\\

\hline {\bf $\mu_3$} & 4& -1& -1& -1& -1& 1& 1& 1& 1& 1& 4& 0& 0&
0& 0\\

\hline {\bf $\mu_4$} & 4& -1& -1& -1& -1& 1& 1& -1& -1& 1& 4& 0&
0& 0& 0\\

\hline {\bf $\mu_5$} & 5& 0& 0& 0& 0& -1& -1& 1& 1& -1& 5& 1& -1&
-1& 1\\

\hline {\bf $\mu_6$} & 5& 0& 0& 0& 0& -1& -1& -1& -1& -1& 5& 1& 1&
1& 1\\

\hline {\bf $\mu_7$} & 6& 1& 1& 1& 1& 0& 0& 0& 0& 0& 6& -2& 0& 0&
-2\\

\hline {\bf $\mu_8$} & 6& 1& 1& -1& -1& 0& 0& 0& 0& 0& -2& 2& 0&
0& -2\\

\hline {\bf $\mu_9$} & 6& 1& 1& -1& -1& 0& 0& 0& 0& 0& -2& 2& 0&
0& -2\\

\hline {\bf $\mu_{10}$} & 6& 1& 1& -1& -1& 0& 0& 0& 0& 0& -2& -2&
2i& -2i& 2\\

\hline {\bf $\mu_{11}$} & 6& 1& 1& -1& -1& 0& 0& 0& 0& 0& -2& -2&
-2i& 2i& 2\\

\hline {\bf $\mu_{12}$} & 10& 0& 0& 0& 0& 1& 1& -1& 1& -1& 2& 2&
0& 0& -2\\

\hline {\bf $\mu_{13}$}& 10& 0& 0& 0& 0& 1& 1& 1& -1& -1& 2& 2& 0&
0& -2\\

\hline {\bf $\mu_{14}$} & 10& 0& 0& 0& 0& 1& 1& -1& 1& -1& 2& -2&
0& 0& 2\\

\hline {\bf $\mu_{15}$} & 10& 0& 0& 0& 0& 1& 1& 1& -1& -1& 2& -2&
0& 0& 2\\

\hline {\bf $\mu_{16}$} & 12& 2& -2& 0& 0& 0& 0& 0& 0& 0& 4& 0& 0&
0& 0\\

\hline {\bf $\mu_{17}$} & 12& 2& -2& 0& 0& 0& 0& 0& 0& 0& 4& 0& 0&
0& 0\\

\hline {\bf $\mu_{18}$} & 15& 0& 0& 0& 0& 0& 0& 0& 0& 0& -1& -1&
-1& -1& -1\\

\hline {\bf $\mu_{19}$} & 15& 0& 0& 0& 0& 0& 0& 0& 0& 0& -1& 3& 1&
1& 3\\

\hline {\bf $\mu_{20}$} & 15& 0& 0& 0& 0& 0& 0& 0& 0& 0& -1& 3&
-1& -1& 3\\

\hline {\bf $\mu_{21}$} & 15& 0& 0& 0& 0& 0& 0& 0& 0& 0& -1& -1&
1& 1& -1\\

\hline {\bf $\mu_{22}$} & 20& 0& 0& 0& 0& -1& -1& -1& 1& 1& 4& 0&
0& 0& 0\\

\hline {\bf $\mu_{23}$} & 20& 0& 0& 0& 0& -1& -1& 1& -1& 1& 4& 0&
0& 0& 0\\

\hline {\bf $\mu_{24}$} & 20& 0& 0& 0& 0& 2& -2& 0& 0& 0& -4& 0&
0& 0& 0\\

\hline {\bf $\mu_{25}$} & 20& 0& 0& 0& 0& 2& -2& 0& 0& 0& -4& 0&
0& 0& 0\\

\hline {\bf $\mu_{26}$} & 24& -1& -1& 1& 1& 0& 0& 0& 0& 0& -8& 0&
0& 0& 0\\

\hline {\bf $\mu_{27}$} & 24& -1& 1& A& -A& 0& 0&
      0& 0& 0& 8& 0& 0& 0& 0\\

\hline {\bf $\mu_{28}$} & 24& -1& 1& -A& A& 0&  0& 0& 0& 0& 8& 0&
0& 0& 0\\

\hline {\bf $\mu_{29}$} & 30& 0& 0& 0& 0& 0& 0& 0& 0& 0& -2& -2&
0& 0& -2\\

\hline {\bf $\mu_{30}$} & 40& 0& 0& 0& 0& -2& 2& 0& 0& 0& -8& 0&
0& 0& 0\\

\hline
\end{tabular}
\end{center}
\caption{Character table of $M_{24}^{s_3}$
(i).}\label{tablacarM24^s3i}
\end{table}}}

{\tiny{
\begin{table}[t]
\begin{center}
\begin{tabular}{|c|c|c|c|c|c|c|c|c|c|c|c|c|c|c|c|}
\hline {\bf $k$} & 16 & 17& 18& 19& 20& 21& 22& 23 & 24 & 25 & 26&
27& 28& 29 & 30\\
\hline \hline {\bf $|y_k|$} &   4& 2& 4& 2& 4&
2& 4& 4& 2& 2& 4& 2& 4& 2& 2 \\

\hline \hline

\hline {\bf $\mu_1$} & 1& 1& 1& 1& 1& 1& 1&
1& 1& 1& 1& 1& 1& 1& 1\\

\hline {\bf $\mu_2$} &  -1& 1& -1& 1& -1& 1& 1&
      1& 1& 1& -1& 1& -1& 1& -1\\

\hline {\bf $\mu_3$} &  -2& 4& -2& 4& -2& 4& 0&
      0& 4& 0& 0& 0& 0& 0& -2\\

\hline {\bf $\mu_4$} &  2& 4& 2& 4& 2& 4& 0&
      0& 4& 0& 0& 0& 0& 0& 2\\

\hline {\bf $\mu_5$} &  1& 5& 1& 5& 1& 5& 1& 1&
      5& 1& -1& 1& -1& 1& 1\\

\hline {\bf $\mu_6$} &  -1& 5& -1& 5& -1& 5& 1&
      1& 5& 1& 1& 1& 1& 1& -1\\

\hline {\bf $\mu_7$} &  0& 6& 0& 6& 0& 6& -2& -2&
      6& -2& 0& -2& 0& -2& 0\\

\hline {\bf $\mu_8$} &  0& 6& 0& -2& 0& 2& 0& 0&
      -6& 2& 2& 2& -2& -2& 0\\

\hline {\bf $\mu_9$} &  0& 6& 0& -2& 0& 2& 0& 0&
      -6& 2& -2& 2& 2& -2& 0 \\

\hline {\bf $\mu_{10}$} &  0& 6& 0& -2&
      0& 2& 0& 0& -6& -2& 0& -2& 0& 2& 0\\

\hline {\bf $\mu_{11}$} &  0& 6& 0& -2&
      0& 2& 0& 0& -6& -2& 0& -2& 0& 2& 0 \\

\hline {\bf $\mu_{12}$} &  2& 10& 2& 2& -2& -2& 0&
      0& -10& -2& 0& -2& 0& 2& -2\\

\hline {\bf $\mu_{13}$}&  -2& 10& -2& 2& 2& -2& 0&
      0& -10& -2& 0& -2& 0& 2& 2\\

\hline {\bf $\mu_{14}$} &  0& 10& -4& 2& 0& -2& 0&
      0& -10& 2& 0& 2& 0& -2& 4 \\

\hline {\bf $\mu_{15}$} &  0& 10& 4& 2& 0& -2& 0&
      0& -10& 2& 0& 2& 0& -2& -4\\

\hline {\bf $\mu_{16}$} &  0& -12& 0& -4& 0& 0& 2&
      -2& 0& 4& 0& -4& 0& 0& 0\\

\hline {\bf $\mu_{17}$} &  0& -12& 0& -4& 0& 0& -2&
      2& 0& -4& 0& 4& 0& 0& 0\\

\hline {\bf $\mu_{18}$}  &  -1& 15& 3& -1& -1& -1&
      -1& -1& 15& 3& 1& 3& 1& 3& 3\\

\hline {\bf $\mu_{19}$} &  -1& 15& 3& -1& -1& -1& -1&
      -1& 15& -1& -1& -1& -1& -1& 3\\

\hline {\bf $\mu_{20}$} &  1& 15& -3& -1& 1& -1&
      -1& -1& 15& -1& 1& -1& 1& -1& -3\\

\hline {\bf $\mu_{21}$} &  1& 15& -3& -1& 1& -1&
      -1& -1& 15& 3& -1& 3& -1& 3& -3\\

\hline {\bf $\mu_{22}$} &  -2& 20& 2& 4& 2& -4& 0&
      0& -20& 0& 0& 0& 0& 0& -2\\

\hline {\bf $\mu_{23}$} &  2& 20& -2& 4& -2& -4& 0&
      0& -20& 0& 0& 0& 0& 0& 2 \\

\hline {\bf $\mu_{24}$} &  0& -20& 0& 4& 0& 0& -2&
      2& 0& 4& 0& -4& 0& 0& 0\\

\hline {\bf $\mu_{25}$} &  0& -20& 0& 4& 0& 0& 2&
      -2& 0& -4& 0& 4& 0& 0& 0\\

\hline {\bf $\mu_{26}$} &  0& 24& 0& -8& 0& 8& 0&
      0& -24& 0& 0& 0& 0& 0& 0\\

\hline {\bf $\mu_{27}$} & 0& -24& 0& -8& 0& 0& 0& 0& 0& 0& 0& 0& 0& 0& 0\\

\hline {\bf $\mu_{28}$} &  0& -24& 0& -8& 0& 0& 0& 0& 0& 0& 0& 0& 0& 0&  0\\

\hline {\bf $\mu_{29}$} &  0& 30& 0& -2& 0& -2& 2& 2& 30& -2& 0& -2& 0& -2& 0\\

\hline {\bf $\mu_{30}$} &  0& -40& 0& 8& 0& 0& 0& 0&
      0& 0& 0& 0& 0& 0& 0\\

\hline
\end{tabular}
\end{center}
\caption{Character table of $M_{24}^{s_3}$
(ii).}\label{tablacarM24^s3ii}
\end{table}}}

For every $k$, $1 \leq k \leq 30$, we call $\rho_k=(\rho_k,V_k)$
the irreducible representation of $M_{24}^{s_3}$ whose character
is $\mu_k$. We will prove that the Nichols algebra
$\toba(\oc_{s_{3}}, \rho_k)$ is infinite-dimensional, for every
$k$, $1 \leq k \leq 30$. We compute that $s_3\in
\oc_{17}^{M_{24}^{s_3}}$. Now, if $k\neq 16$, $17$, $24$, $25$,
$27$, $28$, $30$, then the result follows from Lemma \ref{odd},
because $q_{s_3s_3}\neq -1$. Assume that $k= 16$, $17$, $24$,
$25$, $27$, $28$ or $30$; thus, $q_{s_3s_3}= -1$. We compute that
$\oc_{s_{3}}\cap M_{24}^{s_{3}}$ has 278 elements and it contains
$\sigma_1:=s_{3}$ and
\begin{align*} \sigma_2:=(& 1, 5)( 2,10)(
3,12)( 4, 8)( 6,11)( 7,24)( 9,19)(13,18)(14,22)(15,23)\\
 (&16,21) (17,20),
\end{align*}
with $\sigma_2 \in \oc_{30}^{M_{24}^{s_3}}$. We choose $g_1:=\id$
and
\begin{align*} g_2&:=( 5,11)( 7,12)( 8,13)(
9,20)(10,23)(14,21)(16,17)(19,22).
\end{align*}
These elements are in $M_{24}$ and they satisfy $\sigma_1 g_1=g_1
\sigma_1$, $\sigma_2 g_1=g_1 \sigma_2$, $\sigma_1 g_2=g_2
\sigma_2$ and $\sigma_2 g_2=g_2 \sigma_1$. Since $12\leq
\deg(\rho_k)$ even and $\mu_k(\oc_{30}^{M_{24}^{s_3}})=0$, for $k=
16$, $17$, $24$, $25$, $27$, $28$ or $30$, we can proceed as in
the previous case. Therefore,
$\dim\toba(\oc_{s_3},\rho_k)=\infty$, for all $k$.\epf

\begin{obs}
We compute that $M_{24}^{s_4}$, $M_{24}^{s_5}$, $M_{24}^{s_7}$,
$M_{24}^{s_9}$, $M_{24}^{s_{10}}$, $M_{24}^{s_{11}}$,
$M_{24}^{s_{12}}$ and $M_{24}^{s_{13}}$ have 17, 18, 26, 20, 15,
15, 21 and 21 conjugacy classes, respectively. Hence, there are
502 possible pairs $(\oc,\rho)$ for $M_{24}$; 492 of them have
infinite-dimensional Nichols algebras, and 10 have negative
braiding.
\end{obs}


\section{Using techniques based on non-abelian
subracks}\label{sec:non-abelian}

In the previous section, we discard the pairs $(\oc,\rho)$ with
$\dim \toba(\oc,\rho)=\infty$ by mean of abelian subracks of
$\oc$. As a result, we showed that remain 23 pairs which give rise
to a braided vector spaces with negative braiding. In this
section, we consider these 23 remained pairs $(\oc,\rho)$ and show
that 16 of them lead to infinite-dimensional Nichols algebras
using non-abelian subracks of $\oc$ -- see Subsection
\ref{subsec:nonabeltech}.

\subsection{The group $M_{11}$} We have 5 remained cases.

\smallbreak

\emph{CASE: $j=10$}. We choose in $\oc_{s_{10}}$ the following
elements: $\sigma_0:=s_{10}$, 
$\sigma_1:=( 1, 6, 8)( 2, 5, 3, 4,10, 9)( 7,11)$ and
$\sigma_2:=\sigma_0 \trid\sigma_1$. It is easy to see that the
family $(\sigma_i)_{i\in \Z_3}$ is of type $\D_3$ in
$\oc_{s_{10}}$. Then $\dim\toba(\oc_{s_{10}},\chi_{(-1)})=\infty$,
by Lemma \ref{coro:dp-cor}, with $p=3$.

\subsection{The group $M_{12}$} We have 4 remained cases.

\smallbreak

\emph{CASE: $j=14$}. We choose in $\oc_{s_{14}}$ the following
elements: $\sigma_1:=s_{14}$, 
$\sigma_2:=( 1, 3, 5,11, 4, 8,10,12)( 7, 9)$,
$\sigma_3:=\sigma_2\trid \sigma_1$, $\sigma_4:=\sigma_3\trid
\sigma_1$, $\sigma_5:=\sigma_4\trid \sigma_1$,
$\sigma_6:=s_{14}^3$ and $\tau_l:=\sigma_l^5$, $1\leq l \leq 6$.
By straightforward computation we can check that the family
$(\sigma_l)_{l=1}^6 \cup (\tau_l)_{l=1}^6$ is of type
$\oct^{(2)}$. 
Then, $\dim\toba(\oc_{s_{14}},\chi_{(-1)})=\infty$ from Lemma
\ref{co:especial2}.

\smallbreak

\emph{CASE: $j=5$}. We choose in $\oc_{s_{5}}$ the following
elements: $\sigma_1:=s_{5}$, 
$\sigma_2:=( 1, 4,12, 7)( 2, 6, 3, 8,11,10, 9, 5)$,
$\sigma_3:=\sigma_2\trid \sigma_1$, $\sigma_4:=\sigma_3\trid
\sigma_1$, $\sigma_5:=\sigma_4\trid \sigma_1$, $\sigma_6:=s_{5}^3$
and $\tau_l:=\sigma_l^5$, $1\leq l \leq 6$. We can check that the
family $(\sigma_l)_{l=1}^6 \cup (\tau_l)_{l=1}^6$ is of type
$\oct^{(2)}$. 
Then, $\dim\toba(\oc_{s_{5}},\chi_{(-1)})=\infty$ from Lemma
\ref{co:especial2}.

\smallbreak

\emph{CASE: $j=2$}. We choose in $\oc_{s_{2}}$ the following
elements: $\sigma_0:=s_2$, 
$\sigma_1:=( 1, 2,12,11, 8,10)( 3, 6, 9)( 4, 5)$ and
$\sigma_2:=\sigma_0 \trid\sigma_1$. It is easy to see that this
family $(\sigma_i)_{i\in \Z_3}$ is of type $\D_3$ in $\oc_{s_2}$.
Then $\dim\toba(\oc_{s_{2}},\chi_{(-1)})=\infty$ from Lemma
\ref{coro:dp-cor}, with $p=3$.

\subsection{The group $M_{22}$} We have one remained case.
We choose in $\oc_{s_{4}}$ the following elements:
$\sigma_1:=s_{4}$, $$\sigma_2:= ( 1, 5,13,10,11, 7,12,22)( 2, 9)(
3,21,19,15,17,18, 6, 8)( 4,14,20,16),$$ $\sigma_3:=\sigma_2\trid
\sigma_1$, $\sigma_4:=\sigma_3\trid \sigma_1$,
$\sigma_5:=\sigma_4\trid \sigma_1$, $\sigma_6:=s_{4}^{-1}$ and
$\tau_l:=\sigma_l^5$, $1\leq l \leq 6$. By straightforward
computation we can check that the family $(\sigma_l)_{l=1}^6 \cup
(\tau_l)_{l=1}^6$ is of type
$\oct^{(2)}$. 
Then, $\dim\toba(\oc_{s_{4}},\chi_{(-1)})=\infty$ from Lemma
\ref{co:especial2}.

In view of Theorem \ref{teorM22} and the previous paragraph we can
state the following result.

\begin{theorem}\label{mainteor}
Any finite-dimensional complex pointed Hopf algebra $H$ with
$G(H)\simeq M_{22}$ is necessarily isomorphic to the group algebra
of $M_{22}$.\qed
\end{theorem}

\subsection{The group $M_{23}$} We have 3 remained cases.

\smallbreak

\emph{CASE: $j=9$}. We choose in $\oc_{s_{9}}$ the following
elements: $\sigma_1:=s_{9}$, $$\sigma_2:= ( 1, 3,
5,20,10,14,13,23)( 2,15, 7, 8)( 4,22,12, 6,17,16,21,11)( 9,19),$$
$\sigma_3:=\sigma_2\trid \sigma_1$, $\sigma_4:=\sigma_3\trid
\sigma_1$, $\sigma_5:=\sigma_4\trid \sigma_1$,
$\sigma_6:=s_{9}^{-1}$ and $\tau_l:=\sigma_l^5$, $1\leq l \leq 6$.
We compute that the family $(\sigma_l)_{l=1}^6 \cup
(\tau_l)_{l=1}^6$ is of type
$\oct^{(2)}$. 
Then, $\dim\toba(\oc_{s_{9}},\chi_{(-1)})=\infty$ from Lemma
\ref{co:especial2}.

\subsection{The group $M_{24}$} We have 10 remained cases.

\smallbreak

\emph{CASE: $j=6$}. The representative $s_{6}$ is
$$( 1, 9,20,17)( 2, 6)( 3,10)( 4, 8)(5,24,19,7)(11,14,18,23)(12,21,13,15)(16,\!22).$$
We choose in $\oc_{s_{6}}$ the following elements:
$\sigma_1:=s_{6}$, $\sigma_2$ to be
$$( 1, 9,20,17)( 2,11)(3,14)( 4,18)( 5,19, 7,24)( 6,10, 8,22)(12,21,15,13)(16,\!23),$$
$\sigma_3:=\sigma_2\trid \sigma_1$, $\sigma_4:=\sigma_3\trid
\sigma_1$, $\sigma_5:=\sigma_4\trid \sigma_1$,
$\sigma_6:=\sigma_2\trid \sigma_3$, $\tau_1:=\sigma_6^{-1}$,
$\tau_2:=\sigma_4^{-1}$, $\tau_3:=\sigma_5^{-1}$,
$\tau_4:=\sigma_2^{-1}$, $\tau_5:=\sigma_3^{-1}$ and
$\tau_6:=\sigma_1^{-1}$. By straightforward computation we can
check that the family $(\sigma_l)_{l=1}^6 \cup (\tau_l)_{l=1}^6$
is of type
$\oct^{(2)}$. 
Now, assume that $\rho=\rho_{2,6}$ or $\rho_{3,6}$; then
$q_{\sigma_1\sigma_1}=-1$ because $s_{6}\in
\oc_{23}^{M_{24}^{s_6}}$. We define $g:=( 3,16)( 5,12)( 6, 8)(
7,15)( 9,17)(13,19)(14,23)(21,24)$; we verify that $g\in M_{24}$
and that $g\trid \sigma_1=\tau_1$. Also we compute that $\tau_1$,
$\sigma_6\in \oc_{24}^{M_{24}^{s_6}}$. Hence, from Table
\ref{tablacarM24^s6ii}, we have that
$\rho(\sigma_6)=\rho(\tau_1)=\rho(g\sigma_6 g)=-1$. Therefore,
$\dim\toba(\oc_{s_{6}},\rho)=\infty$, from Lemma
\ref{teor:aplicAHSch2}.

\smallbreak

\emph{CASE: $j=8$}. The representative $s_8$ is $$( 1, 4,24,14)(
2,21,15, 6)( 3,16, 8,12)( 5,11,23,20)( 7,18,17,13) (
9,10,22,19).$$ We choose in $\oc_{s_{8}}$ the following elements:
$\sigma_0:=s_{8}$, $\sigma_1$ to be
$$( 1, 2,24,15)( 3, 5, 8,23)( 4,19,14,10)( 6,22,21, 9)( 7,16,17,12)(11,13,20,18)$$
and $\sigma_2:=\sigma_0 \trid\sigma_1$. We compute that
$(\sigma_i)_{i\in \Z_3}$ is a family of type $\D_3$ in
$\oc_{s_{8}}$. Now, if $\rho=\rho_{2,8}$ or $\rho_{2,8}$, then
$q_{\sigma_0\sigma_0}=-1$, and
$\dim\toba(\oc_{s_{14}},\rho)=\infty$ by from Lemma
\ref{coro:dp-cor}, with $p=3$.

\smallbreak

\emph{CASE: $j=14$}. 
The representative $s_{14}$ is
$$( 2,10, 4,16,20,15, 6,18)( 3,17,11, 5, 7,12,21,13)(
9,19)(14,24,23,22).$$ We choose in $\oc_{s_{14}}$ the following
elements: $\sigma_1:=s_{14}$, $$\sigma_2:=( 2, 4,18,15,20,
6,16,10)( 3,12,13,21, 7,17, 5,11)( 8, 9)(14,22,24,23),$$
$\sigma_3:=\sigma_2\trid \sigma_1$, $\sigma_4:=\sigma_3\trid
\sigma_1$, $\sigma_5:=\sigma_4\trid \sigma_1$,
$\sigma_6:=s_{14}^{3}$ and $\tau_l:=\sigma_l^5$, $1\leq l \leq 6$.
By straightforward computation we can check that the family
$(\sigma_l)_{l=1}^6 \cup (\tau_l)_{l=1}^6$ is of type
$\oct^{(2)}$. 
Now, if $\rho=\epsilon \otimes \chi_{(-1)}$ or $\sgn \otimes
\chi_{(-1)}$, then $\dim\toba(\oc_{s_{14}},\rho)=\infty$, from
Lemma \ref{co:especial2}.

\smallbreak

\emph{CASE: $j=17$}. 
The representative $s_{17}$ is
$$( 1, 9,12,10,17,14, 3,23, 5,21,19,13)( 2,18)( 4, 8,15,20)(
6,16, 7,24,22,11).$$ We choose in $\oc_{s_{17}}$ the following
elements: $\sigma_1:=s_{17}$, $\sigma_2$ to be $$( 1, 9,13,
3,17,14,10,19, 5,21,23,12)( 2,18)( 4, 8,20,15)( 6,11,
7,16,22,24),$$ $\sigma_3:=\sigma_2\trid \sigma_1$,
$\sigma_4:=\sigma_3\trid \sigma_1$, $\sigma_5:=\sigma_4\trid
\sigma_1$, $\sigma_6:=s_{17}^{7}$ and $\tau_l:=\sigma_l^5$, $1\leq
l \leq 6$. We compute that the family $(\sigma_l)_{l=1}^6 \cup
(\tau_l)_{l=1}^6$ is of type
$\oct^{(2)}$. 
Then, $\dim\toba(\oc_{s_{17}},\chi_{(-1)})=\infty$ from Lemma
\ref{co:especial2}.

\smallbreak

\emph{CASE: $j=18$}.  
The representative $s_{18}$ is
$$( 1,18,20, 4,17, 5,24,13,11,14, 7,23)(
2,10,16,21,22, 8,15,19,12, 6, 9, 3).$$ We choose in $\oc_{s_{18}}$
the following elements: $\sigma_1:=s_{18}$, $\sigma_2$ to be $$(
1, 2,10,14,17,22, 8,18,11,12, 6, 5)( 3,15,20, 4,21, 9,24,13,19,16,
7,23),$$ $\sigma_3:=\sigma_2\trid \sigma_1$,
$\sigma_4:=\sigma_3\trid \sigma_1$, $\sigma_5:=\sigma_4\trid
\sigma_1$, $\sigma_6:=s_{18}^{7}$ and $\tau_l:=\sigma_l^5$, $1\leq
l \leq 6$. We compute that the family $(\sigma_l)_{l=1}^6 \cup
(\tau_l)_{l=1}^6$ is of type
$\oct^{(2)}$. 
Then, $\dim\toba(\oc_{s_{18}},\chi_{(-1)})=\infty$ from Lemma
\ref{co:especial2}.

\smallbreak

\emph{CASE: $j=19$}. 
The representative $s_{19}$ is
$$( 1, 5,20,23,19, 2,18, 7,17, 9,21,24, 6,12)( 3,14,
8,15,13,11,16)( 4,22).$$ We choose in $\oc_{s_{19}}$ the following
elements: $\sigma_0:=s_{19}$,
$$\sigma_1:=( 1,14,20,15,19,11,18, 3,17, 8,21,13, 6,16)( 2,12, 7, 5, 9,23,24)( 4,10)$$
and $\sigma_2:=\sigma_0 \trid\sigma_1$. We compute that
$(\sigma_i)_{i\in \Z_3}$ is a family of type $\D_3$ in
$\oc_{s_{19}}$. Then $\dim\toba(\oc_{s_{19}},\chi_{(-1)})=\infty$,
from Lemma \ref{coro:dp-cor}, with $p=3$ and $k=9$..

\smallbreak

\emph{CASE: $j=20$}. 
The representative is $s_{20}=s_{19}^{-1}$. Now, the family
$(\sigma_i^{-1})_{i \in \Z_3}$, with $\sigma_i$ as in the case
$j=19$ above, is of type $\D_3$. Then
$\dim\toba(\oc_{s_{20}},\chi_{(-1)})=\infty$, from Lemma
\ref{coro:dp-cor}, with $p=3$ and $k=9$.

In view of Theorem \ref{teorM24} and the previous paragraph we can
state the following result.

\begin{theorem}\label{mainteor2}
Any finite-dimensional complex pointed Hopf algebra $H$ with
$G(H)\simeq M_{24}$ is necessarily isomorphic to the group algebra
of $M_{24}$.\qed
\end{theorem}

\begin{obs}
In the 7 remained cases, that appear in Table \ref{maintabla}, we
compute that there is no family of type $\oct^{(2)}$ nor $\D_p$,
for any odd prime integer $p$, inside the respective conjugacy
classes. Hence, in these cases we cannot decide whether the
dimension of $\toba(\oc,\rho)$ is infinite or not with the methods
available today.
\end{obs}

\subsection*{Acknowledgement}
The author is grateful to N. Andruskiewitsch for his valuable
advice and encouragement with the presentation of the paper and to
I. Heckenberger for pointed out to me a more convenient form of
Lemma \ref{Hecke} and important suggestions. Also, I want to thank
to L. Vendramin for interesting discussions and for his help about
using GAP.

\end{document}